\documentclass[a4paper,12pt,twoside]{article}
\usepackage[british]{babel}
\usepackage{amssymb,amsmath,amsfonts,verbatim,xkeyval,bm,upgreek}
\usepackage{graphicx}
\usepackage{subfigure}
\usepackage{ wasysym } 

\renewcommand{\cal}[1]{{\mathcal #1}}         

\newtheorem{theorem}{Theorem}[section]
\newtheorem{lemma}[theorem]{Lemma}
\newtheorem{proposition}[theorem]{Proposition}

\newtheorem{definition}[theorem]{Definition}
\newenvironment{proof}[1][Proof]{\begin{trivlist}
\item[\hskip \labelsep {\bfseries #1}]}{\end{trivlist}
{\nobreak\quad\nobreak\hfill\nobreak{$\blacksquare$}}}

\def\lie#1{\Lscr_{#1}}

\def\Lie#1{\Lscr_{#1}}
\def\reali{\mathbb{R}}
\def\complessi{\mathbb{C}}
\def\toro{\mathbb{T}}
\def\interi{\mathbb{Z}}
\def\naturali{\mathbb{N}}

\def\poisson#1#2{\lbrace#1,#2\rbrace}

\def\dist{\mathop{\rm dist}}

\def\det{\mathop{\rm det}}

\def\vet#1{{\bm #1}}
\def\epsilon{{\varepsilon}}
\def\theta{{\vartheta}}

\def\imunit{{i}}
\let\flusso=\phi       
\def\phi{{\varphi}}
\def\Ascr{{\cal A}}
\def\Bscr{{\cal B}}

\def\Dscr{{\cal D}}

\def\Fscr{{\cal F}}
\def\Gscr{{\cal G}}
\def\Hscr{{\cal H}}

\def\Kscr{{\cal K}}
\def\Lscr{{\cal L}}
\def\Mscr{{\cal M}}
\def\Nscr{{\cal N}}
\def\Oscr{{\cal O}}
\def\Pscr{{\cal P}}

\def\Rscr{{\cal R}}

\def\Uscr{{\cal U}}
\def\Vscr{{\cal V}}
\def\Wscr{{\cal W}}
\def\Xscr{{\cal X}}

\def\imunit{{\bf i}}
\def\Pset{{\cal P}}
\def\PPset{\widehat{\cal P}}
\def\scalprod#1#2{#1\cdot#2}
\def\rmI{{\rm I}}
\def\rmII{{\rm II}}

\def\parder#1#2{{\partial#1 \over \partial#2}}
\def\realpart{\mathop{\rm Re}\nolimits}
\def\imaginary{\mathop{\rm Im}\nolimits}
\def\Ebarra{{\bar E}}
\def\ee{{e}}

\def\diff{\mathord{\rm d}}

\def\lt{<}
\title{\bf Normal form for lower dimensional elliptic tori in Hamiltonian systems\thanks{{\it 2020
      Mathematics Subject Classification.}  Primary: 37J40; Secondary:
     70H08, 70H15. {\it Key words and phrases:} lower dimensional invariant tori, KAM theory, normal form
    methods, perturbation theory for Hamiltonian systems.}}

\author{
{\bf CHIARA CARACCIOLO}\\
{\small Dipartimento di Matematica, 
Universit\`a degli Studi di Milano,}\\
{\small via Saldini 50, 20133\ ---\ Milano (Italy).}\\
{\small e-mail:
  {\tt chiara.caracciolo@unimi.it}}\\
}

\begin{document}
\maketitle


\selectlanguage{british}
\thispagestyle{empty}

\begin{abstract}
We give a proof of the convergence of an algorithm for the
construction of lower dimensional elliptic tori in nearly integrable
Hamiltonian systems. The existence of such invariant tori is proved by
leading the Hamiltonian to a suitable normal form. In particular, we
adapt the procedure described in a previous work by Giorgilli and
co-workers, where the construction was made so as to be used in the
context of the planetary problem. We extend the proof of the
convergence to the cases in which the two sets of frequencies,
describing the motion along the torus and the transverse oscillations,
have the same order of magnitude.
\end{abstract}
\bigskip

\markboth{C. Caracciolo}{{Normal form for lower dimensional elliptic tori $\ldots$}}

\maketitle

\section*{Foreword}
In a previous work in collaboration with Locatelli
(see~\cite{Car-Loc-2021}), we have emphasized the key role played by
the invariant elliptic tori in the dynamics of the FPUT
non-linear chains with a small number of nodes.  In fact, in that
article we have described a suitable algorithm constructing the normal
forms corresponding to those lower-dimensional manifolds. Moreover, we
have shown that most of the initial conditions of the same type of
those considered in the original work written by Fermi, Pasta and Ulam
(see~\cite{FPU-1955}) stay in the effective stability region around
elliptic tori. These are keystone concepts to which Prof. Antonio
Giorgilli has repeatedly given fundamental contributions during his
scientific career: dynamics in FPUT models, perturbative
methods constructing normal forms, effective stability,
lower-dimensional tori and very slow diffusion in the neighborhood of
an invariant manifold (see,
e.g.,~\cite{Ber-Gal-Gio-2004},~\cite{Giorgilli-1995},~\cite{Giorgilli-1997},~\cite{Gio-Loc-San-2014},~\cite{Giorgilli-1997.4}).
Since this paper is devoted to a rigorous proof of the convergence of
such a constructive algorithm, it looks very natural to dedicate this
work of mine to him, so as to include it in this Volume in honor of
Antonio Giorgilli.

I think that my personal growth during my Ph.D. studies can be seen as
an ``indirect fruit" of the huge work he did in the foundation and
consolidation of a research group that is well recognized at an
international level for its contributions in the field of Hamiltonian
perturbation theory. Indeed, I could enjoy his marvelous notes, his
very explanatory articles and I can work by using the powerful
algebraic manipulator (namely, {\it X$\rho$\'o$\nu o\varsigma$}) he
has developed.  Last but not least, even if my direct contacts with
him have not been very frequent, he has been the mentor of my
supervisor, that is used to repeat: ``Basically, I am trying to
transmit to you what I have learned from Antonio". I think it is more
than enough to be grateful to him.

\section{Introduction}

Lower dimensional elliptic tori can be seen as an extension of stable
periodic solutions, because they are invariant manifolds on which the
motion is quasi-periodic but characterized by a number of frequencies
less than the number of degrees of freedom of the system; furthermore,
the linear approximation of the local expansion of the equations of
motion in the remaining degrees of freedom possesses an elliptic
equilibrium point in correspondence of the origin.  For what concerns
the persistence of elliptic tori in nearly-integrable Hamiltonian
systems of finite dimension, a perturbation theory was firstly
proposed by Melnikov (see~\cite{Melnikov-1965}), even if a complete
proof was lacking; the first proofs were given by
Moser~(\cite{Moser-1967}) and Eliasson~(\cite{Eliasson-1988});
in~\cite{Poschel-1989}, P\"oschel gave another proof of the existence
of lower dimensional elliptic tori using a quadratic perturbative
method and including the case of tori embedded in infinite dimensional
Hamiltonian systems. These results have been further improved
in~\cite{Poschel-1996}, in order to make them suitable to applications
in nonlinear PDEs. More recently, in~\cite{Bia-Chi-Val-2003}
and~\cite{Bia-Chi-Val-2006}, the authors proved the existence of
elliptic tori in two kinds of planetary problems. However, analytical
results concerning the persistence of these invariant objects tend to
be unusable when dealing with the application to specific dynamical
systems, being the hypothesis on the smallness of the perturbation
usually unrealistically small.  Moreover, when dealing with physical
problems, the perturbation is not a free parameter, but it is
given. Therefore, despite the numerical evidence of stable motions,
the existence of invariant manifolds cannot be guaranteed by the mean
of these analytical results, because the problem is not so close to an
integrable one.

For what concerns the KAM theorem, several different proof schemes
have been introduced in the last 50 years with the aim of improving
the applicability to realistic physical problems. In this context,
approaches based on explicit techniques that can be translated in an
algorithm, which can be coded in any programming language, are in a
better position.  In the present work, the attention will be focused
on the particular kind of proofs based on a classical convergence
scheme with respect to the small parameter $\epsilon$
(see~\cite{Gio-Loc-1997} and~\cite{Gio-Loc-1997.1}). The key idea of
these approaches is that explicit computations of proper normal forms
can significantly reduce the size of the perturbation so as to allow
the application of analytical results, otherwise
inaccessible. Moreover, they can be effectively complemented with
computer-assisted estimates (see~\cite{Cel-Gio-Loc-2000} and, for
computer-assisted rigorous results based on a different technique,
e.g.,~\cite{Cel-Chi-1987}).  In the same spirit of what was done for
the existence of full dimensional invariant tori,
in~\cite{San-Loc-Gio-2011} the authors developed a constructive
algorithm for showing the existence of lower dimensional elliptic tori
in the planetary problem. The construction of elliptic tori was
performed by giving the Hamiltonian a suitable normal form, using an
infinite sequence of near-to-the-identity canonical transformations
defined by the Lie series operator. In the subsequent
paper~\cite{Gio-Loc-San-2014}, the complete proof of the convergence
of such a procedure was included. Though, that method and the proof
were designed with the aim of being applied to the study of planetary
systems, where typically the frequencies of the motion on the lower
dimensional torus are faster with respect to the transverse
frequencies, that are indeed generated by the small perturbation.

We adapt that algorithm and, more specifically, the proof so as to
apply it to a general class of Hamiltonians for which the frequencies
of the torus are of the same order as the frequencies describing the
limit of the small oscillations that are transverse to the torus. On
one hand such a difference does not affect the general structure of
the algorithm, on the other hand it has a role in the proof of the
convergence of the construction of the normal form, in particular in
the discussion of the so called {\it geometry of the
  resonances}. Moreover, let us remark that this new formulation of
the theorem includes also the cases in which the two sets of
frequencies are of different order, as the one considered
in~\cite{Gio-Loc-San-2014}. The main obstacle to the convergence is
the introduction of possible small divisors, that are integer
combinations of the two sets of frequencies, i.e., the frequencies of
the invariant torus and the ones related to the motion transverse to
the torus. Furthermore, the algorithm itself introduces small
corrections on both the frequencies, thus making the non-resonance
condition to be verified at each step different and somehow
unpredictable. Since the transverse frequencies have the same weight
with respect to the others and they depend on the frequencies of the
torus, we need to introduce some hypotheses in order to assure the
convergence for an initial set of frequencies of positive measure.
Let us stress that some of the hyphoteses used here could be
weakened. For example, in~\cite{Gio-Loc-San-2014} the authors assumed
the weaker \emph{condition $\vet \tau$} instead of the usual
Diophantine condition in the estimates for the analitic part
(see~\cite{Gio-Mar-2010} and~\cite{Gio-Loc-San-2015} for a more
detailed discussion about the \emph{condition $\vet \tau$} and its
link with Bruno numbers). However, while the set of Diophantine
numbers has positive Lebesgue measure, the one satisfying
\emph{condition $\vet \tau$} that are not Diophantine has zero
Lebesgue measure. Since the result we are presenting here cannot
ensure the convergence for a specific set of initial frequencies and,
therefore, we cannot ensure that the frequencies for which it will
work satisfy such a condition, we decided to not include it here.
In~\cite{Ber-Bia-2011}, the authors showed how also the Melnikov
conditions can be imposed just to the final frequencies. Both these
kind of improvements could be considered by substituting the measure
argument in favour of a statement for specific final
frequencies. Something in this direction has already been done in the
case of full dimensional tori in~\cite{San-Dan}, where the set of
initial frequencies for which the theorem is valid is determined
\emph{a posteriori}.

\subsection{Outline of the algorithm and of the proof}
\label{sec:outline}
By following the approach of P\"oschel~\cite{Poschel-1989}, let us
consider a family of initial Hamiltonians parametrized with respect to
$\vet \omega^{(0)}$ and written as
\begin{equation}
\label{frm:H0-rozza}
\vcenter{\openup1\jot\halign{ \hbox {\hfil $\displaystyle {#}$} &\hbox
    {\hfil $\displaystyle {#}$\hfil} &\hbox {$\displaystyle
      {#}$\hfil}\cr \Hscr^{(0)}(\vet p, \vet q, \vet x,\vet y;\vet
    \omega^{(0)}) =& \vet \omega^{(0)}\cdot \vet p+
    \sum_{j=1}^{n_2}\left[ \frac{\Omega_j^{(0)}(\vet
        \omega^{(0)})}{2}\left(x_j^2+y_j^2\right)\right] \cr & +
    \epsilon\Fscr_0(\vet q;\vet \omega^{(0)}) +\epsilon\Fscr_1(\vet
    q,\vet x,\vet y;\vet \omega^{(0)})+\epsilon\Fscr_2(\vet p,\vet
    q,\vet x,\vet y;\vet \omega^{(0)}) \cr & + \Fscr_{\ge 3}^{\rm (av.)}(\vet p, \vet x,
    \vet y; \vet \omega^{(0)}) + \epsilon \Fscr_{\ge 3}^{\rm (n.av.)}(\vet p,\vet q,\vet x,\vet y;\vet
    \omega^{(0)})\, , \cr }}
\end{equation}
where the variables, split in action-angle variables $(\vet p,\,\vet
q)$ and canonical variables $(\vet x,\,\vet y)$, belong to
neighborhoods of the origin. We collect the Hamiltonian terms in the
functions $\Fscr_j$, where the subscript $j$ denotes the quantity $2
\,{\rm deg}(\vet p)+{\rm deg}(\vet x, \vet y)$. Moreover, the terms of
total degree $j\ge 3$ are split in the averaged terms $\Fscr_{ \ge
  3}^{\rm (av.)}$ and in the part with null average with respect to
the angles $\vet q$, i.e., $\epsilon\Fscr_{\ge 3}^{\rm
  (n. av.)}$. As usual, $\epsilon$ denotes a small parameter.  Here,
$\vet \omega^{(0)}$ is considered to be a parameter with values in
some open domain. Under the hypothesis of the nondegeneracy of the
frequency-action map related to the integrable approximation, we can
expand a given analytic quasi-integrable Hamiltonian as in
formula~\eqref{frm:H0-rozza} by doing expansions for different values
of the actions $\vet p$, that therefore correspond to different
frequency vectors $\vet\omega^{(0)}$. An example of how to define such
a family of Hamiltonians in the case of the FPUT model can be
found in~\cite{Car-Loc-2021}. The search for an elliptic torus
proceeds by applying to the Hamiltonian a perturbative procedure that
is a natural adaptation of the original scheme of Kolmogorov. Let us
give a sketch of it.

We look for a
canonical transformation that gives the Hamiltonian above the normal
form
$$ \vcenter{\openup1\jot\halign{ \hbox {\hfil $\displaystyle {#}$}
    &\hbox {\hfil $\displaystyle {#}$\hfil} &\hbox {$\displaystyle
      {#}$\hfil}\cr \Hscr^{(\infty)}(\vet P,\vet Q,\vet X,\vet Y;\vet
    \omega^{(0)}) &= &\vet \omega^{(\infty)}\cdot \vet P+
    \sum_{j=1}^{n_2}\frac{\Omega_j^{(\infty)}}{2}\left(X_j^2+Y_j^2\right)+
    o\big(\|\vet P\|+\|(\vet X,\vet Y)\|^2\big)\ , \cr }}
$$ where $\vet \omega^{(\infty)}=\vet \omega^{(\infty)}(\vet
\omega^{(0)})$ and $\vet \Omega^{(\infty)}=\vet \Omega^{(\infty)}(\vet
\omega^{(0)})\,$. For this purpose, the transformation should remove all
terms $\Fscr_0(\vet q;\vet \omega^{(0)})$, $\Fscr_1(\vet q,\vet x,\vet
y;\vet \omega^{(0)})$ and $\Fscr_2(\vet p,\vet q,\vet x,\vet y;\vet
\omega^{(0)})$.  Indeed, as a straightforward consequence of the normal form above, the
torus $\vet P=\vet 0$, $\vet X=\vet Y=\vet 0$ is invariant and
elliptic, and it carries a quasi-periodic motion with frequencies $\vet
\omega^{(\infty)}$. 

The proof is divided in two parts: analytic and geometric.  The
analytic part requires an infinite sequence of near to the identity
canonical transformations through the Lie series (see, e.g.,~\cite{Grobner-60}
and~\cite{Giorgilli-2003} for a self-consistent introduction) defined as
\begin{equation}
\label{def:lie-serie}
 \exp \lie {f} (\cdot)  = \sum_{j=0}^{+\infty}\frac{1}{j!}\Lie{f}^j (\cdot) \,,
\end{equation}
where $\lie f(\cdot) = \poisson f \cdot$ is the Lie derivative with
respect to a generic dynamic function $f$, being $\poisson \cdot
\cdot$ the Poisson bracket\footnote{ Let $f$ and $g$ be dynamical
  functions in the canonical variables $(\vet p, \vet q)$, therefore
  the Poisson bracket between them is defined as $ \poisson f g =
  \sum_{j=1}^n \left(\frac{\partial f}{\partial q_j} \frac{\partial
    g}{\partial p_j} - \frac{\partial f}{\partial p_j} \frac{\partial
    g}{\partial q_j}\right).  $}. An accurate control of the
accumulation of the small divisors is needed in order to assure the
convergence of the procedure.

\noindent
Starting with $\Hscr^{(0)}$ as in~\eqref{frm:H0-rozza}, we
introduce an infinite sequence of Hamiltonians
\begin{equation}
\vcenter{\openup1\jot\halign{ \hbox {\hfil $\displaystyle {#}$} &\hbox
    {\hfil $\displaystyle {#}$\hfil} &\hbox {$\displaystyle
      {#}$\hfil}\cr \Hscr^{(r)}&(\vet p,\vet q,\vet x,\vet y;\vet
    \omega^{(0)}) =  \vet \omega^{(r)}\cdot \vet p+
    \sum_{j=1}^{n_2}\left[ \frac{\Omega_j^{(r)}(\vet
        \omega^{(0)})}{2}\left(x_j^2+y_j^2\right)\right] \cr &
    +\epsilon^{r+1}\Fscr_0^{(r)}(\vet q;\vet \omega^{(0)}) 
    +\epsilon^{r+1}\Fscr_1^{(r)}(\vet q,\vet x,\vet y;\vet
    \omega^{(0)})+\epsilon^{r+1}\Fscr_2^{(r)}(\vet p,\vet q,\vet
    x,\vet y;\vet \omega^{(0)}) \cr & + \Fscr_{\ge 3}^{(\rm{av.})(r)}(\vet p, \vet x, \vet y; \vet \omega^{(0)})+\epsilon \Fscr_{\ge 3}^{(\rm{n.av.})(r)}(\vet
    p,\vet q,\vet x,\vet y;\vet \omega^{(0)})\ ,  \cr }}
\label{frm:Hr-rozza}
\end{equation}
such that at every step the size of the perturbing terms is reduced,
as indicated by the factor $\epsilon^{r+1}$. During the iteration of
the algorithm, we introduce possible small divisors, that are integers
combinations of the frequencies $\vet \omega^{(r)}$ and $\vet
\Omega^{(r)}$, which are at each step different. This is an additional
issue to the control of their contribution. Usually, the accumulation
of small divisors is controlled by the mean of a perturbation method
that exhibits a fast convergence (quadratic method), as in the
original paper of Kolmogorov. Here, we follow the classical procedure
that works step-by-step in powers of the parameter $\epsilon$, with
the aim of producing an effective constructive algorithm. The
accumulation of small divisors can be geometrically controlled, as it
follows some strict ``selection rules''. The method is essentially the
same that has been used for the proof of Kolmogorov theorem
in~\cite{Gio-Loc-1997} and~\cite{Gio-Loc-1997.1} and described
in~\cite{Gio-Loc-1999}.  The analytic procedure is followed by a
geometric argument concerning the estimate of the measure of a
suitable set of non-resonant frequencies, based on an adaptation of
the approach described in~\cite{Poschel-1989}.  The scheme is actually
a revisiting of the geometric construction of non-resonant domains
that can be found, e.g.,
in~\cite{Arnold-1963.1},~\cite{Nekhoroshev-1977},~\cite{Nekhoroshev-1979},~\cite{Eliasson-1988}
and~\cite{Giorgilli-1997.4}.

\subsection{Statement of the result}\label{sbs:intro-3}
We come now to a formal statement of our main result.  Let us consider
a $2(n_1+n_2)$-dimensional phase space endowed with $n_1$ pairs of
action-angle coordinates $(\vet p,\vet q)\in\Oscr_1\times\toro^{n_1}$
and $n_2$ pairwise conjugated canonical variables $(\vet x,\vet
y)\in\Oscr_2\subseteq\reali^{2n_2}$, where both
$\Oscr_1\subseteq\reali^{n_1}$ and $\Oscr_2$ are open sets containing
the origin.  We also introduce an open set $\Uscr\subset\reali^{n_1}$
and the frequency vector $\vet \omega^{(0)}\in\Uscr$ which plays the
role of a parameter.

\begin{theorem}\label{thm:main-theorem}Consider a
family of analytic real Hamiltonians as in~\eqref{frm:H0-rozza}, parametrized
by the $n_1$-dimensional frequency vector $\vet \omega^{(0)}$, with
$\epsilon$ playing the usual role of small parameter.  Let us assume
that

\item{(a)} the frequencies $\vet \Omega_j^{(0)}:\Uscr\to\reali$ are
  analytic functions of $\vet \omega^{(0)}\in\Uscr\,$; similarly
  $\Fscr_0\,$, $\Fscr_1\,$, $\Fscr_2$ and $\Fscr_{\ge 3}$ are
  analytic functions of $(\vet p,\vet q,\vet x,\vet y;\vet
  \omega^{(0)})\in\Oscr_1\times\toro^{n_1}\times\Oscr_2\times\Uscr\,$;

\item{(b)} $\Omega_i^{(0)}(\vet \omega^{(0)})\neq\Omega_j^{(0)}(\vet
  \omega^{(0)})$ and $\Omega_{k}^{(0)}(\vet \omega^{(0)})\neq 0$ for
  $\vet \omega^{(0)}\in\Uscr$ and $1\le i<j\le n_2$, $1\le k\le
  n_2\,$;

\item{(c)} the functions $\Fscr_j$ are such that $j = 2p+m$,
  where $p$ is the degree in $\vet p$, while $m$ is the degree in
  $(\vet x, \vet y)$; in particular, $\Fscr_{\ge 3}$ contains terms with degree
  $j$ at least $3$;

\item{(d)}   for some $E>0$ one has
    \begin{equation*}
    \sup_{(\vet p,\vet q,\vet x,\vet y;\vet
          \omega^{(0)})\in\Oscr_1\times\toro^{n_1}\times\Oscr_2\times\Uscr}
        \big|\Fscr_j(\vet p,\vet q,\vet x,\vet y;\vet
        \omega^{(0)})\big| \leq E \qquad \forall\ j\ge 0.
    \end{equation*}

\noindent
Then, there is a positive $\epsilon^{\star}$ such that for $0\le
\epsilon<\epsilon^{\star}$ the following statement holds true: there
exists a non-resonant set $\Uscr^{(\infty)}\subset\Uscr$ of positive
Lebesgue measure and with the measure of $\,
\Uscr\setminus\Uscr^{(\infty)}$ tending to zero for $\epsilon\to 0$
for bounded $\Uscr$, such that for each $\vet
\omega^{(0)}\in\Uscr^{(\infty)}$ there are functions $\vet
  \omega^{(\infty)}=\vet \omega^{(\infty)}(\vet \omega^{(0)})$, $\vet
  \Omega^{(\infty)}=\vet \Omega^{(\infty)}(\vet \omega^{(0)})\,$ and an analytic canonical transformation $(\vet p,\vet
q,\vet x,\vet y)=\psi_{\vet \omega^{(0)}}^{(\infty)}(\vet P,\vet
Q,\vet X,\vet Y)$ leading the Hamiltonian in the normal form
$$ \vcenter{\openup1\jot\halign{ \hbox {\hfil $\displaystyle {#}$}
    &\hbox {\hfil $\displaystyle {#}$\hfil} &\hbox {$\displaystyle
      {#}$\hfil}\cr \Hscr^{(\infty)}(\vet P,\vet Q,\vet X,\vet Y;\vet
    \omega^{(0)}) &= &\vet \omega^{(\infty)}\cdot \vet P+
    \sum_{j=1}^{n_2}\frac{\Omega_j^{(\infty)}}{2}\left(X_j^2+Y_j^2\right)+
    o\big(\|\vet P\|+\|(\vet X,\vet Y)\|^2\big)\ .\cr }}
$$
\end{theorem}
We remark that the main difference with respect to the result reported
in~\cite{Gio-Loc-San-2014} is the lack of the small parameter
$\epsilon$ in front of the terms which describe the motion transverse
to the $n_1$-dimensional torus, as these transverse oscillation are of
the same order of the quasi-periodic motion. For the sake of
simplicity, in the statement above we omitted any explicit estimate of
the threshold value $\epsilon^\star$, which can be found
in~\eqref{frm:soglia-finale}. In turn, this estimate depends on a set
of parameters introduced in lemma~\ref{lem:Ham-expansion}, in
proposition~\ref{pro:geometric-prop} and in
formul{\ae}~\eqref{frm:def-M},~\eqref{frm:soglia-analitica}
and~\eqref{def:bstorto}.

The whole proof requires a lot of technical estimates and it is
organized as follows.  In section~\ref{sec:algorithm} we describe the
formal algorithm. Section~\ref{sec:proof} is devoted to the proof of
theorem~\ref{thm:main-theorem} and it is composed of two different
parts: the analytic part of the proof contains the quantitative
estimates that lead to the proof of convergence, with a special
emphasis to the control of small divisors; in the geometric part of
the proof, we show that our procedure applies to a set of initial
frequencies of large relative measure. Finally, in
section~\ref{sec:final-discussions} we comment the result and the
possible applications.

\section{Formal algorithm}\label{sec:algorithm}
This section is devoted to the formal algorithm that takes a
Hamiltonian~\eqref{frm:H0-rozza} of the family $\Hscr^{(0)}$,
parametrized by the frequency vector $\vet\omega^{(0)}$, and brings it
into a normal form.  Therefore, the section includes all the
formul{\ae} that will be used in order to establish the convergence of
the normalization process. We found convenient to use the complex
variables $(\vet z, i\bar{\vet z})$, defined by the canonical
transformation $\vet z=(\vet x+ i \vet y)/\sqrt{2}$, in order to deal
with the transverse oscillations.

\subsection{Initial settings and strategy
of the formal algorithm}\label{sbs:algo-1} For some fixed positive
integer $K$, we introduce the distinguished classes of functions
$\PPset_{\hat m,\hat \ell}^{sK}\,$, with integers $\hat m,\,\hat
\ell,\,s\ge 0\,$, which can be written as
\begin{equation}
g(\vet p,\vet q,\vet z, \imunit\bar {\vet z}) = \sum_{{\scriptstyle{\vet
      m\in\naturali^{n_1}}}\atop{\scriptstyle{|\vet m|=\hat {\vet
        m}}}} \,\sum_{{\scriptstyle{(\vet \ell,\bar {\vet
        \ell})\in\naturali^{2n_2}}}\atop{\scriptstyle{|\vet
      \ell|+|\bar {\vet \ell}|=\hat{ \vet \ell}}}}
\,\sum_{{\scriptstyle{{\vet
        k\in\interi^{n_1}}}\atop{\scriptstyle{|\vet k|\le sK}}}}
c_{\vet m,\vet \ell,\bar {\vet \ell},\vet k} \,\vet p^{\vet m} \vet
z^{\vet \ell}(\imunit\bar {\vet z})^{\bar {\vet \ell}} \exp(\imunit \vet k\cdot
\vet q)\ ,
\label{frm:esempio-g-in-PPset}
\end{equation}
with coefficients $c_{\vet m,\vet \ell,\bar {\vet \ell},\vet
  k}\in\complessi$.  Here we denote by $|\cdot|$ the $\ell_1$-norm and
we adopt the multi-index notation, i.e., $\vet p^{\vet
  m}=\prod_{j=1}^{n_1} p_j^{m_j}$.  We say that
$g\in\Pset_{\ell}^{sK}$ in case
$$ g\in\bigcup_{{\scriptstyle {\hat m\ge 0\,,\,\hat l\ge 0}}\atop
  {\scriptstyle {2\hat m+\hat l=\ell}}}\PPset_{\hat m,\hat\ell}^{sK} \ .
$$

\noindent
Finally, we denote by $\langle g\rangle_{\vet \vartheta} =
\int_{\toro^n}\diff\vartheta_1\ldots\diff\vartheta_n \,g/(2\pi)^n$ the
average of a function $g$ with respect to the angles $\vet
\vartheta\in\toro^n$.  We shall also omit the dependence of the
function from the variables, unless it has some special meaning.

In the following lemma we state a relevant algebraic property what we are going to use several times.
\begin{lemma}\label{lem:pb-classes-functions}
Let $g\in\Pset_{\ell}^{sK}$
and $g^{\prime}\in\Pset_{\ell^{\prime}}^{s^{\prime}K}$ for some
$\ell,\,s,\,\ell^{\prime},\,s^{\prime}\ge 0$ and $K>0\,$. Then
$\{g,g^{\prime}\}\in\Pset_{\ell+\ell^{\prime}-2}^{(s+s^{\prime})K}\,$.
\end{lemma}
\noindent
The straightforward proof is left to the reader.

We start with the Hamiltonian in the form

\begin{equation}
\vcenter{\openup1\jot\halign{ \hbox {\hfil $\displaystyle {#}$} &\hbox
    {\hfil $\displaystyle {#}$\hfil} &\hbox {$\displaystyle
      {#}$\hfil}\cr H^{(0)} &= \scalprod{\vet \omega^{(0)}}{\vet p} +
    \sum_{j=1}^{n_2}\Omega^{(0)}_{j}z_{j}{\bar z}_{j}
    +\sum_{\ell>2}\sum_{s\geq 0} \epsilon^{s} f_{\ell}^{(0,s)} \cr
    &\quad+\sum_{s\geq 1} \epsilon^{s} f_{0}^{(0,s)}+ \sum_{s\geq 1}
    \epsilon^{s} f_{1}^{(0,s)}+ \sum_{s\geq 1} \epsilon^{s}
    f_{2}^{(0,s)}\ , \cr }}
\label{frm:H(0)}
\end{equation}
where $f_{\ell}^{(0,s)}\in\Pset_{\ell}^{sK}$. The
Hamiltonian~\eqref{frm:H0-rozza} may be written in the
form~\eqref{frm:H(0)} (see section~\ref{sec:proof}).

Starting from $H^{(0)}$, we introduce an infinite sequence of Hamiltonians
$\left\{H^{(r)}\right\}_{r\geq 0}$ with $H^{(r)}$ in normal form up to
order $r\,$, according to the following description.  We transform $H^{(r-1)}$
into $H^{(r)}$ via a near-the-identity canonical transformation
generated by a composition of three Lie series of the form
\begin{equation}
H^{(r)} = H^{(r-1)} \circ \exp\left( \epsilon^r\Lie{\chi_0^{(r)}} \right)\circ
\exp\left( \epsilon^r\Lie{\chi_1^{(r)}} \right)\circ \exp\left(
\epsilon^r\Lie{\chi_2^{(r)}} \right) \ ,
\label{frm:serie-Lie-passo-r}
\end{equation}
where 
$\chi_0^{(r)}(\vet q)\in\Pset_{0}^{rK}\,$, $\chi_1^{(r)}(\vet q,\vet z,
\imunit\bar{\vet z})\in\Pset_{1}^{rK}\,$ and $\chi_2^{(r)}(\vet p,\vet q,\vet
z,\imunit \bar{\vet z})\in\Pset_{2}^{rK}$.
The generating functions $\chi_0^{(r)}$, $\chi_1^{(r)}$ and
$\chi_2^{(r)}$ are determined by homological equations.  At every
normalization step the frequencies $\vet \omega^{(r)}$ and $\vet
\Omega^{(r)}$ may change by a small quantity (see
formul{\ae}~\eqref{frm:chgfreq.toro} and~\eqref{frm:chgfreq.trasv}).

The small divisors are controlled by introducing a non-resonance
condition up to a finite order $rK$, namely
\begin{equation}
\,\min_{{\scriptstyle{\vet k\in\interi^{n_1}\,,\,0 <|\vet k|\leq
      rK}\atop\scriptstyle{\vet \ell\in\interi^{n_2}\,,\,0\leq
      |\vet \ell|\leq2}}}\big| \scalprod{\vet k}{\vet \omega^{(r-1)}(\vet
  \omega^{(0)})} +\scalprod{\vet \ell}{\vet \Omega^{(r-1)}(\vet
  \omega^{(0)})}\big| \ge \frac{\gamma}{(rK)^\tau}\ ,
\label{frm:nonres_a}
\end{equation}
and 
\begin{equation}
\,\min_{{\scriptstyle{0< |\vet \ell| \le 2}}}
\big|\vet \ell \cdot \vet \Omega^{(r-1)}(\vet \omega^{(0)}) \big| \ge \gamma\ ,
\label{frm:nonres_b}
\end{equation}
where $\gamma>0$ and $\tau\ge n_1-1$.  For $|\vet \ell|=0$
condition~\eqref{frm:nonres_a} is the usual condition of strong
non-resonance, while for $|\vet \ell|= 1,2$ it is usually referred to
as the first and second Melnikov condition, respectively.

We come to the description of the generic $r$-th normalization step.
Let us write the Hamiltonian $H^{(r-1)}$ as

\begin{equation}
\vcenter{\openup1\jot\halign{ \hbox {\hfil $\displaystyle {#}$} &\hbox
    {\hfil $\displaystyle {#}$\hfil} &\hbox {$\displaystyle
      {#}$\hfil}\cr H^{(r-1)} &= \scalprod{\vet \omega^{(r-1)}}{\vet
      p} + \sum_{j=1}^{n_2}\Omega^{(r-1)}_{j}z_{j}{\bar z}_{j}
    +\sum_{\ell>2}\sum_{s\geq 0} \epsilon^{s} f_{\ell}^{(r-1,s)} \cr
    &\quad+\sum_{s\geq r} \epsilon^{s} f_{0}^{(r-1,s)}+ \sum_{s\geq r}
    \epsilon^{s} f_{1}^{(r-1,s)}+ \sum_{s\geq r} \epsilon^{s}
    f_{2}^{(r-1,s)}\ , \cr }}
\label{frm:H(r-1)}
\end{equation}
where $f_{\ell}^{(r-1,s)}\in\Pset_{\ell}^{sK}$.  In the expansion above,
the functions $f_{\ell}^{(r-1,s)}$ may depend analytically on
$\epsilon$, the relevant information being that they carry a common
factor $\epsilon^s$.  Such an expansion is clearly not unique, but
this is harmless.

\subsection{First stage of the normalization step}\label{sbs:step1}
Our aim is to remove the term $f_{0}^{(r-1,r)}$.  We determine the
generating function $\chi^{(r)}_{0}$ by solving the homological
equation
\begin{equation}
\Lie{\chi^{(r)}_{0}} \left(\scalprod{\vet \omega^{(r-1)}}{\vet
  p}\right) + f_{0}^{(r-1,r)} - \langle f_{0}^{(r-1,r)}\rangle_{\vet
  q} = 0\ .
\label{frm:chi0r}
\end{equation}
Considering the Taylor-Fourier expansion
\begin{equation}
f_{0}^{(r-1,r)}(\vet q)=
\sum_{{\scriptstyle{|\vet k|\le rK}}}
c_{\vet 0,\vet 0,\vet 0,\vet k}^{(r-1)}\exp(\imunit\vet  k\cdot \vet  q)\ ,
\label{frm:espansione-f_0^(r-1,r)}
\end{equation}
we readily get
\begin{equation}
\chi^{(r)}_{0}(\vet q)=\sum_{{\scriptstyle{0<|\vet k|\le rK}}}
\frac{c_{\vet 0,\vet 0,\vet 0,\vet k}^{(r-1)}}{\imunit\scalprod{\vet k}{\vet \omega^{(r-1)}}}
\exp(\imunit \vet k\cdot \vet q)\ .
\label{frm:espansione-chi0r}
\end{equation}
The denominators are not zero in view of the non-resonance
condition~\eqref{frm:nonres_a} with $|\vet \ell|=0\,$.  The new
Hamiltonian is determined as the Lie series with generating function
$\epsilon^r \chi^{(r)}_{0}$, namely

\begin{equation}
\vcenter{\openup1\jot\halign{ \hbox {\hfil $\displaystyle {#}$} &\hbox
    {\hfil $\displaystyle {#}$\hfil} &\hbox {$\displaystyle
      {#}$\hfil}\cr H^{(\rmI;r)} &=
    \exp\left(\epsilon^{r}\Lie{\chi^{(r)}_{0}}\right)H^{(r-1)} \cr &=
    \scalprod{\vet \omega^{(r-1)}}{\vet p} +
    \sum_{j=1}^{n_2}\Omega^{(r-1)}_{j}z_{j}{\bar z}_{j}
    +\sum_{\ell>2}\sum_{s\geq 0} \epsilon^{s} f_{\ell}^{(\rmI;r,s)}
    \cr &\quad+\sum_{s\geq r} \epsilon^{s} f_{0}^{(\rmI;r,s)}+
    \sum_{s\geq r} \epsilon^{s} f_{1}^{(\rmI;r,s)}+ \sum_{s\geq r}
    \epsilon^{s} f_{2}^{(\rmI;r,s)}\ .  \cr }}
\label{frm:H(I;r)-espansione}
\end{equation}
The functions $f_{\ell}^{(\rmI;r,s)}$ are recursively defined as
\begin{equation}
\vcenter{\openup1\jot \halign{
    $\displaystyle\hfil#$&$\displaystyle{}#\hfil$&$\displaystyle#\hfil$\cr
    f_{0}^{(\rmI;r,r)} &= 0\ ,&\cr \noalign{\smallskip}
    f_{0}^{(\rmI;r,r+m)} &= f_{0}^{(r-1,r+m)} &\quad\hbox{for }
    0<m<r\,,\cr \noalign{\smallskip} f_{\ell}^{(\rmI;r,s)}
    &=\sum_{j=0}^{\lfloor s/r\rfloor} \frac{1}{j!}
    \Lie{\chi^{(r)}_{0}}^{j}f^{(r-1,s-jr)}_{\ell+2j}& \quad\hbox{for
    }{\vtop{\hbox{${\ell =0\,,\ s\ge 2r\ \hbox{ or }
            \ \ell=1,\,2\,,\ s\ge r}$}
\vskip-2pt\hbox{\hskip-5pt$\hbox{or }\ \ell\ge 3\,,\ s\ge 0\ ,$}}}\cr
}}
\label{frm:fI}
\end{equation}
with $f_{\ell}^{(\rmI;r,s)}\in\Pset_{\ell}^{sK}$.  The constant term
$c_{\vet 0,\vet 0,\vet 0,\vet 0}^{(r-1)}= \langle
f_{0}^{(r-1,r)}\rangle_{\vet q}$ has been omitted. Indeed, it does not
affect the Hamilton equations, but it is just the level of the energy
of the solution we are going to construct.

\subsection{Second stage of the normalization step}\label{sbs:step2}
We remove $f_{1}^{(\rmI;r,r)}$ from~\eqref{frm:H(I;r)-espansione}.  We
determine a second generating function $\chi^{(r)}_{1}$ by solving the
homological equation
\begin{equation}
\Lie{\chi^{(r)}_{1}}\biggl(\scalprod{\vet \omega^{(r-1)}}{\vet p}+
\sum_{j=1}^{n_2}\Omega^{(r-1)}_{j}z_{j}{\bar z}_{j}\biggr)
+f_{1}^{(\rmI;r,r)}= 0\ .
\label{frm:chi1r}
\end{equation}
Again, if we consider the Taylor-Fourier expansion
\begin{equation}
f_{1}^{(\rmI;r,r)}(\vet q,\vet z,\imunit \bar {\vet z})=\sum_{|\vet
  \ell|+|\bar {\vet \ell}|=1}\, \sum_{{\scriptstyle{|\vet k|\le
      rK}}}\, c_{\vet 0, \vet \ell,\bar{\vet \ell},\vet
  k}^{(\rmI;r)}\vet z^{\vet \ell}(\imunit\bar {\vet z})^{\bar {\vet
    \ell}}\exp(\imunit\scalprod{\vet k}{\vet q})\ ,
\label{frm:espansione-f_1^(I;r,r)}
\end{equation}
then the solution of the homological equation is
\begin{equation}
\chi^{(r)}_{1}(\vet q,\vet z, \imunit\bar{\vet z})=\sum_{|\vet
  \ell|+|\bar {\vet \ell}|=1}\, \sum_{{\scriptstyle{|\vet k|\le rK}}}
\,\frac{c_{\vet 0,\vet \ell,\bar{ \vet \ell},\vet k}^{(\rmI;r)}\,\vet
  z^{\vet \ell}(\imunit\bar{\vet z})^{\bar {\vet \ell}}
  \exp(\imunit\vet k\cdot \vet q)}{\imunit\big[\scalprod{\vet k}{\vet
      \omega^{(r-1)}}+ \scalprod{(\vet \ell-\bar{\vet \ell})}{\vet
      \Omega^{(r-1)}}\big]}\ ,
\label{frm:espansione-chi1r}
\end{equation}
where the divisors cannot vanish in view of condition~\eqref{frm:nonres_a}
and of the condition $\Omega_i^{(r-1)} \neq 0$.

The new Hamiltonian is calculated as
\begin{equation}
H^{(\rmII;r)} = \exp\left(\epsilon^{r}\Lie{\chi^{(r)}_{1}}\right) H^{(\rmI;r)}\ ,
\label{frm:H(II;r)}
\end{equation}
and is given the form~\eqref{frm:H(I;r)-espansione}, replacing
the upper index ${\rm I}$ by ${\rm II}\,$, with
\begin{equation}
\vcenter{\openup1\jot \halign{
    $\displaystyle\hfil#$&$\displaystyle{}#\hfil$&$\displaystyle#\hfil$\cr
    f_{\ell}^{(\rmII;r,r)} &= 0 &\quad\hbox{for } \ell=0,1\,,\cr
    \noalign{\smallskip} f_{\ell}^{(\rmII;r,r+m)} &=
    f_{\ell}^{(\rmI;r,r+m)} &\quad\hbox{for } \ell=0,1\,,\ 0<m<r\,,\cr
    \noalign{\smallskip} f_{0}^{(\rmII;r,2r)} &= f_{0}^{(\rmI;r,2r)} +
    \frac{1}{2} \Lie{\chi_{1}^{(r)}}f_{1}^{(\rmI;r,r)}\ ,&\cr
    \noalign{\smallskip} f_{0}^{(\rmII;r,2r+m)} &=
    f_{0}^{(\rmI;r,2r+m)} + \Lie{\chi_{1}^{(r)}}f_{1}^{(\rmI;r,r+m)}
    &\quad\hbox{for } 0<m<r\,,\cr \noalign{\smallskip}
    f_{\ell}^{(\rmII;r,s)} &= \sum_{j=0}^{\lfloor s/r\rfloor}
    \frac{1}{j!} \Lie{\chi^{(r)}_{1}}^{j}f^{(\rmI;r,s-jr)}_{\ell+j}&
    \quad\hbox{for }{\vtop{\hbox{${\ell =0\,,\ s\ge 3r\ \hbox{ or }
            \ \ell=1\,,\ s\ge 2r}$}
\vskip-2pt\hbox{\hskip-5pt$\hbox{or }\ \ell=2\,,\ s\ge r\ \hbox{ or }\ \ell\ge 3\,,\ s\ge 0\,.$}}}\cr
}}
\label{frm:fII}
\end{equation}

\subsection{Third stage of the normalization step}\label{sbs:step3}
We remove here at the same time the $\vet q$-dependent part and the
average one of $f_{2}^{(\rmII;r,r)}$.

Considering again the Taylor-Fourier expansion

\begin{equation}
\vcenter{\openup1\jot\halign{ \hbox {\hfil $\displaystyle {#}$} &\hbox
    {\hfil $\displaystyle {#}$\hfil} &\hbox {$\displaystyle
      {#}$\hfil}\cr f_{2}^{(\rmII;r,r)}(\vet p,\vet q,\vet
    z,\imunit\bar{\vet z}) &=\sum_{|\vet m|=1}\,
    \sum_{{\scriptstyle{|\vet k|\le rK}}}\, c_{\vet m,\vet 0,\vet
      0,\vet k}^{(\rmII;r)}\vet p^{\vet m}\exp(\imunit\scalprod{\vet
      k}{\vet q}) \cr &\phantom{=}+\sum_{|\vet \ell|+|\bar {\vet
        \ell}|=2}\, \sum_{{\scriptstyle{|\vet k|\le rK}}}\, c_{\vet
      0,\vet \ell,\bar{\vet \ell},\vet k}^{(\rmII;r)}\vet z^{\vet
      \ell}(\imunit\bar {\vet z})^{\bar {\vet \ell}} \exp(\imunit
    \scalprod{\vet k}{\vet q})\ , \cr }}
\label{frm:espansione-f_2^(II;r,r)}
\end{equation}
we determine the generating function $\chi^{(r)}_2$ by solving the
homological equation
\begin{equation}
\Lie{\chi^{(r)}_{2}} \biggl(\scalprod{\vet \omega^{(r-1)}}{\vet p}+
\sum_{j=1}^{n_2} \Omega^{(r-1)}_{j}z_{j}{\bar z}_{j}\biggr)
+ f_{2}^{(\rmII;r,r)} = Z^{(r)}\ .
\label{frm:chi2r}
\end{equation}
with
\begin{equation}
Z^{(r)}= \sum_{|\vet m|=1}c_{\vet m,\vet 0,\vet 0,\vet 0}^{(\rmII;r)}
\vet p^{\vet m} + \sum_{|\ell|=1}c_{\vet 0,\vet \ell,\vet \ell,\vet
  0}^{({\rmII};r)}\vet z^{\vet \ell}(\imunit\bar {\vet z})^{\vet
  \ell}\ ,
\label{frm:struttura-Z^(r)}
\end{equation}
where the coefficients $c_{\vet m,\vet 0,\vet 0,\vet 0}^{(\rmII;r)} $
and $c_{\vet 0,\vet \ell,\vet \ell,\vet 0}^{({\rmII};r)}\,$ are the
same as those related to the monomials $\vet p^{\vet m}$ and $\vet
z^{\vet \ell}(\imunit\bar{\vet z})^{\vet \ell}$ appearing in the
expansion of $ f_{2}^{(\rmII;r,r)}$
in~\eqref{frm:espansione-f_2^(II;r,r)}.

\noindent
We get
\begin{equation}
\vcenter{\openup1\jot\halign{ \hbox {\hfil $\displaystyle {#}$} &\hbox
    {\hfil $\displaystyle {#}$\hfil} &\hbox {$\displaystyle
      {#}$\hfil}\cr \chi^{(r)}_{2}(\vet p,\vet q,\vet z, \imunit\bar
    {\vet z}) &=\sum_{|\vet m|=1}\,\sum_{{\scriptstyle{0<|\vet k|\le
          rK}}} \,\frac{c_{\vet m,\vet 0,\vet 0,\vet
        k}^{(\rmII;r)}\vet p^{\vet m}\exp(\imunit \vet k\cdot \vet q)}
    {\imunit\scalprod{\vet k}{\vet \omega^{(r-1)}}} \cr
    &\phantom{=}+\sum_{|\vet \ell|+|\bar {\vet \ell}|=2}\,
    \sum_{{\scriptstyle{0<|\vet k|\le rK}}} \,\frac{c_{\vet 0,\vet
        \ell,\bar {\vet \ell},\vet k}^{(\rmII;r)}\vet z^{\vet
        \ell}(\imunit\bar {\vet z})^{\bar {\vet \ell}}
      \exp(\imunit\vet k\cdot \vet q)}{\imunit \big[\scalprod{\vet
          k}{\vet \omega^{(r-1)}}+ \scalprod{(\vet \ell-\bar{\vet
            \ell})}{\vet \Omega^{(r-1)}}\big]} \cr
    &\phantom{=}+\sum_{{\scriptstyle{|\vet \ell|+|\bar {\vet
            \ell}|=2}}\atop{\scriptstyle{\vet \ell\neq\bar {\vet
            \ell}}}} \,\frac{c_{\vet 0,\vet \ell,\bar {\vet \ell},\vet
        0}^{({\rmII};r)}}{\imunit\scalprod{(\vet \ell-\bar{\vet
          \ell})}{\vet \Omega^{(r-1)}}} \vet z^{\vet \ell}(\imunit\bar
        {\vet z})^{\bar {\vet \ell}} \ , \cr }}
\label{frm:espansione-chi2r}
\end{equation}
where the divisors cannot vanish in view of
conditions~\eqref{frm:nonres_a} and~\eqref{frm:nonres_b}.

The terms of the Hamiltonian can now be calculated as
\begin{equation}
\vcenter{\openup1\jot \halign{
    $\displaystyle\hfil#$&$\displaystyle{}#\hfil$&$\displaystyle#\hfil$\cr
    f_{\ell}^{(r,r)} &= 0 &\quad\hbox{for } \ell=0,1,2\,,\cr
    \noalign{\smallskip} f_{\ell}^{(r,s)} &= \sum_{j=0}^{\lfloor
      s/r\rfloor-1}\frac{1}{j!} \Lie{\chi^{(r)}_{2}}^{j}
    f^{(\rmII;r,s-jr)}_{\ell}&\quad\hbox{for } \ell=0,1 \,,\ s>r\,,\cr
    \noalign{\smallskip} f_{2}^{(r,jr)} &=
    \frac{j-1}{j!}\Lie{\chi^{(r)}_{2}}^{j-1}
    f^{(\rmII;r,r)}_{2}+\sum_{i=0}^{j-2}\frac{1}{i!}
    \Lie{\chi^{(r)}_{2}}^{i} f^{(\rmII;r,(j-i)r)}_{2}&\quad\hbox{for }
    j\geq2\,,\cr \noalign{\smallskip} f_{2}^{(r,jr+m)} &=
    \sum_{i=0}^{j-1}\frac{1}{i!}  \Lie{\chi^{(r)}_{2}}^{i}
    f^{(\rmII;r,(j-i)r+m)}_{2}&\quad\hbox{for } j\geq 1
    \,,\ 0<m<r\,,\cr \noalign{\smallskip} f_{\ell}^{(r,s)} &=
    \sum_{j=0}^{\lfloor s/r\rfloor} \frac{1}{j!}
    \Lie{\chi^{(r)}_{2}}^{j} f^{(\rmII;r,s-jr)}_{\ell}& \quad\hbox{for
    } \ell\geq 3\,,\ s\geq 0\,.  \cr }}
\label{frm:f}
\end{equation}
The average terms that cannot be eliminated contained in $f_2^{({\rm
    II}; r,r)}$ are recollected in the frequencies $\vet \omega^{(r)}$
and $\vet \Omega^{(r)}$, which are therefore corrected as:
\begin{equation}
 \omega_{j}^{(r)}=
 \omega_{j}^{(r-1)}+\epsilon^{r}\parder{Z^{(r)}}{p_j}
\qquad
\hbox{for } j=1,\,\ldots,\,n_1
\label{frm:chgfreq.toro}
\end{equation}
and 
\begin{equation}
\Omega_{j}^{(r)}=\Omega_{j}^{(r-1)}+\epsilon^r
\frac{\partial^2 Z^{(r)}}{\partial z_j\partial (\imunit{\bar z}_j)}
\qquad
\hbox{for } j=1,\,\ldots,\,n_2\ 
\label{frm:chgfreq.trasv}
\end{equation}
with $Z^{(r)}$ as in~\eqref{frm:struttura-Z^(r)}.

The transformed Hamiltonian is given the
form~\eqref{frm:H(r-1)}, replacing $r-1$ with $r$, namely
\begin{equation}
\vcenter{\openup1\jot\halign{
 \hbox {\hfil $\displaystyle {#}$}
&\hbox {\hfil $\displaystyle {#}$\hfil}
&\hbox {$\displaystyle {#}$\hfil}\cr
H^{(r)} &= \scalprod{\vet \omega^{(r)}}{\vet p} +
\sum_{j=1}^{n_2}\Omega^{(r)}_{j}z_{j}{\bar z}_{j}
+\sum_{\ell>2}\sum_{s\geq 0} \epsilon^{s} f_{\ell}^{(r,s)}
\cr
&\quad+\sum_{s\geq r+1} \epsilon^{s} f_{0}^{(r,s)}+
\sum_{s\geq r+1} \epsilon^{s} f_{1}^{(r,s)}+
\sum_{s\geq r+1} \epsilon^{s} f_{2}^{(r,s)}\ ,
\cr
}}
\label{frm:H(r)-espansione}
\end{equation}
eventually with a change of the frequencies $\vet \omega^{(r)}$ and
$\vet \Omega^{(r)}$.

Let us add a few considerations. First, we remark that our formulation of the algorithm works for both real and complex Hamiltonians. However, if the
expansion~\eqref{frm:H(0)} contains only real functions, then all
terms of type $\scalprod{\vet \omega^{(r)}}{\vet p}\,$,
$\sum_{j=1}^{n_2}\Omega^{(r)}_{j}z_{j}{\bar z}_{j}$ and
$f_{\ell}^{(r,s)}$ generated by the algorithm are real too, as it can be easily
checked.
Finally, the Hamiltonian $H^{(r)}$ in~\eqref{frm:H(r)-espansione} has
the same form of $H^{(r-1)}$, so that the induction step can be
iterated provided the conditions~\eqref{frm:nonres_a} and~\eqref{frm:nonres_b}
hold true with $r+1$ in place of $r\,$.

\section{Proof of the convergence of the algorithm}\label{sec:proof}
\subsection{Analytic settings}
We introduce the complex domains $\Dscr_{\rho, \sigma, R, h}
= \Gscr_{\rho} \times \toro^{n_1}_{\sigma} \times \Bscr_{R} \times \Wscr_h\,$,
where $\Gscr_{\rho}\subset\complessi^{n_1}$ and
$\Bscr_R\subset\complessi^{n_2}\times\complessi^{n_2}$ are open balls
centered at the origin with radii $\rho$ and $R$,
respectively, $\Wscr$ is a subset of $\reali^{n_1}$, while the
subscripts $\sigma$ and $h$ denote the usual complex
extensions of real domains (see~\cite{Giorgilli-2003}), i.e., 
\begin{equation}
\Gscr_{\rho}=\big\{\vet z\in\complessi^{n_1}:\max_{1\le j\le
n_1}|z_j|<\rho\big\}, 
\end{equation}
\begin{equation}
\toro^{n_1}_{\sigma}=\big\{\vet q\in\complessi^{n_1}:\realpart
q_j\in\toro,\break\ \max_{1\le j\le n_1}|\imaginary q_j|<\sigma\big\}\,, 
\end{equation}
\begin{equation}
\Bscr_R=\{\vet z\in\complessi^{2n_2}: \max_{1\le j\le 2n_2} |z_j|<R\,\}
\end{equation}
and 
\begin{equation}
\Wscr_h = \big\{\vet z\in\complessi^{n_1}: \exists\ \vet \omega\in\Wscr\,,
\ \max_{1\le j\le n_1}|z_j-\omega_j|<h\big\}.
\end{equation} 
Let us consider a generic analytic function $g:\Dscr_{\rho,
\sigma,R, h}\to\complessi$,
\begin{equation}
g(\vet p,\vet q,\vet z,\imunit\bar {\vet z}; \vet \omega) =
\sum_{{\scriptstyle{\vet k\in\interi^{n_1}}}} g_{\vet k}(\vet p,\vet z,\imunit\bar{\vet  z};\vet \omega)
\exp(\imunit \vet k\cdot \vet q)\ ,
\label{frm:funz}
\end{equation}
where $g_{\vet k}:\Gscr_{\rho}\times \Bscr_{R} \times \Wscr_h\to\complessi\,$.
We define the weighted Fourier norm
$$
\|g\|_{\rho,
\sigma,R, h}=\sum_{{\scriptstyle{\vet k\in\interi^{n_1}}}}
\big|g_{\vet k}\big|_{\rho,R,h}\exp(|\vet k|\sigma)\ ,
$$
where
\begin{equation}
\left|g_{\vet k}\right|_{\rho,R,h}=
\sup_{{\scriptstyle{p\in\Gscr_{\rho}}}
      \atop{{\scriptstyle{(\vet z,\imunit \bar {\vet z})\in \Bscr_{R}}}
      \atop{\scriptstyle{\vet\omega\in \Wscr_h}}}}
\big|g_{\vet k}(\vet p,\vet z,\imunit\bar{\vet  z};\vet \omega)\big|\ .
\label{frm:norma}
\end{equation}
It is convenient to introduce the Lipschitz constant
for the Jacobian of the function
$\vet \Omega^{(0)}:\,\Wscr^{(0)}_{h_0}\to\complessi^{n_2}$, defined as
\begin{equation}
\left| \frac{\partial\vet \Omega^{(0)}}{\partial\vet \omega^{(0)}} \right|_{\infty;\Wscr^{(0)}_{h_0}}=
\,\sup_{{\scriptstyle{\vet \omega^{(0)}\in \Wscr^{(0)}_{h_0}}}}\,
\sup_{{{\scriptstyle{\vet \beta\neq 0}}}\atop{\scriptstyle{\vet \beta+\vet \omega^{(0)}\in \Wscr^{(0)}_{h_0}}}}\,
\max_{1\le i\le n_2 \atop{1\le j\le n_1}}\frac{\left|\Omega^{(0)}_i(\vet \omega^{(0)}+\vet \beta)-
\Omega^{(0)}_i(\vet \omega^{(0)})\right|}{\left|\beta_j\right|}\ .
\label{frm:cost-Lipschitz}
\end{equation}

We notice that the dependence on the parameter $\vet \omega^{(0)}$ plays no role
in the following, thus we shorten the notation by ignoring it.
\begin{lemma}\label{lem:Ham-expansion}
Let us assume the same hypotheses of theorem~\ref{thm:main-theorem}
over the family of Hamiltonians $\Hscr^{(0)}$. Then, there exist
positive parameters $\rho\,$, $R\,$, $\sigma\,$, $h_0\,$, $\gamma\,$,
$\tau\,$, $J_0\,$, $\Theta_0\,$, $\Ebarra\,$, a compact set
$\Wscr^{(0)}\subset\reali^{n_1}$ and a positive integer value $K$ such
that the canonical change of coordinates $(\vet p,\vet q,\vet
z,\imunit\bar{\vet z})\mapsto(\vet p,\vet q,\vet x,\vet y)$ transforms
$\Hscr^{(0)}$ in the Hamiltonian
$H^{(0)}:\Dscr_{\rho,\sigma,R}\times\Wscr^{(0)}_{h_0}\to\complessi$
described by the expansion~\eqref{frm:H(0)}, where both $\vet
\Omega^{(0)}(\vet \omega^{(0)})$ and all the terms of type
$f_{\ell}^{(0,s)}$ are real analytic functions of $\vet
\omega^{(0)}\in\Wscr^{(0)}_{h_0}\,$. Moreover, the following
properties are satisfied:

\item{(a')} the initial set $\Wscr^{(0)}$ of frequencies  is
  {\it non-resonant} up to the finite order $K\,$, namely every
  $\vet \omega^{(0)}\in\Wscr^{(0)}$ satisfies
  $$
  \,\min_{\scriptstyle{\vet k\in\interi^{n_1}\,,\,0<|\vet k|\leq K}\atop
    \scriptstyle{\vet \ell\in\interi^{n_2}\,,\,0\leq |\vet \ell|\leq2}} \big|
  \scalprod{\vet k}{\vet \omega^{(0)}} +
  \scalprod{\vet \ell}{\vet \Omega^{(0)}(\vet \omega^{(0)})}\big|\,>
  \frac{2\gamma}{K^\tau}
  $$
  and
  $$
  \,\min_{\scriptstyle{\,\vet \ell\in\interi^{n_2}\,,\,0<|\vet \ell|\leq2}}
  \big| \scalprod{\vet \ell}{\vet \Omega^{(0)}(\vet \omega^{(0)})} \big|\,>
  2 \gamma \,;
  $$
\item{(b')} the Jacobian of
  $\vet \Omega^{(0)}(\vet \omega^{(0)})$ is uniformly bounded in the extended domain
  $\Wscr^{(0)}_{h_0}\,$, namely
  $\big| \partial \vet \Omega^{(0)}/\partial \vet \omega^{(0)}\big|_{\infty;\Wscr^{(0)}_{h_0}} \leq
  J_0 < \infty\,$; moreover, it is such that
  $$
  {\rm dist}(\vet k,\Kscr_{\vet \ell}^{(0)})\ge 2\Theta_0
  \qquad
  \forall\ \vet k\in\interi^{n_1}\setminus\{\vet 0\}\,,\ \vet \ell\in\interi^{n_2},\ 0<|\vet \ell|\leq2\ ,
  $$
  being ${\rm dist}$ the usual euclidean distance from a vector to
  a set and $\Kscr_{\vet \ell}^{(0)}$ the closed convex hull of the gradient set
  $$
  \Gscr_{\vet \ell}^{(0)}=
  \left\{\partial_{\vet \omega^{(0)}}\big[\scalprod{\vet \ell}{\vet \Omega^{(0)}\big(\vet \omega^{(0)}\big)}\big]
  \>:\ \vet \omega^{(0)}\in\Wscr^{(0)}\right\}\ ;
  $$

\item{(c')} $f_{\ell}^{(0,s)}\in\Pset_{\ell}^{sK}$;

\item{(d')} the following upper bounds hold true:
  $$
  \left\|f_{\ell}^{(0,s)}\right\|_{\rho,\sigma,R}\le \frac{\Ebarra}{2^{\ell}}\ .
  $$
\end{lemma}
\noindent
We give a sketch of the proof. The Hamiltonian terms should be split
in different parts, each of them with a finite expansion. For any
fixed value of index $\ell\,$, standard arguments on the Fourier decay
of the coefficients allow us to determine suitable values of the
parameters $K$ and $\sigma\,$, such that for $s\ge 0$ the norms of the
functions $f_{\ell}^{(0,s)}$ are uniformly bounded (see, e.g., the
proof of lemma~5.2 in~\cite{Giorgilli-2003}). For a fixed $K$, we can
determine $\gamma>0\,$, $\tau> n_1-1\,$, $h_0>0$ and a compact set
$\Wscr\subset\Uscr$ such that a Diophantine inequality in~(a') is
satisfied, being $\Uscr$ the initial set of frequency vectors $\vet
\omega^{(0)}$ in the hypotheses of theorem~\ref{thm:main-theorem}. The
first property in (b') is a straightforward consequence of the
analyticity of $\Omega_j^{(0)}$ on the domain $\Uscr$ and the second
property is satisfied if we choose $\Wscr^{(0)}$ small enough. Indeed,
by reducing $\Wscr^{(0)}$ we could reduce also the gradient set
$\Gscr_{\vet \ell}^{(0)}$ in view of the analiticity of $\vet
\Omega^{(0)}$ in $\Wscr^{(0)}$.  Property~(c') is a consequence of the
fact that we can organize different terms of the Taylor-Fourier
expansion of the Hamiltonian. Finally, the inequality at point~(d') is
satisfied for appropriate values of $\rho$ and $R\,$ determined by
standard arguments on the Taylor expansions of homogeneous
polynomials\footnote{An example of how to estimate $\bar E$ can
    be found in lemma~5.2 of~\cite{Giorgilli-2003}. }.  Taking into
account that the usual sup-norm is bounded by the weighted Fourier
one, for $\epsilon<1$ we have
$$
\sum_{\ell\,,\,s}\epsilon^s
\|f_{\ell}^{(0,s)}\|_{\rho,\sigma,R}\le 2\Ebarra/(1-\epsilon)\ .
$$
This implies that the Hamiltonian $H^{(0)}$ is {\it analytic} in
$\Dscr_{\rho,\sigma,R}\times\Wscr_{h_0}\,$.

We give now some estimates for the Lie series.  We shorten
the notation by writing $\|\cdot\|_{\alpha}\,$ in place of the norm
$\|\cdot\|_{\alpha(\rho,\sigma,R)}$, where $\alpha$ is a positive
number.
\begin{lemma}\label{lem:Lie-der-estimate}
Let $d$ and $d^{\prime}$ be real numbers such that $d>0\,$, $d^{\prime}\ge 0$ and
$d+d^{\prime}<1\,$; let $\Xscr$ and $g$ be two analytic functions on
$\Dscr_{(1-d^{\prime})(\rho,\sigma,R)}$
having finite norms $\|\Xscr\|_{1-d^{\prime}}$ and
$\|g\|_{1-d^{\prime}}\,$, respectively. Then, for $j\ge 1$, we have
\begin{equation}
\frac{1}{j!}\left\|\Lie{\Xscr}^{j}g\right\|_{1-d-d^{\prime}}
\le\frac{1}{\ee^{2}}
\left(\frac{2e}{\rho\sigma}+\frac{e^2}{R^{2}}\right)^{j}
\frac{1}{d^{2j}}\|\Xscr\|^{j}_{1-d^{\prime}}\|g\|_{1-d^{\prime}}\ .
\label{frm:stimalie}
\end{equation}
\end{lemma}
\noindent
Estimates similar to~\eqref{frm:stimalie} are contained, e.g., in~\cite{Giorgilli-2003}. Nevertheless, a little additional
work is needed in order to adapt them to the present context.

\subsection{Analytic part}\label{sec:analytic-2}
In this section, we translate our formal algorithm into a recursive
scheme of estimates on the norms of the functions involved in the
normalization process.

\subsubsection{Small divisors}\label{sbs:indices}
It is well known that the accumulation of the small divisors is one of the major obstacles to the convergence of perturbative proof schemes. In order to control it, firstly, we introduce lists of indices $S = \{s_1, \ldots, s_j\}$ that may contain repeated elements and we give some definitions that will be useful later.
\begin{definition}\label{def:dioph-freq}
For any fixed $\vet \omega^{(0)}\,$, we say that
the sequence of frequency vectors $$\big\{\big(\vet \omega^{(r)}(\vet \omega^{(0)})\,,\,
\vet \Omega^{(r)}(\vet \omega^{(0)})\big)\big\}_{r\ge 0}$$ is
Diophantine, if there are two constants $\gamma>0$ and $\tau\ge n_1-1$
 such that~\eqref{frm:nonres_a} and~\eqref{frm:nonres_b}
are satisfied, for all $r\ge 1$.
\end{definition}

\begin{definition}
Let $S$ be a list of indices. We define the evaluation operator $\Vscr$ on $S$ as follows:
\begin{equation}
 \Vscr(S) = \prod_{s\in S} {s^{4+\tau}}
\label{def:eval-op}
\end{equation}
where $\tau$ is related to the Diophantine condition in definition~\ref{def:dioph-freq}.

\noindent
Moreover, let $S_1, \ldots, S_n$ be lists of indices, the operator $\Mscr\Ascr\Xscr$ extracts the list of indices that maximizes the evaluation operator, i.e.,
\begin{equation}
\Mscr\Ascr\Xscr(S_1, \ldots, S_n) = S_m,\quad   where  \ 1 \le m \le n \ \ and\ \  \Vscr( S_m) \ge \Vscr(S_j) \quad  \forall \ 1 \le j \le n.
\label{def:max-list}
\end{equation}
\noindent
Finally, we introduce the operator $\Nscr_k$ on a list $S$, which counts the numbers of elements in $S$ with value in a quadratic range, that is
\begin{equation}
\Nscr_k(S) = \# \{ s \in S : 2^k \le s < 2^{k+1} \} \quad \forall \, k \ge 0 \, . 
\label{def:counter-list}
\end{equation}
\label{def:list-indices}
\end{definition}

Now, let us describe the mechanism of
accumulation of the small divisors in a rather informal way.  The
key remark is that the small divisors are created by the
solution of a homological equation and accumulate through Lie
derivatives.  For instance, when we solve~\eqref{frm:chi0r} with
$r=1$, in $\chi_0^{(1)}$ we introduce a divisor which can be estimated by $a_1 = \gamma/K^\tau$ according to condition~\eqref{frm:nonres_a}.  We can therefore associate to $\chi_0^{(1)}$ the list $\{1\}$.
Similarly, at order $r$ if the list $\{j_1,\ldots,j_w\}$ is associated
to $f_0^{(r-1,r)}$, then $\chi_0^{(r)}$ possesses the associated list
$\{j_1,\ldots,j_w,r\}$.  Concerning the Poisson bracket, suppose that
the functions $\chi$ and $f$ possess the associated lists
$\{j_1,\ldots,j_w\}$ and $\{j'_1,\ldots,j'_v\}$, respectively.  Then
$\Lie\chi f$ possesses the list
$$
\{j_1,\ldots,j_w\} \cup \{j'_1,\ldots,j'_v\} = 
\{j_1,\ldots,j_w,j'_1,\ldots,j'_v\}\ .
$$
For a sum of several Lie derivatives, we select the greater list,
according to the definition of the operator $\Mscr \Ascr \Xscr$ in~\eqref{def:max-list}.  
Thus, the process of
accumulation of divisors is described via the union of lists.  An
heuristic analysis of the mechanism of accumulation, with the aim of
identifying the worst one, could be made by unfolding all the
recursive definitions of the functions $f_{\ldots}^{(\ldots)}$
in~\eqref{frm:fI}, \eqref{frm:fII} and~\eqref{frm:f}, in consequence of the three transformations of coordinates introduced at every normalization step. The process is substantially similar to  the simpler cases of the Birkhoff normal form or the Kolmogorov normal form.  Using  as a reference the description in~\cite{Gio-Loc-1999} and~\cite{Loc-Mel-1998}, we define a valid selection rule that allows us to control geometrically the product of small divisors.

A second source of divergence is connected with
the restriction of domains required by the estimate of multiple
Poisson brackets.  The restrictions are controlled by a sequence $d_r$
that should converge to some positive $d$ small enough.  We define the sequences
$\{d_r\}_{r\geq0}$ and $\{\delta_r\}_{r\ge 1}$ as
\begin{equation}
d_{0}=0\ ,\qquad
d_{r}=d_{r-1}+3\delta_{r}\ ,\qquad
\delta_{r}=\frac{1}{2\pi^{2} r^{2}}\ .
\label{frm:dr-and-deltar}
\end{equation}
At the $r$-th step of the algorithm, we find an Hamiltonian $H^{(r)}$
which is analytic in $\Dscr_{(1-d_r)(\rho,\sigma,R)}\,$, namely in a domain restricted by $\delta_r$.  Nevertheless, the parameters
can be chosen in a way that $\lim_{r\to\infty}d_{r}=1/4$ and, therefore, the sequence
of domains converges to a compact set with no empty interior.

Finally, in order to give quantitative estimates of the accumulation of small
divisors, we are going to use the operator $\Nscr_k(\cdot)$ that counts how many indices of a specific range we have in a list.

\subsubsection{Convergence of the algorithm under
non-resonance conditions}\label{sbs:convergence}
It is convenient to introduce the constant
\begin{equation}
M=\max\left\{1\,,\,\frac{4\pi^4 \Ebarra K^\tau}{ \gamma}\left(\frac{2e}{\rho\sigma}+\frac{e^2}{R^{2}}\right)\,\right\}\ .
\label{frm:def-M}
\end{equation}
This will allow us to produce uniform estimates in the parameters.
The number of summands
in the recursive formul{\ae}~\eqref{frm:fI}, \eqref{frm:fII}
and~\eqref{frm:f} will be estimated by the three sequences of
integers $\{\nu_{r,s}\}_{r\ge 0\,,\,s\ge 0}$,
$\{\nu_{r,s}^{(\rm I)}\}_{r\ge 1\,,\,s\ge 0}$ and
$\{\nu_{r,s}^{(\rm II)}\}_{r\ge 1\,,\,s\ge 0}$ defined as
\begin{equation}
\vcenter{\openup1\jot 
\halign{
$\displaystyle\hfil#$&$\displaystyle{}#\hfil$&$\displaystyle#\hfil$\cr
\nu_{0,s} &= 1
&\quad\hbox{for } s\ge 0\,,
\cr
\nu_{r,s}^{(\rm I)} &= \sum_{j=0}^{\lfloor s/r \rfloor} \nu_{r-1,r}^{j}\nu_{r-1,s-jr}
&\quad\hbox{for } r\ge 1\,,\ s\ge 0\,,
\cr
\nu_{r,s}^{(\rm II)} &= \sum_{j=0}^{\lfloor s/r \rfloor}
(\nu_{r,r}^{(\rm I)})^{j}\nu_{r,s-jr}^{(\rm I)}
&\quad\hbox{for } r\ge 1\,,\ s\ge 0\,,
\cr
\nu_{r,s} &= \sum_{j=0}^{\lfloor s/r \rfloor}
(\nu_{r,r}^{(\rm II)})^{j}\nu_{r,s-jr}^{(\rm II)}
&\quad\hbox{for } r\ge 1\,,\ s\ge 0\,.
\cr
}}
\label{frm:seqnu}
\end{equation}

\begin{lemma}\label{lem:iterative-estimates}
Consider a Hamiltonian
$H^{(0)}$of type~\eqref{frm:H(0)}.  
 Assume that the
non-resonance hypotheses~\eqref{frm:nonres_a} and~\eqref{frm:nonres_b} are
satisfied up to order $r\ge 1$, for some $\gamma>0$ and $\tau\ge n_1-1$.
Then, the formal algorithm described in section~\ref{sec:algorithm} can be iterated up to introduce  $H^{(r-1)}$ expanded as in~\eqref{frm:H(r-1)} with $f_\ell^{(r-1,s)} \in \Pscr_{\ell}^{sK}$ and also the $r$-th normalization step can be performed. Moreover, we assume that the terms of the Hamiltonian $H^{(r-1)}$ are bounded so that
\begin{equation}
\vcenter{\openup1\jot 
\halign{
$\displaystyle\hfil#$&$\displaystyle{}#\hfil$&$\displaystyle#\hfil$\cr
\|f_{\ell}^{( r-1,s)}\|_{1-d_{r-1}}&\leq \frac{\Ebarra M^{3s-\zeta_\ell}}{2^{\ell}}\,  
\Vscr \left(\Hscr_\ell^{(r-1,s)}\right)
\nu_{r-1,s}
\quad\,&{\hbox{for }0\le\ell\le 2\,,\ s \ge r\,,} \ 
\cr
& &\hbox{or }{\vtop{\hbox{${\ell\ge 3\,,\ s\ge 0\,,}$}
}}\cr
}}
\label{frm:f_l^(r-1,s)-lemmone}
\end{equation}
where $\zeta_\ell = 3-\ell$ if $\ell\le 3 $, $\zeta_\ell =0 $ otherwise, and $\Vscr$ is the operator defined in~\eqref{def:eval-op}. 
\noindent
Therefore, the norms of the generating functions of the $r$-th step are subject to the following limitations:
\begin{equation}
\vcenter{\openup1\jot 
\halign{
$\displaystyle\hfil#$&$\displaystyle{}#\hfil$&$\displaystyle#\hfil$\cr
\left(\frac{2e}{\rho\sigma}+\frac{e^2}{R^{2}}\right)\frac{1}{\delta_r^2}
\,\|\chi_{0}^{(r)}\|_{1-d_{r-1}}&\leq M^{3r-2}\,  
\Vscr\left(\Gscr_0^{(r)}\right) \,\nu_{r-1,r} \ , 
\cr
\left(\frac{2e}{\rho\sigma}+\frac{e^2}{R^{2}}\right)\frac{1}{\delta_r^2}
\,\|\chi_{1}^{(r)}\|_{1-d_{r-1}-\delta_{r}}&\leq M^{3r-1}\,  
\Vscr\left(\Gscr_1^{(r)}\right)
\,\nu_{r,r}^{(\rm I)}\ ,
\cr
\left(\frac{2e}{\rho\sigma}+\frac{e^2}{R^{2}}\right)\frac{1}{\delta_r^2}
\,\|\chi_{2}^{(r)}\|_{1-d_{r-1}-2\delta_{r}}&\leq M^{3r} 
\Vscr\left(\Gscr_2^{(r)}\right)
\, \nu_{r,r}^{(\rm II)}\ .
\cr
}}
\label{frm:generatrici-lemmone}
\end{equation}

\noindent
Furthermore, the terms appearing in the expansion of the new Hamiltonian
$H^{(r)}$ in~\eqref{frm:H(r)-espansione} are bounded by
\begin{equation}
\vcenter{\openup1\jot 
\halign{
$\displaystyle\hfil#$&$\displaystyle{}#\hfil$&$\displaystyle#\hfil$\cr
\|f_{\ell}^{( r,s)}\|_{1-d_{r}}&\leq \frac{\Ebarra M^{3s-\zeta_\ell}}{2^{\ell}}\,  
\Vscr \left(\Hscr_\ell^{(r,s)}\right)
\nu_{r,s}\
\quad\,&{\hbox{for }0\le\ell\le 2\,,\, s>r\ } 
\cr & &
\hbox{or }{\vtop{\hbox{${\ell\ge 3\,,\ s\ge 0\,.}$}
}}\cr
}}
\label{frm:f_l^(r,s)-lemmone}
\end{equation}
\noindent
The lists of positive integer numbers $\big\{ \Hscr_\ell^{(r,s)}\big\}$,and $\big\{ \Gscr_j^{(r)}\big\}_{0\le j \le 2}$ are recursively defined as follows:
\begin{equation}
\vcenter{\openup1\jot 
\halign{
$\displaystyle\hfil#$&$\displaystyle{}#\hfil$&$\displaystyle#\hfil$\cr
\Hscr_\ell^{(0,s)} &= \emptyset \qquad &
{\hbox {for}} \ s \ge 1\ , \ \ell \ge 0,
\cr
\Hscr_\ell^{({\rm I} ;\,r,s)} &= { \Mscr \Ascr \Xscr}_{ 0 \le j \le \lfloor s/r\rfloor } \left\{   \cup_{i=1}^j   \Gscr_0^{(r)} \cup \Hscr_{\ell+2j}^{(r-1,s-jr)} \right\}
\ \ & {\hbox {for }}\ \ell \ge 0, \ s \ge 2 r 
\cr & & \ {\hbox {or }} \ \ell \ge 1, \ s \ge 0,
\cr
\Hscr_0^{({\rm I} ;\,r,s)} &= \Hscr_0^{(r-1,s)}  \quad &{\hbox {for }} \  r < s < 2r, 
\cr
\Hscr_\ell^{({\rm II};\,r,s)} &= { \Mscr \Ascr \Xscr}_{ 0 \le j \le \lfloor s/r\rfloor }\left\{   \cup_{i=1}^j   \Gscr_1^{(r)} \cup \Hscr_{\ell+j}^{({\rm I};\,r,s-jr)} \right\} 
& {\hbox {for }} \ \ell=0, \ s  \ge 3r 
\cr & & \ {\hbox {or }} \ \ell=1, \ s\ge 2r 
\cr & & \ {\hbox {or }} \ \ell \ge 2, \ s \ge r,
\cr
\Hscr_\ell^{({\rm II} ;\,r,s)} &= { \Mscr \Ascr \Xscr}_{ 0 \le j \le \lfloor s/r\rfloor-1 }\left\{   \cup_{i=1}^j   \Gscr_1^{(r)} \cup \Hscr_{\ell+j}^{({\rm I};\,r,s-jr)} \right\}\ &  {\hbox {for }} \ \ell=0, \ r < s < 3r
\cr & & \ {\hbox {or }}\ \ell=1, \ r < s < 2r, 
\cr
\Hscr_\ell^{(r,s)} &= { \Mscr \Ascr \Xscr}_{ 0 \le j \le \lfloor s/r\rfloor } \left\{  \cup_{i=1}^j   \Gscr_2^{(r)} \cup \Hscr_\ell^{({\rm II};\,r,s-jr)} \right\} \ & {\hbox {for }} \ \ell \ge 3, \ s \ge 0,
\cr
\Hscr_\ell^{(r,kr+m)} &={ \Mscr \Ascr \Xscr}_{ 0 \le j \le k-1 } \left\{  \cup_{i=1}^j   \Gscr_2^{(r)} \cup \Hscr_\ell^{({\rm II};\, r,(k-j)r+m)} \right\} \quad
& {\hbox {for }} \ \ k \ge 1, \ 0 <m <r, 
\cr &&  \ell=0,1,2,
\cr
\Hscr_\ell^{(r,kr)} &=  { \Mscr \Ascr \Xscr}_{ 0 \le j \le k-2 }\left\{  \cup_{i=1}^j   \Gscr_2^{(r)} \cup \Hscr_\ell^{({\rm II};\,r,(k-j)r)} \right\} \ & {\hbox {for }} \ \ k \ge 2, \ \ell=0,1,
\cr
\Hscr_2^{(r,kr)} &=   { \Mscr \Ascr \Xscr}_{ 0 \le j \le k-1 } \left\{  \cup_{i=1}^j   \Gscr_2^{(r)} \cup \Hscr_2^{({\rm II};\,r,(k-j)r)} \right\} \ & {\hbox {for }} \ \ k \ge 1,
\cr
\Gscr_0^{(r)} &=  \Hscr_0^{(r-1,r)} \cup \{r\} \,,
\cr
\Gscr_1^{(r)} &=   \Hscr_1^{({\rm I}; \, r,r)} \cup \{r\} \,,
\cr
\Gscr_2^{(r)} &=  \Hscr_2^{({\rm II};\, r,r)} \cup \{r\} \,.
\cr
}}
\label{frm:H_l^(r,s)-lemmone}
\end{equation}

\noindent 
Finally, for $r\ge 1\,$, the variations of the frequencies induced by
the $r$-th normalization step, are bounded by
\begin{equation}
\max_{{1\le i\le n_1}\atop{1\le j\le n_2}}
\left\{\frac{1}{\sigma}\big|\omega_{i}^{(r)}-\omega_{i}^{(r-1)}\big|\,,\,
\big|\Omega_{j}^{(r)}-\Omega_{j}^{(r-1)}\big|\right\}\le \gamma
\epsilon^{r}M^{3r} \Vscr\left(\Gscr_2^{(r)}\right)
\nu_{r,r} \, .
\label{frm:stima-variazione-frequenze}
\end{equation}
\end{lemma}

\noindent
The proof is made by induction on $r$ and it is deferred to appendix~\ref{app:stimeiterative}.

\begin{lemma}\label{lem:small-div-acc}
The following rules hold true for the list of indices defined in~\eqref{frm:H_l^(r,s)-lemmone} for all $r \ge 1$ and $s \ge r$ when $\ell= 0, \,1, \,2$, or $s \ge 0$ if $\ell \ge 3$:
\begin{equation}
\Nscr_k\left(\Hscr_\ell^{(r,s)}\right)  \le \left\{
\vcenter{\openup1\jot 
\halign{
$\displaystyle#\hfil$&\quad$\displaystyle{}#\hfil$&\quad$\displaystyle#\hfil$\cr
 3 \left\lfloor \frac {s}{2^k} \right\rfloor  &  {\rm for} \ k < \lfloor \log_2{r} \rfloor \>, 
\cr
 3 \left\lfloor \frac {s}{2^k} \right\rfloor -3 +\ell  &  {\rm for} \ \ell \le 3, \, k =  \lfloor \log_2{r} \rfloor \>,
\cr
 3 \left\lfloor \frac {s}{2^k} \right\rfloor  &  {\rm for} \  \ell > 3, \, k =  \lfloor \log_2{r} \rfloor \>,
\cr
0 & {\rm for}  \ k > \lfloor \log_2{r} \rfloor \>,
\cr
}}
\right .
\label{frm:Nk_Ham}
\end{equation}
\begin{equation}
\Nscr_k\left(\Gscr_j^{(r)}\right) \le \left\{
\vcenter{\openup1\jot 
\halign{
$\displaystyle#\hfil$&\quad$\displaystyle{}#\hfil$&\qquad$\displaystyle#\hfil$\cr
  3 \left\lfloor \frac {r}{2^k} \right\rfloor  & {\rm for} \ k< \lfloor \log_2{r} \rfloor \>,
\cr
j+1 & {\rm for}  \  k=\left\lfloor \log_2{r}\right\rfloor \>, & \forall \ j=0,\,1, \,2 \ .
\cr
0 & {\rm for}  \ k > \lfloor \log_2{r} \rfloor \>,
\cr
}}
\right .
\label{frm:Nk_gen}
\end{equation}
\end{lemma}
The proof is made by induction and it is deferred to appendix~\ref{app:accumulodiv}.

We can estimate the sequence $\{\nu_{r,s}\}_{r\ge 0\,,\,s\ge 0}$, counting the number of summands generated by the algorithm, as in the following lemma. 
\begin{lemma}\label{lem:nu-estimate}
The sequence of positive integer numbers
$\{\nu_{r,s}\}_{r\ge 0\,,\,s\ge 0}$ defined
in~\eqref{frm:seqnu} is bounded by
$$
\nu_{r,s}\le\nu_{s,s}\leq 2^{8s}
\qquad
\hbox{for } r\ge 0\,,\ s\ge 0\ .
$$
\end{lemma}
The proof is available in~\cite{Gio-Loc-San-2014}.
\begin{lemma}\label{lem:stime-valut-insiemi}
Consider a Hamiltonian $H^{(0)}$of type~\eqref{frm:H(0)}.   Assume that the non-resonance hypotheses~\eqref{frm:nonres_a} and~\eqref{frm:nonres_b} are
satisfied up to order $r\ge 1$, with $\gamma>0$ and $\tau\ge n_1-1$. Then the formal algorithm described in section~\ref{sec:algorithm} can be iterated up to introduce  $H^{(r-1)}$ expanded as in~\eqref{frm:H(r-1)}, with $f_\ell^{(r-1,s)} \in \Pscr_{\ell}^{sK}$, and also the $r$-th normalization step can be performed.
Moreover, let us suppose that conditions~(c') and~(d') of lemma~\ref{lem:Ham-expansion} are satisfied.

\noindent
Therefore, at  step $r$ of the normalization procedure, the generating functions are bounded by: 
\begin{equation}
\left(\frac{2e}{\rho\sigma}+\frac{e^2}{R^{2}}\right)\frac{1}{\delta_r^2}
\,\|\chi_j^{(r)}\|_{1-d_{r-1}-j\delta_r} \le \Ascr^r  \qquad {\textit for} \ j=0,1,2;
\label{frm:bound-generatrici}
\end{equation}
the functions $f_\ell^{(r,s)}$, with $s>r$ for $ 0 \le \ell \le 2$, $s \ge 0$ for $\ell\ge 3$, are bounded by:
\begin{equation}
\label{frm:bound-fl}
\| f_\ell^{(r,s)} \| _{1-d_r}\le  \frac {\Ebarra}{2^\ell} \Ascr^s\,;
\end{equation}
finally, the sequences $\{\vet \omega^{(r)}\}_{r \ge 0}$ and $\{\vet \Omega^{(r)}\}_{r \ge 0}$ satisfy:
\begin{equation}
\max_{1\le i\le n_1}
\big|\omega_{i}^{(r)}-\omega_{i}^{(r-1)}\big|\le
\gamma\sigma\big(\epsilon\Ascr\big)^r\ ,
\qquad
\max_{1\le j\le n_2}
\big|\Omega_{j}^{(r)}-\Omega_{j}^{(r-1)}\big|\le 
\gamma\big(\epsilon \Ascr\big)^r\ ,
\label{frm:stima-variazione-frequenze-+-esplicita}
\end{equation}
where $\Ascr = M^3 2^{12(4+\tau)}2^8$ and $M$ is defined as in~\eqref{frm:def-M}.

\end{lemma}
\begin{proof}
Thanks to the hypotheses~(c') and~(d') of lemma~\ref{lem:Ham-expansion}, it is straightforward that the terms of the Hamiltonian at step $0$ statisfy~\eqref{frm:f_l^(r-1,s)-lemmone}. 
Then, lemma~\ref{lem:iterative-estimates} applies and we can iterate it up to $r$ times. Therefore, since we can use the estimates in~\eqref{frm:generatrici-lemmone},~\eqref{frm:f_l^(r,s)-lemmone} and~\eqref{frm:stima-variazione-frequenze}, in order to prove the new inequalities in the statement we need to provide an upper bound for $\Vscr(\cdot)$. Using the definition~\eqref{def:eval-op} and the rules for the list of indices in~\eqref{frm:Nk_gen}, we have the following estimate:
$$
\Vscr\left(\Gscr_0^{(r)}\right) = \prod_{s \in \Gscr_0^{(r)} } s^{4+\tau} \le 
 \prod_{k = 0}^\infty \left( 2^{(k+1)(4+\tau)} \right)^{\Nscr_k\left(\Gscr_0^{(r)}\right)}\,.
$$
It is convenient to work with $\log_2 \Vscr\big(\Gscr_0^{(r)}\big)$:
\begin{equation}
\vcenter{\openup1\jot 
\halign{
$\displaystyle\hfil#$&$\displaystyle{}#\hfil$&$\displaystyle#\hfil$\cr
\log_2{\Vscr\left(\Gscr_0^{(r)}\right)} = & \log_2{\left[\prod_{s \in \Gscr_0^{(r)} } s^{4+\tau}\right]}
\le \log_2{\left[ \prod_{k = 0}^\infty \left( 2^{(k+1)(4+\tau)} \right)^{\Nscr_k\left(\Gscr_0^{(r)}\right)}\right]}
\cr
\le &
\sum_{k = 0}^\infty \Nscr_k\left(\Gscr_0^{(r)}\right) ((4+\tau)(k+1)) 
\cr
 \le & 3r (4+\tau)  \sum_{k = 0}^\infty \frac{(k+1)}{2^k}
 \le 12r(4+\tau) \,.
\cr
}}
\label{frm:stima_Vscr}
\end{equation}
Thus, it follows that $\Vscr \big(\Gscr_0^{(r)}\big) \le 2^{12r(4+\tau)} $; the same inequality is valid for $\Vscr\big(\Gscr_1^{(r)}\big)$ and $\Vscr\big(\Gscr_2^{(r)}\big)$.
For the generating functions, it suffices to observe that we can estimate $\nu_{r-1,r}$, $\nu_{r,r}^{(\rmI)}$ and $\nu_{r,r}^{(\rmII)}$  with $\nu_{r,r}$ and finally use the estimate provided in lemma~\ref{lem:nu-estimate}. This allows to fully justify~\eqref{frm:bound-generatrici}.

\noindent
The proof of~\eqref{frm:bound-fl} is very similar; indeed, we have an analogous estimate for $\log_2 \Vscr\big(\Hscr_\ell^{(r,s)}\big)$ as that in~\eqref{frm:stima_Vscr}, by replacing $r$ with $s$. Using that $\nu_{r,s} \le \nu_{s,s}$, we can verify~\eqref{frm:bound-fl} .

\noindent
We have now to prove~\eqref{frm:stima-variazione-frequenze-+-esplicita}. Using again~\eqref{frm:stima-variazione-frequenze} and the estimates for $\Vscr(\Gscr_2^{(r)})$ and $\nu_{r,r}$, the inequality directly follows.
This concludes the proof.
\end{proof}

We summarize the analytic estimates in the following
\begin{proposition}\label{pro:analytic-prop}
Consider an analytic Hamiltonian
$H^{(0)}:\Dscr_{\rho,\sigma,R}\to\complessi\,$ expanded as
in~\eqref{frm:H(0)}, satisfying hypotheses~(c')--(d') of
lemma~\ref{lem:Ham-expansion}.

\noindent
Assume moreover the following hypotheses:

\item{(e')} the sequences of frequency vectors $\{\vet \omega^{(r)}\}_{r\ge 0}$
and $\{\vet \Omega^{(r)}\}_{r\ge 0}$ are Diophantine in the sense of definition~\ref{def:dioph-freq}.

\item{(f')} the parameter $\epsilon$ is smaller than the ``analytic
threshold value'' $\epsilon^{\star}_{{\rm an}}\,$, being
\begin{equation}
\epsilon^{\star}_{{\rm an}}= \frac{1}{\Ascr}
\qquad \hbox{with} \qquad
\Ascr=M^3  2^{56+12 \tau}\,,
\label{frm:soglia-analitica}
\end{equation}
where $M$ is defined in~\eqref{frm:def-M} and $\tau$ is related to the Diophantine condition in~\ref{def:dioph-freq}.

\noindent
Then, there exists an analytic canonical transformation
$\Upphi_{\vet \omega^{(0)}}^{(\infty)}:\Dscr_{1/2(\rho,\sigma,R)}\to
\Dscr_{3/4(\rho,\sigma,R)}$ such that the
Hamiltonian $H^{(\infty)}=H^{(0)}\circ\Upphi_{\vet \omega^{(0)}}^{(\infty)}$ is in
normal form, i.e.,
\begin{equation}
H^{(\infty)} = \scalprod{\vet \omega^{(\infty)}}{\vet p} +
\sum_{j=1}^{n_2}\Omega^{(\infty)}_{j}z_{j}{\bar z}_{j}
+\sum_{\ell>2}\sum_{s\geq 0} \epsilon^{s} f_{\ell}^{(\infty,s)}\ .
\label{frm:H(infty)-espansione}
\end{equation}
Furthermore, the norms of the functions $f_{\ell}^{(\infty,s)}\in\Pset_{\ell}^{sK}$ are bounded by
\begin{equation}
\left\|f_{\ell}^{(\infty,s)}\right\|_{3/4}\le
\frac{\Ebarra}{2^{\ell}}\Ascr^s
\qquad\hbox{for } \ell\ge 3\,,\ s\ge 0\,,
\label{frm:f_l^(infty,s)}
\end{equation}
and the sequences $\{\vet \omega^{(r)}\}_{r\ge 0}$ and $\{\vet \Omega^{(r)}\}_{r\ge 0}$ converge to the limits 
$\vet \omega^{(\infty)}$ and $\vet \Omega^{(\infty)}$, respectively.
\end{proposition}
\noindent
The rest of the section contains a sketch of the proof, which depends
on all the previous lemmas.

First, thanks to the hypothesis~(e') of non-resonance  and the condition~(c'),~(d') of lemma~\ref{lem:Ham-expansion}, the estimates in lemma~\ref{lem:stime-valut-insiemi} hold true at every step $r$.
Remarking that $\epsilon\Ascr<1$ we have that $H^{(r)}$,
written as in~\eqref{frm:H(r)-espansione}, is analytic on
$\Dscr_{(1-d_r)(\rho,\sigma,R)}\supset\Dscr_{3/4(\rho,\sigma,R)}\,$.
With similar calculations, in view of
condition~(f'), since the inequalities
in~\eqref{frm:stima-variazione-frequenze-+-esplicita} hold true at every step $r$, we can conclude that the sequences $\{\vet \omega^{(r)}\}_{r\ge 0}$ and $\{\vet \Omega^{(r)}\}_{r\ge 0}$ are Cauchy sequences, so they have limits.

Let us now focus on the difference $H^{(r)}-H^{(r-1)}$. Using~\eqref{frm:H(r)-espansione} and \eqref{frm:H(r-1)} and the fact that $f_\ell^{(r,r)}=0$ for $\ell=0,1,2$ and $f_\ell^{(r-1,s)} = f_\ell^{(r,s)} \ \forall \ \ell >2$, $0 \le s <r$ , we get
\begin{equation}
\vcenter{\openup1\jot 
\halign{
$\displaystyle#\hfil$&\quad$\displaystyle{}#\hfil$&\qquad$\displaystyle#\hfil$\cr
H^{(r)}-H^{(r-1)}=&\scalprod{\big(\vet \omega^{(r)}-\vet \omega^{(r-1)}\big)}{\vet p} +
\sum_{j=1}^{n_2}\big(\Omega^{(r)}_{j}-\Omega^{(r-1)}_{j}\big)
z_{j}{\bar z}_{j}
\cr
& +\sum_{\ell\ge 0}\sum_{s\geq r} \epsilon^{s} 
\left(f_{\ell}^{(r,s)}-f_{\ell}^{(r-1,s)}\right) \,.
\cr
}}
\label{frm:H(r)menoH(r-1)-espansione}
\end{equation}
Therefore, in view of the inequalities in~\eqref{frm:bound-fl}--\eqref{frm:stima-variazione-frequenze-+-esplicita},
we have
\begin{equation}
\left\|H^{(r)}-H^{(r-1)}\right\|_{3/4}\le
\left(n_1\gamma \sigma\rho+n_2\gamma R^2
+\frac{4\Ebarra}{1-\epsilon\Ascr}\right)\big(\epsilon\Ascr\big)^r
\ .
\label{frm:stima-H(r)menoH(r-1)}
\end{equation}
Recalling that the sup-norm is bounded by the weighted Fourier one,
this proves that $\{H^{(r)}\}_{r\ge 0}$ is a Cauchy sequence of
analytic Hamiltonians admitting a limit $H^{(\infty)}$, since the
r.h.s. of the estimate above tends to zero for $r\to\infty$.
Moreover,
formul{\ae}~\eqref{frm:H(r)menoH(r-1)-espansione}--\eqref{frm:stima-H(r)menoH(r-1)}
imply that also $\{\vet\omega^{(r)}\}_{r\ge 0}\,$, $\{\vet \Omega^{(r)}\}_{r\ge
0}$ and $\{f_{\ell}^{(r,s)}\}_{r\ge 0}$ for $\ell\ge 3\,,\ s\ge 0$ are
Cauchy sequences, which in turn implies that
$H^{(\infty)}=\lim_{r\to\infty}H^{(r)}$ has the
form~\eqref{frm:H(infty)-espansione}, with
$f_{\ell}^{(\infty,s)}\in\Pset_{\ell}^{sK}$ bounded as
in~\eqref{frm:f_l^(infty,s)}, in view of
inequality~\eqref{frm:bound-fl}.

It remains to prove that the canonical transformation
$\Upphi_{\vet \omega^{(0)}}^{(\infty)}:\Dscr_{1/2(\rho,\sigma,R)}\to
\Dscr_{3/4(\rho,\sigma,R)}$ is analytic.  The proof is based on
standard arguments in the Lie series theory, that we recall here,
referring to subsection~4.3 of~\cite{Gio-Loc-1997} for more details.
Let us denote by $\phi^{(r)}$ the canonical change of coordinates
induced by the $r$-th step, i.e.,
\begin{equation}
\phi^{(r)}(\vet p,\vet q,\vet z,i\bar{\vet  z})=
\exp\left( \epsilon^r\Lie{\chi_0^{(r)}} \right) \circ
\exp\left( \epsilon^r\Lie{\chi_1^{(r)}} \right)\circ
\exp\left( \epsilon^r\Lie{\chi_2^{(r)}} \right)
(\vet p,\vet q,\vet z,\imunit\bar {\vet z})\ .
\label{frm:def-phi(r)}
\end{equation}
Using inequalities in~\eqref{frm:generatrici-lemmone},
one can easily verify that the following inequalities hold true:
\begin{equation}
\vcenter{\openup1\jot 
\halign{
$\displaystyle\hfil#$&$\displaystyle{}#\hfil$&$\displaystyle#\hfil$\cr
\max_{1\le j\le n_1}
\Big\|\exp\left( \epsilon^r\Lie{\chi_0^{(r)}} \right) p_j-p_j\Big\|_{3/4 - d_{r-1} - \delta_r}
& < &\delta_r \rho\, ,
\cr
\max_{1\le j\le n_1}
\Big\|\exp\left( \epsilon^r\Lie{\chi_1^{(r)}} \right) p_j-p_j\Big\|_{3/4 - d_{r-1} - 2\delta_r}
& < &\delta_r \rho\, ,
\cr
\max_{1\le j\le n2}
\Big\|\exp\left( \epsilon^r\Lie{\chi_1^{(r)}} \right) z_j-z_j\Big\|_{3/4 - d_{r-1} - 2\delta_r}
&< &\delta_r R\, ,
\cr
\max_{1\le j\le n_1}
\Big\|\exp\left( \epsilon^r\Lie{\chi_2^{(r)}} \right) p_j-p_j\Big\|_{3/4 - d_r}
& < &\delta_r \rho\, , 
\cr
\max_{1\le j\le n_1}
\Big\|\exp\left( \epsilon^r\Lie{\chi_2^{(r)}} \right) q_j-q_j\Big\|_{3/4 - d_r}
& < &\delta_r \sigma\, , 
\cr
\max_{1\le j\le n_2}
\Big\|\exp\left( \epsilon^r\Lie{\chi_2^{(r)}} \right) z_j-z_j\Big\|_{3/4 - d_r}
&<&\delta_r R\, .
\cr
}}
\label{frm:scartamento}
\end{equation}
For example, let us consider the generating function $\chi_0^{(r)}$, that depends only on $\vet q$. 
Then,
\begin{align*}
\Big\|\exp\left( \epsilon^r\Lie{\chi_0^{(r)}} \right) p_j-p_j\Big\|_{3/4- d_{r-1}- \delta_r} 
& \le 
\sum_{s\ge1} \frac{\epsilon^{sr}}{s!} \Big\| \Lie{\chi_0^{(r)}}^{s-1} \frac{\partial \chi_0^{(r)}}{\partial q_j} \Big\| 
\\
& \le \rho \delta_r\sum_{s\ge1} \frac{\epsilon^{sr}}{s} \left(M^{3r} \Vscr\left(\Gscr_0^{(r)}\right) \nu_{r-1,r}\right)^s .
\end{align*}
Thanks to the condition on $\epsilon$, the sum converges and we obtain the corresponding estimate in~\eqref{frm:scartamento}.
Estimates for the other Lie series appearing in~\eqref{frm:scartamento} can be deduced analogously using the inequalities in~\eqref{frm:generatrici-lemmone} and in~\eqref{frm:stimalie}.
Therefore, one has
$\phi^{(r)}(\Dscr_{(1/2+d_{r-1})(\rho,\sigma,R)})\subset
\Dscr_{(1/2+d_{r})(\rho,\sigma,R)}$ and defining
$\Phi^{(r)}=\phi^{(1)}\circ\ldots\circ\phi^{(r)}$ one has
$\Phi^{(r)}(\Dscr_{1/2(\rho,\sigma,R)})\subset\Dscr_{3/4(\rho,\sigma,R)}\,$.
By repeatedly using the so-called exchange theorem for Lie
series, one immediately obtains that
$H^{(r)}=H^{(0)}\circ\Phi^{(r)}$. By using
estimate~\eqref{frm:scartamento}, we can prove that the
canonical transformation
$\Phi_{\vet \omega^{(0)}}^{(\infty)}=\lim_{r\to\infty}\Phi^{(r)}$ is well
defined in $\Dscr_{1/2(\rho,\sigma,R)}\,$. Finally, we get
$H^{(0)}\circ\Phi_{\vet \omega^{(0)}}^{(\infty)}=\lim_{r\to\infty}H^{(0)}\circ\Phi^{(r)}=
\lim_{r\to\infty}H^{(r)}=H^{(\infty)}\,$.

Actually, with some additional effort, we could prove that
$\Phi_{\vet \omega^{(0)}}^{(\infty)}$ differs from the identity just for
terms of order $\Oscr(\epsilon)\,$. As a final comment, we adopt the
symbol $\Phi_{\vet \omega^{(0)}}^{(\infty)}\,$ in order to emphasize the
parametric dependence of that canonical transformation on the initial
frequency $\vet \omega^{(0)}\,$, which is relevant in the next section.

\subsection{Geometric part: measure of the resonant regions}\label{sec:measure}
The aim of this section is to show that there exists a set of
frequencies $\vet \omega^{(0)}$ of relatively big measure, for which our
procedure converges.  To this end we consider the sequence
$\big\{(\vet \omega^{(r)}\,,\,\vet \Omega^{(r)})\big\}_{r\ge 0}$ as function of
the parameter $\vet \omega^{(0)}$ recursively defined
by~\eqref{frm:chgfreq.toro} and~\eqref{frm:chgfreq.trasv}.  We start
from the compact set $\Wscr^{(0)}\subset\reali^{n_1}$ and its
complex extension $\Wscr^{(0)}_{h_0}$ satisfying the assumptions in
lemma~\ref{lem:Ham-expansion}. Then, we construct a sequence of
complex extended domains
$\Wscr^{(0)}_{h_0}\supseteq\Wscr^{(1)}_{h_1}\supseteq\Wscr^{(2)}_{h_2}
\supseteq\ldots\,$, where $\{h_{r}\}_{r\ge 0}$ is a positive
non-increasing sequence of real numbers, with the following
properties: $\vet \omega^{(r)}(\vet \omega^{(0)})$ is an analytic function
admitting an inverse $\phi^{(r)}$ mapping the domain
$\Wscr^{(r)}_{h_{r}}\,$ to $\Wscr^{(0)}_{h_{0}}\,$.  We must prove the
sequence of functions $\phi^{(r)}$ converges to $\phi^{(\infty)}$
mapping $\Wscr^{(\infty)}_{h_{\infty}}\,$ to $\Wscr^{(0)}_{h_{0}}\,$,
the image having large relative measure.

Let us describe the process in some more detail. For $r\ge 1$ and some
fixed positive values of the parameters $\gamma\, ,\tau\in\reali$ and
$K\in\naturali\,$, we define the sequence of real domains
$\{\Wscr^{(r)}\}_{r\ge 0}\,$ as follows: at each step $r$, we remove
from $\Wscr^{(r-1)}$ all the resonant regions related to the new small
divisors appearing in the formal algorithm described in section~\ref{sec:algorithm}. Therefore,
\begin{equation}
\Wscr^{(r)}=\Wscr^{(r-1)}\,\backslash\,\Rscr^{(r)}
\ ,
\qquad{\rm with}\quad
\Rscr^{(r)} = \bigcup_{{rK<|k|\leq (r+1)K}\atop{|\ell|\leq 2}}\Rscr^{(r)}_{k,\ell}\ ,
\label{frm:WWWr}
\end{equation}
where each so called resonant strip is defined as follows:
\begin{equation}
\Rscr_{\vet k,\vet \ell}^{(r)}=\left\{\vet \omega\in\Wscr^{(r-1)}:
\left| \vet k\cdot\vet \omega+
\vet \ell\,\cdot\vet \Omega^{(r)}\circ\phi^{(r)}(\vet \omega)\right|\lt
\frac{2\gamma}{\big((r+1)K)^{\tau}}\right\}\ .
\label{frm:strisciarisonantekl}
\end{equation}
The complex extensions $\Wscr^{(r)}_{h_{r}}\,$ are technically
necessary in order to work with Cauchy estimates, but it is crucial
that the measure is evaluated on the real domain.
It will be also convenient to introduce the functions
$\delta\vet \omega^{(r)}$ and $\Delta\vet \Omega^{(r)}$ defined as
\begin{equation}
\delta\vet \omega^{(r)}=\vet \omega^{(r)}\circ\phi^{(r-1)}-{\rm Id}\ ,
\qquad
\Delta\vet \Omega^{(r)}=\vet \Omega^{(r)}\circ\phi^{(r-1)}-\vet \Omega^{(r-1)}\circ\phi^{(r-1)}\ .
\label{frm:delta-Delta-frequenze}
\end{equation}

The following part is an adaptation of the
approach described by P\"oschel in~\cite{Poschel-1989}. First, we report here lemma~D.1 in~\cite{Poschel-1989}, since it will be used in the proof of the proposition.
\begin{lemma}\label{lem:D1-poschel}
Let $f$ be a real analytic function from $\Wscr_h$ to $\complessi^n$, with $\Wscr_h$ the complex estension of a domain $\Wscr\subset \reali^n$.
If 
\begin{equation}
|f - {\rm Id}|_{\infty; \Wscr_h}\le \delta \le  \frac{h}{4} \,,
\end{equation}
then $f$ admits a real analytic inverse $\phi$ on $\Wscr_{\frac h 4}$, that satisfies the following inequalities
\begin{equation}
|\phi- {\rm Id}|_{\infty, \Wscr_{\frac h 4} }\le \delta \quad and \quad |\partial_\omega (\phi- {\rm Id})|_{\infty, \Wscr_{\frac h 4}}  \le  \frac{4 \delta }{h}\,,
\end{equation}
where the norm $\big|\cdot\big|_{\infty;\Wscr^{(0)}_{\frac{h_0}{4}}}$ is defined in a way analogous
to~\eqref{frm:cost-Lipschitz}.
\end{lemma}

\begin{proposition}\label{pro:geometric-prop}
Let us consider 
the family~\eqref{frm:H(0)} of Hamiltonians $H^{(0)}$ parametrized
by the $n_1$-dimensional frequency vector $\vet \omega^{(0)}$. Assume that
there exist positive parameters $\sigma\,$, $\gamma\,$, $\tau\,$, $\bar b\,$, $J_0$,
a positive integer $K$ and a compact set $\Wscr\subset\reali^{n_1}$
such that the function $\vet \Omega^{(0)}:\Wscr_{h_0}\to\complessi^{n_2}$ is
analytic and satisfies the properties~(a')--(b') of
lemma~\ref{lem:Ham-expansion}.  Define the sequence
$\{h_{r}\}_{r\ge 0}$ of radii of the complex extensions as
\begin{equation}
h_0 = \min\left\{ \eta, \frac{1}{e\left(J_0+\frac{1}{\sigma}\right)}\right\}\cdot \frac{\gamma}{8 K^\tau}
\qquad
\hbox{and}
\qquad
h_r = \frac{h_{r-1}}{2^{\tau+2}}
\ \ \hbox{for } r\ge 1,
\label{frm:hr}
\end{equation}
with $\eta=\min\{1/K\,,\,\sigma\}\,$.

\noindent
Considering the sequence of Hamiltonians $\{H^{(r)}\}_{r\ge 0}\,$,
formally defined by the algorithm in section~\ref{sec:algorithm},
let us assume that the functions
$\vet \omega^{(1)},\,\vet \Omega^{(1)},\,\ldots\,,\,\vet \omega^{(r)},\,\vet \Omega^{(r)}$
satisfy the following hypotheses up to a fixed normalization step $r\ge 0$:

\item{(g')} the function $\vet \omega^{(s)}(\vet \omega^{(0)})$ has analytic inverse
$\phi^{(s)}$ on $\Wscr^{(s)}_{h_{s}}\,$, for $0\le s \le r-1$, where
$\phi^{(0)}={\rm Id}$ and the domains are recursively defined as
extensions of those given by
formul{\ae}~\eqref{frm:WWWr}--\eqref{frm:strisciarisonantekl},
starting from $\Wscr^{(0)}$;

\item{(h')} both
$\vet \omega^{(s)}\circ\phi^{(s-1)}:\,\Wscr^{(s-1)}_{h_{s-1}}\to\complessi^{n_1}$
and
$\vet \Omega^{(s)}\circ\phi^{(s-1)}:\,\Wscr^{(s-1)}_{h_{s-1}}\to\complessi^{n_2}$
are analytic functions, for $1\le s \le r$ ;

\item{(i')} for $1\le s \le r\,$, there exist positive
parameters $\epsilon\,$, $\sigma$ and $\Ascr\ge 1$ satisfying

\begin{equation}
\max_{1\le j\le n_1}\sup_{\vet\omega\in \Wscr^{(s-1)}_{h_{s-1}}}
\big|\delta\omega^{(s)}_j(\vet\omega)\big|\le
\gamma\sigma\big(\epsilon\Ascr\big)^s\ ,
\quad
\max_{1\le j\le n_2}\sup_{\vet \omega\in \Wscr^{(s-1)}_{h_{s-1}}}
\big|\Delta\Omega^{(s)}_j(\vet\omega)\big|\le
\gamma(\epsilon\Ascr)^s\ ,
\label{frm:stima-variazione-frequenze-per-parte-geom}
\end{equation}
\item{} where $\delta\vet \omega^{(s)}$ and $\Delta\vet \Omega^{(s)}$ are defined
as in~\eqref{frm:delta-Delta-frequenze}, while $\Ascr$ is given in~\eqref{frm:soglia-analitica};

\item{(j')} the parameter $\epsilon$ is smaller than the ``geometric
threshold value''

\begin{equation}
\epsilon^{*}_{{\rm ge}} = \frac{1}{2^{\tau+3}\Ascr}\,
\min\left\{1\,,\,\frac{h_0}{ \sigma \gamma
 }\,\right\}.
\label{frm:soglia-geometrica}
\end{equation}
Then, the function $\vet \omega^{(r)}(\vet \omega^{(0)})$ admits an analytic
inverse
$\phi^{(r)}:\Wscr^{(r)}_{h_{r}} \to \Wscr^{(0)}_{h_{0}}$
on its domain of definition and satisfies the inclusion relation
$\phi^{(r)}\big(\Wscr^{(r)}_{h_{r}}\big)
\subset\phi^{(r-1)}\big(\Wscr^{(r-1)}_{h_{r-1}}\big)\,$. Moreover, the following
non-resonance inequalities hold true:

\begin{equation}
\vcenter{\openup1\jot\halign{
 \hbox {\hfil $\displaystyle {#}$}
&\hbox {\hfil $\displaystyle {#}$\hfil}
&\hbox {$\displaystyle {#}$\hfil}\cr
\min_{\scriptstyle{\vet k\in\interi^{n_1}\,,\,0<|\vet k|\leq (r+1)K}
\atop\scriptstyle{\vet \ell\in\interi^{n_2}\,,\,0\leq |\vet \ell|\leq2}}
\,\inf_{\vet \omega\in\Wscr^{(r)}_{h_{r}}}
\left| \vet k\cdot\vet \omega
+ \vet \ell\,\cdot\vet \Omega^{(r)}\big(\phi^{(r)}(\vet \omega)\big)\right|
&\ge &
\frac{\gamma}{\big((r+1)K\big)^{\tau}}\ ,
\cr
\min_{\vet \ell\in\interi^{n_2}\,,\,0<|\vet \ell|\leq2}\,\inf_{\vet \omega\in\Wscr^{(r)}_{h_{r}}}
\big| \scalprod{\vet \ell}{\vet \Omega^{(r)}\big(\phi^{(r)}(\vet \omega)\big)} \big|
&\ge &
\gamma\ .
\cr
}}
\label{frm:nonres-cond-diof}
\end{equation}
Finally, the Lipschitz constants related to the Jacobians of the
functions $\phi^{(r)}$ and $\vet \Omega^{(r)}\circ\phi^{(r)}$ are uniformly
bounded as
\begin{equation}
\vcenter{\openup1\jot\halign{
 \hbox {\hfil $\displaystyle {#}$}
&\hbox {\hfil $\displaystyle {#}$\hfil}
&\hbox {$\displaystyle {#}$\hfil}\cr
\left|\frac{\partial\big(\phi^{(r)}-{\rm Id}\big)}{\partial\vet \omega}
\right|_{\infty;\Wscr^{(r)}_{h_r}}&\leq& \frac{16(\epsilon \Ascr)}{\Bscr} \,,
\cr
\left|\frac{\partial\big(\vet \Omega^{(r)}\circ\phi^{(r)}\big)}{\partial\vet \omega}
\right|_{\infty;\Wscr^{(r)}_{h_r}}&\le& e \left(J_0 + \frac{1}{\sigma}\right)\,,
\cr
}}
\label{frm:diseq:Lipschitz-definitive}
\end{equation}
where $\Bscr$ is defined as
\begin{equation}
\Bscr:= \min\left\{ \frac{h_0}{\sigma\gamma}\,, \frac{1}{2^\tau}\right\} \,.
\label{def:bstorto}
\end{equation}
\end{proposition}
\begin{proof}
The proof proceeds by induction. By
hypotheses, $\vet \Omega^{(0)}(\vet \omega)$ is an analytic function on the
complex extended domain $\Wscr^{(0)}_{h_0}$ and its Jacobian is
uniformly bounded in $\Wscr^{(0)}_{h_0}\,$, namely
$\big|\partial\vet \Omega^{(0)}/\partial\vet \omega\big|_{\infty;\Wscr^{(0)}_{h_0}}\leq
J_0\,$, where $\big|\cdot\big|_{\infty;\Wscr^{(0)}_{h_0}}$ is defined
in~\eqref{frm:cost-Lipschitz}. Thus, starting from the
inequality at point~(a') of
lemma~\ref{lem:Ham-expansion}, if 
\begin{equation}
h_0 \leq \frac{1}{\max\{ \frac{K}{2}, \frac{1}{\sigma}\} + J_0}
\frac{\gamma}{4K^{\tau}} \ ,
\label{diseq:h0-piccolo-quanto-basta}
\end{equation}
then the non-resonance conditions~\eqref{frm:nonres_a}
and~\eqref{frm:nonres_b} are satisfied in the complexified domain
$\Wscr^{(0)}_{h_0}$ for $0<|\vet k|\leq K\,$, $|\vet \ell|\leq 2$.  More precisely, if $\vet \omega \in \Wscr^{(0)}$ and $|\bar{\vet \omega}|< h_0$, we have
\begin{equation}
\vcenter{\openup1\jot 
\halign{
$\displaystyle\hfil#$&$\displaystyle{}#\hfil$&$\displaystyle#\hfil$\cr
\left| \vet k\cdot(\vet \omega + \bar{\vet \omega}) + \vet \ell\,\cdot\vet \Omega^{(0)}(\vet \omega + \bar{\vet \omega})\right|&\geq
 \left|  \vet k\cdot\vet \omega + \vet \ell\,\cdot\vet \Omega^{(0)}(\vet \omega)\right|
-K h_0 - 2h_0 J_0 
\cr
&\ge
\frac{(2-1/2)\gamma}{K^\tau}
\cr
}}
\label{frm:nonres-cond-diof-step1}
\end{equation}
and, for $0< |\vet \ell| \le 2\,$,
\begin{equation}
\left| \vet \ell \cdot \vet \Omega^{(0)}(\vet \omega + \bar{\vet \omega})\right|\geq
 \left| \vet \ell \cdot \vet \Omega^{(0)}(\vet \omega)\right|
- 2h_0 J_0 \ge
\left(2-\frac{1}{2}\right) \gamma\ .
\label{frm:nonres-cond-diof-step1trasv}
\end{equation}
For $r=0$, the inequalities in~\eqref{frm:nonres-cond-diof} immediately
follow from the previous ones, recalling that $\phi^{(0)}={\rm
Id}$.  

The first change of the frequencies might occur at the end of the
first perturbation step and the transformed fast frequencies read
$$
\vet \omega^{(1)}(\vet \omega^{(0)}) = \vet \omega^{(0)}+
\delta\vet \omega^{(1)}(\vet \omega^{(0)})
= \left( {\rm Id} +\delta\vet \omega^{(1)} \right)(\vet \omega^{(0)})\,,
$$
\noindent
where $\max_{1\le j\le n_1}\sup_{\omega\in \Wscr^{(0)}_{h_{0}}}
\big|\delta\omega^{(1)}_j(\omega)\big|\le
\gamma\sigma\big(\epsilon\Ascr\big)$ in view of the assumption~(i')
in proposition~\ref{pro:geometric-prop}.  If
\begin{equation}
\mu_0 = \frac{4\gamma\sigma(\epsilon\Ascr)}{h_0} \le 1\ ,
\label{def:mu1}
\end{equation}
in view of lemma~\ref{lem:D1-poschel}, the function $({\rm
Id}+\delta\vet \omega^{(0)})$ admits an analytic inverse
$\phi^{(1)}:\Wscr^{(0)}_{\frac{h_0}{4}} \to \Wscr^{(0)}_{\frac{h_0}{2}}$ and, in the
domain $\Wscr^{(0)}_{\frac{h_0}{4}}\,$, the following estimates hold true:
\begin{equation}
\max_{{\scriptstyle{1\le j\le n_2}}}\,\sup_{\vet \omega\in \Wscr^{(0)}_{\frac{h_0}{4}}}\,
\left|\phi^{(1)}_j(\vet \omega)-\omega_j\right|
\leq \gamma \sigma(\epsilon\Ascr)\ ,
\qquad
\left|\frac{\partial\left(\phi^{(1)}-{\rm Id}\right)}{\partial\vet \omega}
\right|_{\infty;\Wscr^{(0)}_{\frac{h_0}{4}}}
\leq \mu_0\ ,
\label{diseq:stime-phi2}
\end{equation}
where the norm $\big|\cdot\big|_{\infty;\Wscr^{(0)}_{\frac{h_0}{4}}}$ on the
Jacobian of the function $\phi^{(1)}-{\rm
Id}:\,\Wscr_{\frac{h_0}{4}}\to\complessi^{n_1}$ is defined in a way analogous
to~\eqref{frm:cost-Lipschitz}.  The growth of the
Lipschitz constants for the function
$\phi^{(1)}-{\rm Id}$ is controlled by setting
\begin{equation}
\bar J_0=0\ ,
\qquad
\bar J_1=\mu_0\ ,
\label{def:Jbarra-da-0-a-2}
\end{equation}
i.e., we have $\big|{\partial\left(\phi^{(1)}-{\rm
Id}\right)}/{\partial\vet \omega} \big|_{\infty;\Wscr^{(0)}_{\frac{h_0}{4}}}\leq \bar J_1\,$.
Recall now that the step $r=1$ includes also the preparation of the
next step $r=2$, namely cutting out the resonant regions
$$
\Rscr_{\vet k,\vet \ell}^{(1)}=\left\{\vet \omega\in\Wscr^{(0)}:
\big| \vet k\cdot\vet \omega+\vet \ell\,\cdot\vet \Omega^{(1)}\circ\phi^{(1)}(\vet \omega)\big|
\leq 2\gamma/(2K)^{\tau}\right\}\ .
$$
Thus we need an upper bound on both the sup-norm of
$\vet \Omega^{(1)}\circ\phi^{(1)}$ and the Lipschitz constant of its
Jacobian.  The new transverse frequencies can be written
as
$$
\vet \Omega^{(1)}(\vet \omega^{(0)})=
\vet \Omega^{(0)}(\vet \omega^{(0)})+
\Delta\vet \Omega^{(1)}(\vet \omega^{(0)})\ ,
$$
where we have $\max_{1\le j\le n_1}\sup_{\vet \omega\in \Wscr^{(0)}_{h_0}}
\big|\Delta\Omega^{(1)}_j(\vet \omega)\big|\leq \gamma(\epsilon\Ascr)$,
in view of the hypothesis~(i')
in proposition~\ref{pro:geometric-prop}.  Moreover,
in the domain $\Wscr^{(0)}_{\frac{h_0}{4}}$ the new transverse
frequencies are functions of the transformed fast frequencies, namely
\begin{equation}
\vet \Omega^{(1)}\big(\phi^{(1)}(\vet \omega^{(1)})\big) =
\vet \Omega^{(0)}\big(\phi^{(1)}(\vet \omega^{(1)})\big) +
\Delta\vet \Omega^{(1)}\big(\phi^{(1)}(\vet \omega^{(1)})\big)\ .
\label{frm:Omegoni-passo2-da-passo1}
\end{equation}
Using this formula we can bound the Jacobian of the function
$\vet \Omega^{(1)}\circ\phi^{(1)}$.  To this end, we need the preliminary
estimate 
$$
\vcenter{\openup1\jot\halign{
 \hbox {\hfil $\displaystyle {#}$}
&\hbox {\hfil $\displaystyle {#}$\hfil}
&\hbox {$\displaystyle {#}$\hfil}\cr
\left|\frac{\partial\big(\vet \Omega^{(1)}\circ\phi^{(1)}-\vet \Omega^{(0)}\big)}
        {\partial\vet \omega}\right|_{\infty;\Wscr^{(0)}_{\frac{h_0}{4}}}
&\leq\left|\frac{\partial\vet \Omega^{(0)}}{\partial\vet \omega}
\right|_{\infty;\Wscr^{(0)}_{\frac{h_0}{2}}}
\left|\frac{\partial\big(\phi^{(1)}-{\rm Id}\big)}{\partial\vet \omega}
\right|_{\infty;\Wscr^{(0)}_{\frac{h_0}{4}}}
\cr
&\phantom{\leq}+
\left|\frac{\partial\Delta\vet \Omega^{(1)}}{\partial\vet \omega}
\right|_{\infty;\Wscr^{(0)}_{\frac{h_0}{2}}}
\,\left|\frac{\partial\phi^{(1)}}{\partial\vet \omega}\right|_{\infty;\Wscr^{(0)}_{\frac{h_0}{4}}}
\cr
&\le J_0 \mu_0 + \frac{2\gamma \epsilon \Ascr}{h_0 }(1+ \mu_0) \le J_0 \mu_0 + \frac{\mu_0}{2\sigma}(1+ \mu_0)\ ,
\cr
}}
$$
where we used formul{\ae}~\eqref{def:mu1}--\eqref{diseq:stime-phi2},
the estimate
$\big|{\partial\vet \Omega^{(0)}}/{\partial\vet \omega}\big|_{\infty;\Wscr^{(0)}_{\frac{h_0}{2}}}\le
J_0$ and Cauchy inequality to estimate
$\big|{\partial\Delta\vet \Omega^{(1)}}/{\partial\vet \omega}\big|_{\infty;\Wscr^{(0)}_{\frac{h_0}{2}}}\,$.
Thus, we can ensure that the Jacobian of the transformed
transverse frequencies is bounded as
\begin{equation}
\left|\frac{\partial\big(\vet \Omega^{(1)}\circ\phi^{(1)}\big)}
{\partial\vet \omega}\right|_{\infty;\Wscr^{(0)}_{\frac{h_0}{4}}}
\leq \left(J_0  + \frac{\mu_0}{2\sigma}\right)(1+\mu_0)=\mathrel{\mathop:} J_1\ .
\label{eq:J1}
\end{equation}
Using again formula~\eqref{frm:Omegoni-passo2-da-passo1}, we can
bound the deterioration of the non-resonance conditions involving the
transverse frequencies. In fact, for $\vet \ell\in\interi^{n_2}$, $|\vet \ell|\le
2\,$, we have 
\begin{equation}
\vcenter{\openup1\jot\halign{
&\hbox {$\displaystyle {#}$\hfil}\cr
&\sup_{\vet \omega\in \Wscr^{(0)}_{\frac{h_0}{4}}}
\left|\vet \ell\cdot\left[\vet \Omega^{(1)}\big(\phi^{(1)}(\vet \omega)\big)-
\vet \Omega^{(0)}(\vet \omega)\right]\right|
\cr
&\qquad\leq 2\max_j
\sup_{\vet \omega\in \Wscr^{(0)}_{\frac{h_0}{4}}}\left| \Omega^{(0)}_j\big(\phi^{(1)}(\vet \omega)\big)-
\Omega^{(0)}_j(\vet \omega)\right|
+2 \max_j \sup_{\vet \omega\in \Wscr^{(0)}_{\frac{h_0}{4}}}
\left|\Delta\Omega^{(1)}_j\big(\phi^{(1)}(\vet \omega)\big)\right|
\cr
&\qquad\leq 2\left|\frac{\partial\vet \Omega^{(0)}}{\partial\vet \omega}
\right|_{\infty;\Wscr^{(0)}_{\frac{h_0}{2}}}
\,\max_j\,\sup_{\vet \omega\in \Wscr^{(0)}_{\frac{h_0}{4}}}\left|\phi^{(1)}_j(\vet \omega)-\omega_j\right|
+2\max_j\,\sup_{\vet \omega\in \Wscr^{(0)}_{\frac{h_0}{2}}}
\left|\Delta\Omega^{(1)}_j(\vet \omega)\right|
\cr
&\qquad\leq 2 J_0 \gamma \sigma(\epsilon\Ascr)+2\gamma(\epsilon\Ascr)
= \frac {\mu_0 h_0}{2}\left(\frac{1}{\sigma}+J_0\right)\ .
\cr
}}
\label{diseq:stima-var-freq-sec}
\end{equation}
Thus, for $\vet \omega\in\Wscr^{(0)}_{\frac{h_0}{4}}\,$, $0<|\vet k|\leq K$ and
$|\vet \ell|\leq 2\,$, we obtain the non-resonance estimate
\begin{equation}
\vcenter{\openup1\jot\halign{
 \hbox {\hfil $\displaystyle {#}$}
&\hbox {$\displaystyle {#}$\hfil}\cr
\left| \vet k\cdot \vet \omega
+ \vet \ell\,\cdot\vet \Omega^{(1)}\big(\phi^{(1)}(\vet \omega)\big)\right|
&\ge
\frac{(2-1/2)\gamma}{K^\tau}
-\sup_{\vet \omega\in \Wscr^{(0)}_{\frac{h_0}{4}}}
\left|\vet \ell\cdot\left[\vet \Omega^{(1)}\big(\phi^{(1)}(\vet \omega)\big)-
\vet \Omega^{(0)}(\vet \omega)\right]\right|
\cr
&\ge \frac{(2-1/2)\gamma}{K^{\tau}}-
\mu_0 \left(\frac{1}{\sigma}+J_0\right)\frac{h_0}{2}
\ge
\frac{(2-1/2-1/4)\gamma}{K^{\tau}}\ ,
\cr
}}
\label{frm:nonres-cond-diof-step2-vecchio-blocco}
\end{equation}
where we started from estimate~\eqref{frm:nonres-cond-diof-step1}
(holding true on all the complex domain $\Wscr^{(0)}_{\frac{h_0}{4}}$), we
used~\eqref{diseq:stima-var-freq-sec}, \eqref{def:mu1}, the
inequality for $h_0$ in~\eqref{diseq:h0-piccolo-quanto-basta}.  Moreover, 
using the estimates in~\eqref{frm:nonres-cond-diof-step1trasv},~\eqref{diseq:stima-var-freq-sec} and the inequality for $h_0$ in~\eqref{diseq:h0-piccolo-quanto-basta}, one can easily obtain the lower bound
\begin{equation}
\left| \vet \ell \cdot \vet \Omega^{(1)}\big(\phi^{(1)}(\vet \omega)\big)\right|
\geq
\left(2-\frac{1}{2}-\frac{1}{4}\right) \gamma\,,
\label{nonres-cond-trasversa-step2}
\end{equation}
uniformly with respect to $\vet \omega\in\Wscr^{(0)}_{\frac{h_0}{4}}\,$, when
$0<|\vet \ell|\le 2$. 
It is now time to consider the new subset of resonant
regions $\Rscr_{\vet k,\vet \ell}^{(1)}$ for $K<|\vet k|\leq 2K\,$, $|\vet \ell|\leq 2\,$.
First, we remove them from the domain, by defining $\Wscr^{(1)}$
according to~\eqref{frm:WWWr}--\eqref{frm:strisciarisonantekl}.
Having required that the new radius of the complex extension is
so small that 
\begin{equation}
h_1 \le \min\left\{\frac{h_0}{4}\,,
\ \frac{1}{\max\{K\,,\,1/\sigma\}+J_1}\,
\frac{\gamma}{8(2K)^{\tau}}\right\}\ ,
\label{def:h1}
\end{equation}
for $K<|\vet k|\leq 2K\,$, $|\vet \ell|\leq 2\,$ and $\vet \omega \in \Wscr^{(1)}$, $|\bar{\vet \omega}| < h_1$, we have
\begin{equation}
\vcenter{\openup1\jot\halign{
 \hbox {\hfil $\displaystyle {#}$}
&\hbox {$\displaystyle {#}$\hfil}\cr
\left| \vet k\cdot (\vet \omega +\bar{\vet \omega})
+ \vet \ell\,\cdot\vet \Omega^{(1)}\big(\phi^{(1)}(\vet \omega +\bar{\vet \omega})\big)\right|
&\ge
\left| \vet k\cdot\vet \omega
+ \vet \ell\,\cdot\vet \Omega^{(1)}\big(\phi^{(1)}(\vet \omega)\big)\right|
\cr
&\phantom{\ge}-2Kh_{1}-2J_1h_{1}
\ge\frac{(2-1/2-1/4)\gamma}{(2K)^{\tau}}\ .
\cr
}}
\label{frm:nonres-cond-diof-step2-nuovo-blocco}
\end{equation}
\noindent
Here we used the
definitions~\eqref{frm:WWWr}--\eqref{frm:strisciarisonantekl} for
the real set $\Wscr^{(1)}$ and subtracted the contribution due to the
complex extension; moreover, we also used
formul{\ae}~\eqref{eq:J1},~\eqref{def:h1}, (recall that
$\Wscr^{(1)}_{h_{1}}\subseteq\Wscr^{(0)}_{\frac{h_0}{4}}\,$, in view
of~\eqref{frm:WWWr} and~\eqref{def:h1}).  The inequalities
in~\eqref{frm:nonres-cond-diof} are justified in view
of~\eqref{frm:nonres-cond-diof-step2-vecchio-blocco},
\eqref{nonres-cond-trasversa-step2}
and~\eqref{frm:nonres-cond-diof-step2-nuovo-blocco}.  This concludes the
proof of the statement for $r=1\,$.

Iterating the procedure for a generic step $r>1$ is now
straightforward, provided the sequence of restrictions of the frequency
domain is suitably selected.  Let us describe such a scheme of estimates in detail. 

We restart from the relation
$$
\vet \omega^{(r)}(\vet \omega^{(0)}) =
\vet \omega^{(r)}\circ\phi^{(r-1)}\circ\vet \omega^{(r-1)}(\vet \omega^{(0)})=
\left( {\rm Id} + \delta\vet \omega^{(r)}\right)\big(\vet \omega^{(r-1)}(\vet \omega^{(0)})\big)\ ,
$$
where $\max_{1\le j\le n_1}\sup_{\vet \omega\in \Wscr^{(r-1)}_{h_{r-1}}}
\big\{\big|\delta\omega^{(r)}_j(\vet \omega)\big|\big\}\le
\gamma\sigma\big(\epsilon\Ascr\big)^r$ in view of assumption~\eqref{frm:stima-variazione-frequenze-per-parte-geom}.  Thus, using again lemma~\ref{lem:D1-poschel}, if
\begin{equation}
\mu_{r-1} = \frac{4\gamma\sigma(\epsilon\Ascr)^r}{h_{r-1}} \le 1\,,
\label{def:mu_r-1}
\end{equation} 
we obtain $\vet \omega^{(r-1)}(\vet \omega^{(0)})$ from
$\vet \omega^{(r)}(\vet \omega^{(0)})\,$ via the function
$\flusso^{(r)}:\Wscr^{(r-1)}_{\frac{h_{r-1}}{4}} \to \Wscr^{(r-1)}_{\frac{h_{r-1}}{2}}$, because 
\begin{equation}
|\vet \omega^{(r)}\circ \phi^{(r-1)}- \rm{Id}| \le \frac{h_{r-1}}{4}\,;
\end{equation}
moreover, the inverse function $\flusso^{(r)}$ satisfies
\begin{equation}
|\flusso^{(r)}- {\rm Id}| \le \gamma \sigma  (\epsilon \Ascr)^r\, \qquad {\rm and }\qquad |\partial(\flusso^{(r)}- {\rm Id})/\partial \vet \omega| \le \mu_{r-1}\,.
\end{equation}
The function $\phi^{(r)}$, namely the inverse of
$\vet \omega^{(r)}(\vet \omega^{(0)})$, is obtained by composition, i.e.,
$\phi^{(r)}=\phi^{(r-1)}\circ\flusso^{(r)}=
\flusso^{(1)}\circ\ldots\circ\flusso^{(r)}$ and, by construction, we have
$\phi^{(r)}\big(\Wscr^{(r-1)}_{\frac{h_{r-1}}{4}}\big)
\subset\phi^{(r-1)}\big(\Wscr^{(r-1)}_{h_{r-1}}\big)\,$ (remark that $\phi^{(1)} = \flusso^{(1)}$). Replacing $\phi^{(1)}$ with $\flusso^{(r)}$,
formul{\ae}~\eqref{def:mu1}--\eqref{diseq:stima-var-freq-sec}
can be suitably adapted to the case $r>1\,$, so as to prove the
inequalities corresponding
to~\eqref{frm:nonres-cond-diof-step2-vecchio-blocco},
\eqref{nonres-cond-trasversa-step2}
and~\eqref{frm:nonres-cond-diof-step2-nuovo-blocco}.  
First, starting from the relation $\phi^{(r)} = \phi^{(r-1)}\circ \flusso^{(r)}$, the following inequalities hold true:
$$
\vcenter{\openup1\jot\halign{
 \hbox {\hfil $\displaystyle {#}$}
&\hbox {$\displaystyle {#}$\hfil}
&\hbox {$\displaystyle {#}$\hfil}\cr
\big|{\partial\left(\phi^{(r)}-{\rm
Id}\right)}/{\partial\vet \omega} \big|_{\infty;\Wscr^{(r-1)}_{\frac{h_{r-1}}{4}}} &\leq & \big|{\partial\flusso^{(r)}}/{\partial\vet \omega} \big|_{\infty;\Wscr^{(r-1)}_{\frac{h_{r-1}}{4}}}
\big|{\partial\left(\phi^{(r-1)}-{\rm
Id}\right)}/{\partial\vet \omega} \big|_{\infty;\Wscr^{(r-1)}_{\frac{h_{r-1}}{4}}} 
\cr & \phantom{\le} &
+ 
\big|{\partial\left(\flusso^{(r)}-{\rm
Id}\right)}/{\partial\vet \omega} \big|_{\infty;\Wscr^{(r-1)}_{\frac{h_{r-1}}{4}}}
\cr
&\le &
(1+\mu_{r-1})\bar J_{r-1}+ \mu_{r-1} \,,
\cr
}}
$$
where we used the estimates $|\partial(\flusso^{(r)}-{\rm Id})/\partial \vet \omega|_{\infty; \Wscr^{(r-1)}_{\frac{h_r-1}{4}}} \le \mu_{r-1}$ (because of lemma~D.1 in~\cite{Poschel-1989}) and  $|\partial(\phi^{(r-1)}-{\rm Id})/\partial\vet  \omega|_{\infty; \Wscr^{(r-1)}_{\frac{h_r-1}{4}}} \le \bar J_{r-1}$ (that is the inductive formula related to the previous step, generalizing~\eqref{diseq:stime-phi2}--\eqref{def:Jbarra-da-0-a-2}).
Then, we can bound 
$\big|{\partial\left(\phi^{(r)}-{\rm
Id}\right)}/{\partial\vet \omega} \big|_{\infty;\Wscr^{(0)}_{\frac{h_{r-1}}{4}}}$ with $\bar J_r: = (1+\mu_{r-1})\bar J_{r-1}+ \mu_{r-1}\,$.\\
Recall now that we have to prepare for the next step, namely we have to cut out the resonant regions defined in~\eqref{frm:strisciarisonantekl}, i.e.,
$$
\Rscr_{\vet k,\vet \ell}^{(r)}=\left\{\vet \omega\in\Wscr^{(r-1)}:
\big| \vet k\cdot\vet \omega+\vet \ell\,\cdot\vet \Omega^{(r)}\circ\phi^{(r)}(\vet \omega)\big|
< 2\gamma/((r+1)K)^{\tau}\right\}\, .
$$
Thus, we need an upper bound on both the sup-norm of
$\vet \Omega^{(r)}\circ\phi^{(r)}$ and the Lipschitz constant of its
Jacobian.  The new transverse frequencies $\vet \Omega^{(r)}\big(\phi^{(r)}(\vet \omega^{(r)}))$ can be written
as
\begin{equation}
\vet \Omega^{(r)}\big(\phi^{(r-1)}\big(\flusso^{(r)}(\vet \omega^{(r)})\big) =
\vet \Omega^{(r-1)}\big(\phi^{(r-1)}\big(\flusso^{(r)}(\vet \omega^{(r)})\big) +
\Delta\vet \Omega^{(r)}\big(\flusso^{(r)}(\vet \omega^{(r)})\big)\ .
\label{frm:Omegoni-passo-r}
\end{equation}
Since we have $\max_{1\le j\le n_1}\sup_{\vet \omega\in
  \Wscr^{(r-1)}_{h_{r-1}}} \big|\Delta\Omega^{(r)}_j(\omega)\big|\leq
\gamma(\epsilon\Ascr)^r$, in view of the hypothesis~(i') in
proposition~\ref{pro:geometric-prop}, using again
lemma~\ref{lem:D1-poschel} so as to estimate
$|\partial(\flusso^{(r)}-{\rm Id})/\partial\vet
\omega|_{\infty; \Wscr^{(r-1)}_{\frac{h_r-1}{4}}}$, the following
inequality~\eqref{eq:Jr} at the previous inductive step $r-1$,
formula~\eqref{frm:Omegoni-passo-r} and the Cauchy inequality, we can
provide an estimate for
\begin{equation}
\vcenter{\openup1\jot\halign{
 \hbox {\hfil $\displaystyle {#}$}
&\hbox {$\displaystyle {#}$\hfil}\cr
\Bigg| &\frac{\partial\big(\vet \Omega^{(r)}\circ\phi^{(r)}-\vet \Omega^{(r-1)}\circ \phi^{(r-1)}\big)} {\partial\vet \omega}\Bigg|_{\infty;\Wscr^{(r-1)}_{\frac{h_{r-1}}{4}}}
\cr
 &\qquad \qquad \qquad = \left|\frac{\partial\big(\vet \Omega^{(r)}\circ\phi^{(r-1)}\circ\flusso^{(r)}-\vet \Omega^{(r-1)}\circ \phi^{(r-1)}\big)}{\partial\vet \omega}\right|_{\infty;\Wscr^{(r-1)}_{\frac{h_{r-1}}{4}}}
\cr
&\qquad \qquad \qquad \leq\left|\frac{\partial(\vet \Omega^{(r-1)}\circ \phi^{(r-1)})}{\partial\vet \omega}
\right|_{\infty;\Wscr^{(r-1)}_{\frac{h_{r-1}}{2}}}
\left|\frac{\partial\big(\flusso^{(r)}-\rm{Id}\big)}{\partial\vet \omega}
\right|_{\infty;\Wscr^{(r-1)}_{\frac{h_{r-1}}{4}}}
\cr
&\qquad \qquad \qquad \phantom{\leq}+
\left|\frac{\partial\Delta\vet \Omega^{(r)}}{\partial\vet \omega}
\right|_{\infty;\Wscr^{(r-1)}_{\frac{h_{r-1}}{2}}}
\,\left|\frac{\partial\flusso^{(r)}}{\partial\vet \omega}\right|_{\infty;\Wscr^{(r-1)}_{\frac{h_{r-1}}{4}}}
\cr
&\qquad \qquad\qquad \le J_{r-1} \mu_{r-1}+ \frac{2\gamma(\epsilon \Ascr)^r}{ h_{r-1} }(1+ \mu_{r-1})
\cr
&\qquad \qquad\qquad 
\le J_{r-1} \mu_{r-1} + \frac{ \mu_{r-1}}{2\sigma}(1+  \mu_{r-1})\,.
\cr
}}
\label{frm:scostamento-omegone}
\end{equation}
Therefore, we can bound the Jacobian of the function
$\vet \Omega^{(r)}\circ\phi^{(r)}$  with 
\begin{equation}
\left|\frac{\partial\big(\vet \Omega^{(r)}\circ\phi^{(r)}\big)}
{\partial\vet \omega}\right|_{\infty;\Wscr^{(r-1)}_{\frac{h_{r-1}}{4}}}
\leq \left(J_{r-1}  + \frac{\mu_{r-1}}{2\sigma}\right)(1+\mu_{r-1})=\mathrel{\mathop:} J_r\ .
\label{eq:Jr}
\end{equation}
Using again formula~\eqref{frm:delta-Delta-frequenze} and the inequality above at the previous inductive step, we can
bound the deterioration of the non-resonance conditions involving the
transverse frequencies. In fact, for $\vet \ell\in\interi^{n_2}$, $|\vet \ell|\le
2\,$, we have 
\begin{equation}
\vcenter{\openup1\jot\halign{
&\hbox {$\displaystyle {#}$\hfil}\cr
&\sup_{\omega\in \Wscr^{(r-1)}_{\frac{h_{r-1}}{4}}}
\left|\vet \ell\cdot\left[\vet \Omega^{(r)}\big(\phi^{(r)}(\vet \omega)\big)-
\vet \Omega^{(r-1)}(\phi^{(r-1)}(\vet \omega))\right]\right|
\cr
&\qquad\leq 2\max_j
\sup_{\vet \omega\in \Wscr^{(r-1)}_{\frac{h_{r-1}}{4}}}\left|\Omega^{(r-1)}_j\big(\phi^{(r-1)}\circ \flusso^{(r)}(\vet \omega)\big)-
\Omega^{(r-1)}_j(\phi^{(r-1)}(\vet \omega))\right|
\cr
& \qquad \phantom{\le }+2 \max_j \sup_{\vet \omega\in \Wscr^{(r-1)}_{\frac{h_{r-1}}{4}}}
\left|\Delta\Omega^{(r)}_j\big(\flusso^{(r)}(\vet \omega)\big)\right|
\cr
&\qquad\leq 2\left|\frac{\partial(\Omega^{(r-1)}\circ \phi^{(r-1)})}{\partial\vet \omega}
\right|_{\infty;\Wscr^{(r-1)}_{\frac{h_{r-1}}{2}}}
\,\max_j\,\sup_{\vet \omega\in \Wscr^{(r-1)}_{\frac{h_{r-1}}{4}}}\left|\flusso^{(r)}_j(\vet \omega)-\omega_j\right|
\cr
&\qquad \phantom{\le } +2\max_j\,\sup_{\vet \omega\in \Wscr^{(r-1)}_{\frac{h_{r-1}}{2}}}
\left|\Delta\Omega^{(r)}_j(\vet \omega)\right|
\cr
&\qquad\leq 2 J_{r-1}\gamma\sigma (\epsilon \Ascr)^r +2\gamma(\epsilon\Ascr)^r
= \frac {\mu_{r-1} h_{r-1}}{2}\left(\frac{1}{\sigma}+J_{r-1}\right)\ .
\cr
}}
\label{diseq:stima-var-freq-sec_passor}
\end{equation}

\noindent
Thus, for $\vet \omega\in\Wscr^{(r-1)}_{\frac{h_{r-1}}{4}}\,$, $0<|\vet k|\leq rK$ and
$|\vet \ell|\leq 2\,$, if $h_{r-1}$ satisfies the following condition
\begin{equation}
h_{r-1} \leq \min\left\{\frac{h_{r-2}}{4}\,,
\ \frac{1}{\max\big\{\frac{rK}{2}\,,\,\frac{1}{\sigma}\big\}+
 J_{r-1}} \frac{\gamma}{2^{r+1}(rK)^\tau}\right\}\,,
\label{diseq:hr-1-piccolo-quanto-basta}
\end{equation}
 defining $c_r = \sum_{s=1}^r \frac{1}{2^s}$  and starting from the non-resonance inequality at step $r-1$
(holding true on all the complex domain $\Wscr^{(r-1)}_{\frac{h_{r-1}}{4}}$),~\eqref{diseq:stima-var-freq-sec_passor} and~\eqref{def:mu_r-1} allows to obtain the non-resonance estimate
\begin{equation}
\vcenter{\openup1\jot\halign{
 \hbox {\hfil $\displaystyle {#}$}
&\hbox {$\displaystyle {#}$\hfil}\cr
\left| \vet k\cdot\vet \omega
+ \vet \ell\,\cdot\vet \Omega^{(r)}\big(\phi^{(r)}(\vet \omega)\big)\right|
&\ge
\frac{(2-c_r)\gamma}{(rK)^\tau}
-\sup_{\vet \omega\in \Wscr^{(r-1)}_{\frac{h_{r-1}}{4}}}
\left|\vet \ell\cdot\left[\vet \Omega^{(r)}\big(\phi^{(r)}(\vet \omega)\big)-
\vet \Omega^{(r-1)}(\phi^{(r-1)}(\vet \omega))\right]\right|
\cr
&\ge\frac{(2-c_r)\gamma}{(rK)^\tau}-
\frac {\mu_{r-1} h_{r-1}}{2}\left(\frac{1}{\sigma}+J_{r-1}\right)
\cr
&\ge
\frac{(2-c_{r+1})\gamma}{(rK)^\tau}\,,
\cr
}}
\label{frm:nonres-cond-diof-stepr-vecchio-blocco}
\end{equation}
where we used this same inequality at the previous inductive step, i.e., 
$$
\left| \vet k\cdot\vet \omega
+ \vet \ell\,\cdot\vet \Omega^{(r-1)}\big(\phi^{(r-1)}(\vet \omega)\big)\right| \ge \frac{(2-c_r)\gamma}{((r-1)K)^\tau}.
$$
\noindent
 Moreover, using the inequalities for the transverse frequencies at step $r-1$,~\eqref{diseq:stima-var-freq-sec_passor} and the estimate for $h_{r-1}$ in~\eqref{diseq:hr-1-piccolo-quanto-basta}, one can easily obtain the lower bound
\begin{equation}
\left| \vet \ell \cdot \vet \Omega^{(r)}\big(\phi^{(r)}(\vet \omega)\big)\right|
\geq
\left(2-c_{r+1}\right) \gamma\,,
\label{nonres-cond-trasversa-stepr}
\end{equation}
uniformly with respect to $\vet \omega\in\Wscr^{(r-1)}_{\frac{h_{r-1}}{4}}\,$, when
$0<|\vet \ell|\le 2$. 
 Now, we have to remove from the domain $\Wscr^{(r-1 )}$ the subset of resonant regions $\Rscr_{\vet k,\vet \ell}^{(r)}$ as they are defined in~\eqref{frm:strisciarisonantekl} for $rK<|\vet k|\leq (r+1)K\,$, $|\vet \ell|\leq 2\,$.
Then, if we require that 
\begin{equation}
h_r \le \min\left\{\frac{h_{r-1}}{4}\,,
\ \frac{1}{\max\big\{\frac{(r+1)K}{2}\,,\,\frac{1}{\sigma}\big\}+J_r}\,
\frac{\gamma}{2^{r+2}((r+1)K)^{\tau}}\right\}\ ,
\label{def:hr}
\end{equation}
for $rK<|\vet k|\leq (r+1)K\,$, $|\vet \ell|\leq 2\,$ and $\vet \omega \in \Wscr^{(r)}$, $|\bar\omega | < h_r$, we have
\begin{equation}
\vcenter{\openup1\jot\halign{
 \hbox {\hfil $\displaystyle {#}$}
&\hbox {$\displaystyle {#}$\hfil}\cr
\left| \vet k\cdot(\vet \omega + \bar{\vet \omega} )
+ \vet \ell\,\cdot\vet \Omega^{(r)}\big(\phi^{(r)}(\vet \omega + \bar{\vet \omega })\big)\right|
&\ge
\left| \vet k\cdot\vet \omega
+\vet  \ell\,\cdot\vet \Omega^{(r)}\big(\phi^{(r)}(\vet \omega)\big)\right|
\cr
&\phantom{\ge}-(r+1)Kh_{r}-2J_r h_{r}
\ge\frac{(2-c_{r+1})\gamma}{((r+1)K)^{\tau}}\ .
\cr
}}
\label{frm:nonres-cond-diof-stepr-nuovo-blocco}
\end{equation}

Therefore, for $r\ge 1\,$, it is convenient to resume all the conditions in the following way:
$$
\vcenter{\openup1\jot\halign{
 \hbox {\hfil $\displaystyle {#}$}
&\hbox {$\displaystyle {#}$\hfil}
&\qquad\hbox {$\displaystyle {#}$\hfil}\cr
\mu_{r-1}&=\frac{4\gamma\sigma(\epsilon\Ascr)^r}{h_{r-1}}\le
\frac{1}{2^r}
\ ,
&{{\rm (i)}} \cr
\bar J_r &= \bar J_{r-1}  (1+\mu_{r-1}) + \mu_{r-1}
\le e^{\tilde \mu_0}-1 \le (\epsilon \Ascr)\left(\frac{8\sigma \gamma}{h_0}+2^{\tau+3}\right) \le\frac{16 \epsilon \Ascr}{\Bscr}\,,
&{{\rm (ii)}} \cr
J_r &= \left(J_{r-1}  + \frac{\mu_{r-1}}{2\sigma}\right)(1+\mu_{r-1})
\le \left(J_0+\frac{\tilde\mu_{0}}{2\sigma}\right)
e^{\tilde\mu_{0}}
\leq e\left(J_0+\frac{1}{\sigma}\right)\ ,
&{{\rm (iii)}} \cr
h_r &\leq \min\left\{\frac{h_{r-1}}{4}\,,
\ \frac{1}{\max\big\{\frac{(r+1)K}{2}\,,\,\frac{1}{\sigma}\big\}+
 J_{r}}\,
\frac{\gamma}{2^{r+2}\big((r+1)K\big)^{\tau}}\right\}\ ,
&{{\rm (iv)}} \cr
}}
$$

\noindent
where $\tilde\mu_{r}=\sum_{s=r}^{\infty}\mu_{s}\,$ and $\Bscr$ is defined in~\eqref{def:bstorto}. 
The inequalities in~(ii) and~(iii) follow from~(i) and~\eqref{frm:hr}, by estimating $\tilde \mu_0$ as follows:
$$
\vcenter{\openup1\jot\halign{
 \hbox {\hfil $\displaystyle {#}$}
&\hbox {$\displaystyle {#}$\hfil}
&\hbox {$\displaystyle {#}$\hfil}\cr
\tilde \mu_0 &= &\sum_{s=0}^\infty \mu_s =  
\sum_{s=0}^\infty \frac{4\sigma\gamma(\epsilon \Ascr)^{s+1}}{h_s} 
 = 
\mu_0 + \epsilon \Ascr
\sum_{s=1}^\infty  \frac{4\sigma \gamma(\epsilon \Ascr)^{s}}{h_{s-1}}2^{\tau +2}
\cr
&
\le &\mu_0 + \epsilon \Ascr 2^{\tau +2} \sum_{s=1}^\infty  \frac{1}{2^s}
=
 ( \epsilon \Ascr )\left(\frac{4 \sigma \gamma}{h_0}+2^{\tau+2}\right) \le 
 \frac{8 \epsilon \Ascr}{\Bscr} \le 
1\,.
\cr
}}
$$
 The sequence $\{h_r\}_{r\ge 0}$ in~\eqref{frm:hr} has
been chosen so as to satisfy all the smallness
conditions~\eqref{diseq:h0-piccolo-quanto-basta}
and~(iv), which are required along this proof. 
The smallness
conditions~\eqref{def:mu1} and (i) on $\epsilon$ are satisfied since
$\epsilon<\epsilon^{*}_{{\rm ge}}\,$.  
Let us see it in details.\\
For $ r=1$, condition (i) is verified.
At step $r$ the condition is verified if
$$
\epsilon^r \le \frac{h_0}{(2^{\tau+2})^{r-1}}\frac{1}{4 \sigma \gamma(2\Ascr)^r} = \frac{1}{(2^{\tau+3}\Ascr)^r } \frac{h_0}{ \sigma \gamma}
2^{\tau}  \,,
$$
that is satisfied when $\epsilon<\epsilon_{\rm ge}^*$, which is defined in~\eqref{frm:soglia-geometrica}.

We have to prove now that $h_0$ and $h_r$, as defined in proposition~\ref{pro:geometric-prop}, satisfy condition (iv). 
Let us start with $h_0$. 
By definition $h_0\le \frac{\gamma}{4K^\tau e\left(J_0+\frac{1}{\sigma}\right)}$, then the inequality with the second term in~\eqref{diseq:h0-piccolo-quanto-basta} is satisfied.
Using that $\max\{\frac{K}{2}, \frac{1}{\sigma}\} \le  \max\{K, \frac{1}{\sigma}\} = \frac{1}{\eta}$, we obtain the following estimate
$$
\frac{1}{\max\{ \frac{K}{2}, \frac{1}{\sigma}\} + J_0} \frac{\gamma}{4K^{\tau}} 
\ge \frac{1}{\frac{1}{\eta} + e\left(J_0+\frac{1}{\sigma}\right)} \frac{\gamma}{4K^{\tau}} 
\ge  \frac{1}{ \max\{\frac{1}{\eta},e\left(J_0+\frac{1}{\sigma}\right)\}} \frac{\gamma}{8K^{\tau}} = h_0\,.
$$
By definition $h_r = \frac{h_{r-1}}{2^{\tau+2}} $, therefore the inequality with the first term in (iv) is satisfied. Therefore, we only need to show the second inequality. We need a lower bound for the following term:
$$
\vcenter{\openup1\jot\halign{
 \hbox {\hfil $\displaystyle {#}$}
&\hbox {\hfil $\displaystyle {#}$\hfil}
&\hbox {$\displaystyle {#}$\hfil}\cr
\frac{1}{\max\big\{\frac{(r+1)K}{2}\,,\,\frac{1}{\sigma}\big\}+
 J_r}\, 
\frac{\gamma}{2^{r+2}\big((r+1)K\big)^{\tau}} 
&\ge&
\frac{1}{\frac{(r+1)}{2}\frac{1}{\eta}+
e\left(J_0+\frac{1}{\sigma}\right)}\,
\frac{\gamma}{2^{r+2}\big((r+1)K\big)^{\tau}} 
\cr
&\ge &
\frac{1}{(r+1) \max\{ \frac{1}{\eta}, e\left(J_0+\frac{1}{\sigma}\right)\}}
\frac{\gamma}{2^{r+2}\big((r+1)K\big)^{\tau}} \,.
\cr
}}
$$
Since $h_r = \frac{h_0}{(2^{\tau+2})^r}$ and $(r+1)^{\tau+1} 2^{r-1} \le (2^r)^{\tau+1} 2^r = (2^{\tau+2})^r$,
the inequality is satisfied.
This concludes the proof.
\end{proof}

\subsubsection{Explicit estimates for the volume of the resonant strips}
\label{subsec:explicit-estimates}

In order to prove theorem~\ref{thm:main-theorem}, we need to focus on a set of tori with ``Diophantine'' frequencies as defined in definition~\ref{def:dioph-freq}.
First, let us denote by $\Kscr_{\vet \ell}^{(r)}$ the closed convex hull of the gradient set
\begin{equation}
\Gscr_{\vet \ell}^{(r)}=
\left\{\partial_{\vet \omega}\big[\scalprod{\vet \ell}{\vet \Omega^{(r)}\big(\phi^{(r)}(\vet \omega)\big)}\big]
\>:\ \vet \omega\in\Wscr^{(r)}\right\}\ .
\label{frm:gradient-set-passor}
\end{equation}
We emphasize that the previous definition is analogous to that of
$\Gscr_{\vet \ell}^{(0)}$ given in the context of hypothesis~(b') in
lemma~\ref{lem:Ham-expansion}. That assumption ensures that the
Euclidean distance from the initial set of frequencies is such that
${\rm dist}(\vet k,\Kscr_{\vet \ell}^{(0)})\geq 2\Theta_0$ for
$\vet k\in\interi^{n_1}\setminus\{\vet 0\}\,$.  For $\epsilon$ small enough, the
closed convex hull $\Kscr_{\vet \ell}^{(r)}$ that includes $\Gscr_{\vet \ell}^{(r)}$ is
just a little bigger than the original one, so that also at the $r$-th
step it is safely away from the integer vectors, i.e., ${\rm
  dist}(\vet k,\Kscr_{\vet \ell}^{(r)})\geq \Theta_0$
$\forall\ \vet  k\in\interi^{n_1}\setminus\{\vet 0\}\,$. In order to prove that, the following further smallness condition is needed:
\begin{equation}
 \epsilon < \frac{\Theta_0}{ 2^{6}\Ascr 
\left(\frac{1}{\sigma}+J_0\right)}\min\left\{ \frac{h_0}{\sigma\gamma}\,, \frac{1}{2^\tau}\right\}=
 \frac{\Theta_0 \Bscr}{ 2^{6} \Ascr 
\left(\frac{1}{\sigma}+J_0\right)}\,.
\label{frm:eps-piccolo-per-guscio-convesso}
\end{equation}

\noindent
Let us see it in details.
In order to control the growth of $\Kscr_{\vet \ell}^{(r)}$, it suffices to show that the following estimate holds true 
\begin{equation}
\left|\partial_{\vet \omega} \left[\vet \ell \cdot \vet \Omega^{(r)}\left(\phi^{(r)}\left(\vet \omega\right)\right)- \vet \ell\cdot \vet  \Omega^{(r-1)}\left(\phi^{(r-1)}\left(\vet \omega \right)\right)\right]\right| 
\le \frac{\Theta_0}{2^r} \quad  \forall \ r \ge 1\,,
\label{frm:delta-omegone}
\end{equation}
with $\vet \omega\in \Wscr^{(r)}$. Since $\Wscr^{(r)} \subset \Wscr^{(r-1)}\subset  \Wscr^{(r-1)}_{\frac{h_{r-1}}{4}}$, using the inequality in~\eqref{frm:scostamento-omegone}, that has been verified in the proof of proposition~\ref{pro:geometric-prop}, we can bound the function as follows: 
\begin{equation}
\vcenter{\openup1\jot\halign{
 \hbox {\hfil $\displaystyle {#}$}
&\hbox {\hfil $\displaystyle {#}$\hfil}
&\hbox {$\displaystyle {#}$\hfil}\cr
\left|\partial_{\vet \omega} \left[\vet \ell \cdot \vet  \Omega^{(r)}\left(\phi^{(r)}\left(\vet \omega\right)\right)-\vet \ell\cdot \vet \Omega^{(r-1)}\left(\phi^{(r-1)}\left(\vet \omega\right)\right)\right]\right|
&\le 
2\left[J_{r-1} \mu_{r-1} + \frac{ \mu_{r-1}}{2\sigma}(1+  \mu_{r-1})\right]\, 
\cr
& \le 2\mu_{r-1}\left[ e\left(J_0 + \frac{1}{\sigma}\right) + \frac{1}{\sigma}\right]
\cr
& \le 8 \mu_{r-1}\left(J_0 + \frac{1}{\sigma}\right) 
\cr
}}
\label{frm:stima-delta-omegone}
\end{equation}
where $J_{r-1}$, that is the Lipschitz constant about the Jacobian of $\vet \Omega^{(r-1)}(\vet \omega^{(0)})$, has been estimated as in~(iii) at the end of the proof of proposition~\ref{pro:geometric-prop}.
Then, if $\epsilon$ satisfies the condition in~\eqref{frm:eps-piccolo-per-guscio-convesso}, by using also the definition and the estimate of $\mu_{r-1}$ in~(i) at the end of the proof of proposition~\ref{pro:geometric-prop} and the definition of $h_{r-1}$ in~\eqref{frm:hr}, for $r \ge 2$
 we have 
$$
\vcenter{\openup1\jot\halign{
 \hbox {\hfil $\displaystyle {#}$}
&\hbox {\hfil $\displaystyle {#}$\hfil}
&\hbox {$\displaystyle {#}$\hfil}\cr
8\mu_{r-1} \left(\frac{1}{\sigma}+ J_0\right) &\le&
 8 \frac{4  \sigma \gamma (\epsilon \Ascr)^{r-1} (\epsilon \Ascr) 2^{\tau+2} }{h_{r-2}} \left( \frac{1}{\sigma}+ J_0\right) 
\cr
&\le&
8 \frac{(\epsilon \Ascr) 2^{\tau+2}}{2^{r-1}} \left( \frac{1}{\sigma}+ J_0\right) 
\cr
&\le&
\frac{(\epsilon \Ascr) 2^{\tau+6} }{2^{r}} \left( \frac{1}{\sigma}+ J_0\right)
\le
 \frac{\Theta_0}{2^r} \,,
\cr
}}
$$
where the condition~\eqref{frm:eps-piccolo-per-guscio-convesso} is used in the last inequality.
For $r=1$, the r.h.s. in the estimate~\eqref{frm:stima-delta-omegone} can be reformulated as 
$$
8\mu_0 \left(\frac{1}{\sigma} + J_0\right) \le 8 \frac{4 \sigma \gamma \epsilon \Ascr}{h_0}  \left(\frac{1}{\sigma} + J_0\right)
\le 2^5   \frac{\sigma \gamma \epsilon \Ascr}{h_0} \left(\frac{1}{\sigma} + J_0 \right)\le \frac{\Theta_0}{2} \,.
$$
where we used again the estimate~(i) at the end of the proposition~\ref{pro:geometric-prop} and the condition~\eqref{frm:eps-piccolo-per-guscio-convesso}. This ends the justification of~\eqref{frm:delta-omegone}.

\noindent
Therefore, we obtain the following lower bound:
$$
\vcenter{\openup1\jot\halign{
 \hbox {\hfil $\displaystyle {#}$}
&\hbox {\hfil $\displaystyle {#}$\hfil}
&\hbox {$\displaystyle {#}$\hfil}\cr
\dist(\vet k, \Kscr_{\vet \ell}^{(r)})  &\ge &\dist(\vet k, \Kscr_{\vet \ell}^{(0)}) -\sum_{i=0}^{r-1}  \max_{\vet \omega}  \left|\partial_{\vet \omega} \left[\vet \ell \cdot \vet \Omega^{(i+1)}\left(\phi^{(i+1)}\left(\vet \omega\right)\right)-\vet \ell\cdot \vet  \Omega^{(i)}\left(\phi^{(i)}\left(\vet \omega\right)\right)\right]\right|  
\cr
&\ge & 2\Theta_0 - \sum_{i=0}^{\infty} \frac{\Theta_0}{2^{i+1}} \ge \Theta_0\,.
\cr
}}
$$
The estimate of the volume of resonant zones that are removed at each
step is based on lemma~8.1 of P\"oschel paper~\cite{Poschel-1989}, that we
report here (the Lebesgue measure is denoted by
${\rm m}(\cdot)$).
\begin{lemma}\label{lem:Poeschel}
If ${\rm dist}(\vet k,\Kscr_{\vet \ell}^{(r)})=\Theta>0$ then
$$
{\rm m}(\Rscr^{(r)}_{\vet k,\vet \ell})\leq \frac{4\gamma}{((r+1)K)^\tau}
\frac{D^{n_1-1}}{\Theta}\ ,
$$
where $D$ is the diameter of $\Wscr^{(0)}$ with respect to the
sup-norm.
\end{lemma}
The volume of the resonant regions must be compared with respect to
the initial set $\Wscr^{(0)}$.  To this end, we estimate the measure
of $\phi^{(r)}(\Rscr^{(r)}_{\vet k,\vet \ell})\,$ in the original coordinates
$\vet \omega^{(0)}$, recalling that $\phi^{(r)}$ is the inverse function of
$\vet \omega^{(r)}(\vet \omega^{(0)})$. Using lemma~\ref{lem:Poeschel} and
assumption~\eqref{frm:eps-piccolo-per-guscio-convesso}, for
$\vet k\in\interi^{n_1}\setminus\{\vet 0\}$ with $rK< |\vet k| \le (r+1)K$ and $|\vet \ell |< 2$ we have
\begin{equation}
{\rm m}\Big(\phi^{(r)}\big(\Rscr^{(r)}_{\vet k,\vet \ell}\big)\Big)\le
\frac{4\gamma D^{n_1-1}}{\Theta_0((r+1)K)^\tau}\,\sup_{\vet \omega\in\Wscr^{(r)}}
\det\left(\frac{\partial\phi^{(r)}}{\partial \vet \omega}\right)\ .
\label{frm:stima-zona-risonante-in-omeghini-iniz}
\end{equation}
Starting from the first inequality
in~\eqref{frm:diseq:Lipschitz-definitive} and using the well known
Gershgorin circle theorem, we control the actual
expansion of the resonant zones due to the stretching of the
frequencies, by verifying that
\begin{equation}
\sup_{\vet \omega\in\Wscr^{(r)}}
\det\left(\frac{\partial\phi^{(r)}}{\partial \vet \omega}\right)
\le 2
\qquad{\rm when}\quad
\epsilon\le\frac{\Bscr\log 2}{  16 \Ascr n_1^2}\ .
\label{frm:stima-dilatazione-frequenze-veloci}
\end{equation}
Using the
inequalities~\eqref{frm:stima-zona-risonante-in-omeghini-iniz}--\eqref{frm:stima-dilatazione-frequenze-veloci},
we easily get a final estimate of the total volume of the
resonant regions included in $\Wscr^{(0)}$:
\begin{equation}
\vcenter{\openup1\jot\halign{
 \hbox {\hfil $\displaystyle {#}$}
&\hbox {\hfil $\displaystyle {#}$\hfil}
&\hbox {$\displaystyle {#}$\hfil}\cr
\sum_{r=1}^{\infty}\,\sum_{{rK<|\vet k|\leq (r+1)K}\atop{|\vet \ell|\leq 2}}
{\rm m}\Big(\phi^{(r)}\big(\Rscr^{(r)}_{\vet k,\vet \ell}\big)\Big) &\leq
&c_{n_2}\sum_{r=1}^{\infty}\,\sum_{rK<|\vet k|\leq (r+1)K}
\frac{8\gamma D^{n_1-1}}{\Theta_0\,((r+1)K)^\tau}
\cr
&\leq &\gamma\ \frac{2^{n_1+3} c_{n_2} D^{n_1-1}}{\Theta_0\,K^{\tau-n_1}}
\sum_{r=2}^{\infty}\frac{1}{r^{\tau-n_1+1}}\ ,
\cr
}}
\label{frm:measure}
\end{equation}
where $c_{n_2}=(n_2+1)(2n_2+1)$ is the maximum number of polynomial
terms having degree $\le 2$ in the transverse variables $(\vet z,
\imunit \bar {\vet z})\,$ and the maximum number of $\vet k$ such that
$rK<|\vet k|\le (r+1)K$ is estimated with $K \cdot
2^{n_1}((r+1)K)^{n_1-1}$.  The last series is convergent if $\tau>n_1$
and the sum is of order $\Oscr(\gamma)$.

\subsection{Proof} \label{sec:proof-th1}
\noindent
The proof is a straightforward combination of lemma~\ref{lem:stime-valut-insiemi} and
propositions~\ref{pro:analytic-prop} and~\ref{pro:geometric-prop}, which
summarize the analytic study of the convergence of our
algorithm and the geometric part, respectively.
We sketch the argument.

According to lemma~\ref{lem:Ham-expansion} the family of Hamiltonians
$\Hscr^{(0)}$, parametrized by the frequency vectors $\vet
\omega^{(0)}$ and defined on a real domain, can be extended to a
complex domain $\Dscr_{\rho,\sigma,R}\times\Wscr_{h_0}$ with suitable
parameters; moreover, their expansions can be written as $H^{(0)}$
in~\eqref{frm:H(0)}.  Possibly modifying the values of parameters
$\gamma$ and $\tau\,$, we can choose $\gamma$ and $\tau>n_1$ such that
the property~(a') of lemma~\ref{lem:Ham-expansion} is still satisfied
and the estimate of the resonant volume in the last row
of~\eqref{frm:measure} is smaller than ${\rm m}(\Wscr_{h_0})$, i.e.,
we impose that
\begin{equation}
\label{frm:res-measure}
\gamma \frac{2^{n_1+3}(2n_2^2+3n_2+1)D^{n_1-1}}{\Theta_0 K^{\tau-n_1}} \sum_{r=2}^{+\infty} \frac{1}{r^{\tau -n_1+1}} < {\rm m}\left( \Wscr^{(0)} \right)\,.
\end{equation}
Let us consider values of the small parameter
$\epsilon$ such that $\epsilon<\epsilon^{\star}$, with
\begin{equation}
\epsilon^{\star}=
\min\left\{\epsilon^{\star}_{{\rm ge}}
\,,\,\frac{\Theta_0 \,  \Bscr}{2^{6} \Ascr 
\left(\frac{1}{\sigma}+J_0\right)}\,,\,\frac{ \Bscr \, \log 2 }{ 16 \Ascr\, n_1^2}\right\}\ ,
\label{frm:soglia-finale}
\end{equation}
where $\epsilon^{\star}_{{\rm ge}}$ is defined
in~\eqref{frm:soglia-geometrica}. Let us remark that if $\epsilon <
\epsilon^*_{\rm ge}$ then the parameter $\epsilon$ is also smaller
enough to deal with the analytic part of the proof, i.e., $\epsilon
<\epsilon^*_{\rm an}$; this can be easily verified by comparing
definitions~\eqref{frm:soglia-analitica}
and~\eqref{frm:soglia-geometrica}.

Recall that the $r$-th normalization step of the formal algorithm
described in section~\ref{sec:algorithm} can be performed if the
non-resonance conditions~\eqref{frm:nonres_a} and~\eqref{frm:nonres_b}
are satisfied up to order $r$.  Assuming the threshold value
$\epsilon^{\star}$ as in~\eqref{frm:soglia-finale},
lemma~\ref{lem:Ham-expansion} and proposition~\ref{pro:geometric-prop}
ensure that the first step can be performed for every frequency vector
$\vet \omega^{(0)}\in\Wscr^{(0)}_{h_0}\,$.

We now proceed by induction.  Let us suppose that $r-1$ steps have
been performed and proposition~\ref{pro:geometric-prop} applies.  In
view of the non-resonance condition~\eqref{frm:nonres-cond-diof} the
$r$-th normalization step can be performed.  We now check that
proposition~\ref{pro:geometric-prop} applies again. By construction,
both $\vet \omega^{(r)}(\vet \omega^{(0)})$ and $\vet \Omega^{(r)}(\vet \omega^{(0)})$ are
analytic functions on
$\phi^{(r-1)}\big(\Wscr^{(r-1)}_{h_{r-1}}\big)\,$. Moreover, 
lemma~\ref{lem:stime-valut-insiemi} applies and the
estimates~\eqref{frm:stima-variazione-frequenze-+-esplicita} of the
shift of the frequencies holds true.  Thus,
proposition~\ref{pro:geometric-prop} can be applied at the $r$-th step
which completes the induction.

We conclude that the non-resonance
conditions~\eqref{frm:nonres-cond-diof} hold true for $r\ge 0$ and
that the sequence of frequency vectors $\big\{\big(\vet
\omega^{(r)}(\vet \omega^{(0)})\,,\, \vet \Omega^{(r)}(\vet
\omega^{(0)})\big)\big\}_{r\ge 0}$ is Diophantine for $\vet
\omega^{(0)}\in\lim_{r\to\infty}\phi^{(r)}\big(\Wscr^{(r)}_{h_r}\big)=
\bigcap_{r=0}^{\infty}\phi^{(r)}\big(\Wscr^{(r)}_{h_r}\big)\,$, where
we used the inclusion relation
$\phi^{(r)}\big(\Wscr^{(r)}_{h_{r}}\big)
\subset\phi^{(r-1)}\big(\Wscr^{(r-1)}_{h_{r-1}}\big)\,$ between open
sets.  Finally, also hypothesis~(f') of
proposition~\ref{pro:analytic-prop} is satisfied because $\epsilon <
\epsilon^* < \epsilon^*_{\rm an}$ and so, for $\vet
\omega^{(0)}\in\bigcap_{r=0}^{\infty}\phi^{(r)}\big(\Wscr^{(r)}_{h_r}\big)$,
there exists an analytic canonical transformation $\Phi_{\vet
  \omega^{(0)}}^{(\infty)}$ which gives the initial Hamiltonian
$H^{(0)}$ the normal form~\eqref{frm:H(infty)-espansione}.  Since
$\bigcap_{r=0}^{\infty}\phi^{(r)}\big(\Wscr^{(r)}_{h_r}\big)$ is a
countable intersection of open sets, it is measurable and
$$
{\rm m}\bigg(\bigcap_{r=0}^{\infty}\phi^{(r)}\big(\Wscr^{(r)}_{h_r}\big)\bigg)
\ge {\rm m}\big(\Wscr^{(0)}\big)-
\sum_{r=1}^{\infty}\,\sum_{{rK<|\vet k|\leq (r+1)K}\atop{|\vet \ell|\leq 2}}
{\rm m}\Big(\phi^{(r)}\big(\Rscr^{(r)}_{\vet k,\vet \ell}\big)\Big)>0\ ,
$$ where we have taken into account the fact that the complex
extension radius $h_r\to 0$ for $r\to\infty$, the
estimate~\eqref{frm:measure} and the initial choice of the parameters
$\gamma$ and $\tau$ at the beginning of the present section, in such a
way that inequality~\eqref{frm:res-measure} is satisfied. This
concludes the argument proving theorem~\ref{thm:main-theorem}.

\section{Conclusions}
\label{sec:final-discussions}
We have discussed the convergence of an algorithm for proving the
existence of lower dimensional elliptic tori for a family of
Hamiltonians parametrized by the frequency vector $\vet \omega$ of the
motion on the torus. Let us comment a little bit more the result.

From a practical point of view, the small change of the frequencies is
a weakness of our method, because in an explicit construction one is
forced to check at every step that the wanted non-resonance conditions
are satisfied.  This cannot be avoided, because the change of the
frequencies is an effect of the perturbation that is not known at the
beginning of the procedure. By suitably controlling the small
corrections on the frequencies introduced by the algorithm, we can
prove the existence of tori for a set $\Uscr^{(\infty)}$ of positive
measure of initial frequencies, but at the same time we cannot obtain
the same result for one specific torus or, equivalently, for a
specific Diophantine frequency vector $\vet \omega$. Indeed, if on one
hand we could add at every step a translation of the actions of the
torus in order to fix the frequency vector, as suggested by Kolmogorov
in the case of full dimensional tori, at the same time we cannot fix
the transverse frequencies $\Omega_j$, that depend generically on
$\vet\omega$. However, the geometric construction above shows that
with high probability the method will work fine at every order; in
particular, it is even less probable to encounter a resonance in the
explicit calculations, because the algorithm is iterated up to a
finite order. In~\cite{San-Dan}, in the case of the classic Kolmogorov
algorithm for constructing full dimensional tori (but without the
translations to fix the frequency vector), the authors present an
approach to determine {\it a posteriori} the initial frequencies in
order to converge to a picked value of $\vet \omega$. With some
additional work, we believe that the same approach could be extended
to the case of lower dimensional elliptic tori.

Let us now comment the threshold value $\epsilon^{\star}$ on the small
parameter $\epsilon$, that is explicitly defined
in~\eqref{frm:soglia-finale}. Although the definition involves many
parameters, we might produce an asymptotic estimate of the volume of
the resonant regions for $\epsilon\to0$.  In view of
inequality~\eqref{frm:measure}, it is $\Oscr(\gamma)$. Moreover, one
can easily show that $\epsilon^{\star}$ is $\Oscr(\gamma^3)$, since
every element in~\eqref{frm:soglia-finale} is $\Oscr(1/\Ascr)$ and
$\Ascr \sim 1/\gamma^3$. Thus the complement of the set of the
invariant elliptic tori, i.e.,
$\Wscr^{(0)}\setminus\bigcap_{r=0}^{\infty}\phi^{(r)}\big(\Wscr^{(r)}_{h_r}\big)\,$,
has a measure estimated by $\Oscr(\epsilon^{1/3})\,$.  This law is
strictly related to the fact that each normalization step requires the
solution of three homological equations and this is related with the
occurrence of the factor $3$ in the accumulation rules of the small
divisors that are stated in lemma~\ref{lem:small-div-acc}. However, an
estimate such that the measure of the complementary set of the
invariant tori should be $\Oscr\left(\epsilon^{1/3}\right)$ is not
optimal.  For instance, this same quantity has been proved to be
smaller than a bound $\Oscr(\epsilon^{b_1})\,$, with $b_1<1/2$,
in~\cite{Bia-Chi-Val-2003} and~\cite{Bia-Chi-Val-2006}, where the
existence of elliptic tori in a couple of planetary problems is
investigated.  Let us emphasize that considering a planetary problem
(as in~\cite{Gio-Loc-San-2014}) makes much easier the control of the
transversal frequencies because they are $\Oscr(\epsilon)$. Here,
careful additional estimates have been introduced, e.g., in
subsection~\ref{subsec:explicit-estimates}, in order to deal with this
difficulty. On the other hand, it is well known that purely analytic
estimates are usually unrealistically small. For this reason we did
not pay attention to produce optimal estimates.

The interest of such an algorithm is indeed related to the fact that,
being fully constructive, it can be used to obtain semi-analytical
results in realistic problems, by explicitly performing a number of
perturbation steps (e.g., with the help of an algebraic manipulator
especially designed for these kind of computations,
see~\cite{Gio-San-Chronos-2012}) and significantly improving both the
applicability threshold and the estimate of the measure. Furthermore,
by implementing a suitable scheme of computer-assisted estimates, the
semi-analytical results could be made also fully rigorous (see,
e.g.,~\cite{Cel-Chi-2007} and~\cite{Loc-Gio-2000} in the case of full
dimensional tori). Computer-assisted proofs of existence of elliptic
lower dimensional tori have been given e.g. in~\cite{Luq-Vil-2011},
where the authors described a constructive algorithm without the
request of using action-angle canonical coordinates nor the
quasi-integrability of the systems. Also the Poincar\'e-Lindstedt
method can be used for computing an approximation of the periodic or
quasi-periodic orbits (see, e.g.,~\cite{Eliasson-1988.1},
\cite{Gallavotti-1994.1}, \cite{Chi-Fal-94} and \cite{Gen-Mas-96.2}).
Still, these methods can be used to construct a particular solution of
the Hamilton equations, but further work is necessary in order to
study the dynamics in its neighborhood. On the contrary, the fact that
the method is based on a construction of a suitable normal form, i.e.,
we produce not just a specific solution but also a local holomorphic
Hamiltonian, makes it an optimal starting point for studying the
dynamics in the neighborhood of the tori.

From a theoretical point of view, a description of the approach to
follow in order to investigate the stability in the neighborhood of
invariant tori can be found in~\cite{Mor-Gio-1995}. In that article,
in the case of full dimensional tori, Morbidelli and Giorgilli showed
that, if $\rho$ is the distance with respect to the torus, therefore
the upper bound on the diffusion rate is superexponentially small, as
it roughly goes as $\exp(\exp(-1/\rho))$. This phenomenon is strictly
related to the fact that the more an invariant solution is robust
(with respect to the perturbation) the more its neighborhood is filled
by other invariant tori. This realizes a sort of freezing of the
diffusion that is made extremely slow (see~\cite{Giorgilli-1997.4}).
A similar phenomenon happens in the case of lower dimensional elliptic
tori, which are surrounded by a stable zone filled by tori of maximal
dimension. The slow diffusion from invariant tori is particularly
useful in those physical problems for which an eventual proof of
stability for perpetual times is much more than necessary, while it is
more appropriate to refer to the {\it effective stability}. This
concept is easy to define in the case, e.g., of a planetary system:
the typical lifetime of a star, despite very long, is obviously
finite; therefore, in this case, it would be enough to show that the
motion is forced to be confined in a small subset of the phase space
for a time that exceeds the expected lifetime of the star.  To this
purpose, a simple tool for providing a lower bound for the time needed
to diffuse from small region in the neighborhood of a lower
dimensional or full dimensional invariant torus is given by the
Birkhoff normal form. Suitable estimates of the remainder of Birkhoff
normal form, can be used in order to investigate the stability in
Nekhoroshev sense (see,
e.g.,~\cite{Mor-Gio-1995},~\cite{Gio-Loc-San-2009}
and~\cite{San-Loc-Gio-2013}), thus for times exponentially long with
respect to the perturbation. In~\cite{Car-Loc-2021}, Birkhoff normal
form is used in this way in order to explain the stability numerically
observed in the FPUT model, that is a chain of particle
interacting with anharmonic potentials. Moreover, also in this case,
the computation of the stability time can be made rigorous with a
computer-assisted scheme of estimates (for an introduction,
see~\cite{Car-Loc-2020}, where it is considered the simpler case of
Birkhoff normal forms in the neighborhood of an elliptic equilibrium
point). Finally, the convenience of combining different normal forms
and of doing expansions nearby elliptic tori is emphasized
in~\cite{Car-Loc-San-Vol-2021}, where an efficient construction of a
one-dimensional elliptic torus turns out to be an advantageous
intermediate step before the construction of a full dimensional KAM
tori in its neighborhood.

\section*{Acknowledgments}
The author acknowledges U. Locatelli for suggesting the topic and for
the help in structuring the proof.  This work was partially supported
by the MIUR-PRIN project 20178CJA2B -- ``New Frontiers of Celestial
Mechanics: theory and Applications''. The author acknowledges also the
MIUR Excellence Department Project awarded to the Department of
Mathematics of the University of Rome ``Tor Vergata'' (CUP
E83C18000100006), the ``Beyond Borders" programme of the University of
Rome ``Tor Vergata" through the project ASTRID (CUP E84I19002250005)
and the ``Progetto Giovani2019'' programme of the National Group of
Mathematical Physics (GNFM-INdAM) through the project “Low-dimensional
Invariant Tori in FPU-like Lattices via Normal Forms”.

\section*{Appendix}
\appendix
\section{Proof of lemma~\ref{lem:iterative-estimates}}
\label{app:stimeiterative}
Let us first check the estimates for the generating function $\chi_0^{(r)}$. Using the inequalities in~\eqref{frm:f_l^(r-1,s)-lemmone}, we come to the following estimate:
$$
\vcenter{\openup1\jot 
\halign{
$\displaystyle\hfil#$&$\displaystyle{}#\hfil$&$\displaystyle#\hfil$\cr
\left(\frac{2e}{\rho\sigma}+\frac{e^2}{R^{2}}\right) \frac{1}{\delta_r^2}\| \chi_0^{(r)} \|_{1- d_{r-1}} &\le 
 \left(\frac{2e}{\rho\sigma}+\frac{e^2}{R^{2}}\right)\frac{1}{\delta_r^2 a_r} \| f_0^{(r-1,r)}\|_{1- d_{r-1}} 
\cr
&\le 
 \left(\frac{2e}{\rho\sigma}+\frac{e^2}{R^{2}}\right) \Ebarra M^{3r-3} \frac{4\pi^4 r^{4+\tau} K^\tau}{ \gamma }  \Vscr\left(\Hscr_0^{(r-1,r)}\right) \nu_{r-1,r} \, .
\cr
}}
$$
By definition of the evaluation operator, $\Vscr\big(\{r\}\big) = r^{4+\tau}$; then, using the property for the product, $\Vscr\big(\{r\}\big) \cdot\Vscr\big(\Hscr_0^{(r-1,r)}\big) = \Vscr\big(\Gscr_0^{(r)}\big) $. Moreover, constants can be limited by $M$, so the first inequality in~\eqref{frm:generatrici-lemmone} holds true. The other estimates that appear in formula~\eqref{frm:generatrici-lemmone} and concern with the generating functions $\chi_1^{(r}$ and $\chi_2^{(r}$ can be verified in the same way.

In order to prove the other estimates, it will be useful to preliminarily give an estimate for the intermediate functions that compose the Hamiltonian in the different steps of the algorithm, i.e., we want to prove that the following inequalities hold true:
\begin{equation}
\vcenter{\openup1\jot 
\halign{
$\displaystyle\hfil#$&$\displaystyle{}#\hfil$&$\displaystyle#\hfil$\cr
\|f_{\ell}^{(*; \, r,s)}\|_{1-d_{r-1}-(\#\delta_r)} \leq &\frac{\Ebarra M^{3s-\zeta_\ell}}{2^{\ell}}\,  
\Vscr \left(\Hscr_\ell^{(*;\, r,s)}\right)
\,\nu_{r,s}^{*} ,
\cr
}}
\label{frm:f_l^(I,II,III;r,s)-lemmone}
\end{equation}
for $0\le\ell\le 2\,,\ s\ge r$ or $\ell\ge 3\,,\ s\ge 0$, where the generic symbol `$*$' on the top of $f_\ell$ can be replaced by \rm{I} or \rm{II}, while in the lower index `$\#$' can assume the values $1$ or $2$ for the corresponding terms. 

\noindent
By definition~\eqref{frm:fI} of $f_0^{({\rm I};\,r,s)}$, it is easy to check that for $r<s<2r$, 
$$
\| f_0^{({\rm I};\,r,s)} \|_{1- d_{r-1}- \delta_r} \le
\Ebarra M^{3s-3}\Vscr\left(\Hscr_0^{(r-1,s)}\right) \nu_{r-1,s} .
$$
Then, using the equation $\Hscr_0^{({\rm I}; r,s)} = \Hscr_0^{(r-1,s)}$ and the inequality $  \nu_{r-1,s} \le  \nu_{r,s}^{\rm (I)} $, the estimate in~\eqref{frm:f_l^(I,II,III;r,s)-lemmone} is proved for what concerns the case with $r<s<2r$.

\noindent
For $\ell=0$, $s\ge 2r $ or $\ell =1,2$, $s \ge r$ or $\ell\ge 3$, $s\ge0$,  we have that
$$
\vcenter{\openup1\jot 
\halign{
$\displaystyle\hfil#$&$\displaystyle{}#\hfil$&$\displaystyle#\hfil$\cr
\|  f_\ell^{({\rm I};r,s)} \|_{1- d_{r-1} - \delta_r} &\le&
\sum_{j=0}^{\lfloor s/r \rfloor}  \frac{1}{j!}  \| \Lie{\chi_0^{(r)}}^j  f_{\ell+2j}^{(r-1,s-jr)} \|_{1- d_{r-1} - \delta_r}
\cr
&\le &\sum_{j=0}^{\lfloor s/r \rfloor} M^{(3r-2)j}\Vscr\left(\Gscr_0^{(r)}\right)^j \nu_{r-1,r}^j \Ebarra \frac{M^{3(s-jr)- \zeta_{\ell+2j}}}{2^{\ell+2j}} \Vscr\left(\Hscr_{\ell+2j}^{(r-1, s-jr)}\right)\nu_{r-1,s-jr} \, .
\cr
}}
$$
The estimate in~\eqref{frm:f_l^(I,II,III;r,s)-lemmone} follows from the definition of $\Hscr_{\ell}^{(\rm{I}; r, s)}$ in the statement, that of $\nu_{r,s}^{(\rm{I})}$ in~\eqref{frm:seqnu} and in view of the elementary inequality $\zeta_{\ell+2j} +2j \ge \zeta_\ell$ $\forall \ j \ge 0$, so that $M^{3s-2j-\zeta_{\ell+2j}} \le M^{3s-\zeta_\ell} $.

\noindent
Since $\chi_1^{(r)}$ satisfies the inequality:
$$
\vcenter{\openup1\jot 
\halign{
$\displaystyle\hfil#$&$\displaystyle{}#\hfil$&$\displaystyle#\hfil$\cr
\left(\frac{2e}{\rho\sigma}+\frac{e^2}{R^{2}}\right)\frac{1}{\delta_r^2}\| \chi_1^{(r)} \|_{1- d_{r-1} - \delta_r} 
\le & \left(\frac{2e}{\rho\sigma}+\frac{e^2}{R^{2}}\right) \frac{1}{\delta_r^2 a_r} \left\| f_1^{({\rmI};r,r)}\right\|_{1- d_{r-1} -\delta_r}
\cr
 \le & M^{3r-2} \left(\frac{2e}{\rho\sigma}+\frac{e^2}{R^{2}}\right) \frac{4 \pi^4 \Ebarra r^{4+\tau}K^\tau}{\gamma }  \Vscr\left(\Hscr_1^{({\rm I}; \, r,r)}\right) \nu_{r,r}^{({\rm I})}\,,
\cr
}}
$$
we can carry out the justifications of the inequalities involving $\chi_1^{(r)}$ in a similar way to the method used for $\chi_0^{(r)}$.

\noindent
For showing the inequalities satisfied by $f_\ell^{({\rm II};r,s)} $, we have to distinguish between different cases. Since $f_\ell^{({\rm II};r,r+m)} = f_\ell^{({\rm I};r,r+m)}$ for $\ell=0,1$ and $0<m<r$, we can use the formula already proved for what concerns the term of index I. Observing that in this case $\Hscr_{\ell}^{({\rm II};r,r+m)} = \Hscr_{\ell}^{({\rm I};r,r+m)}$ and using a simple consequence of the second definition in~\eqref{frm:seqnu}, i.e., $\nu_{r,r+m}^{(\rm II)}= \nu_{r,r+m}^{(\rm I)} +\nu_{r,r}^{(\rm I)}\cdot \nu_{r,m}^{(\rm I)} > \nu_{r,r+m}^{(\rm I)}$, formula~\eqref{frm:f_l^(I,II,III;r,s)-lemmone} is proved, for what concerns the cases with $\ell = 0,1$ and $s=r+m$, where $0<m<r$.

\noindent
For $\ell=0$ and $s \ge 3r$ or $\ell=1$ and $s \ge 2r$ or $\ell\ge 2$ and $s\ge 0$, 
$$
\vcenter{\openup1\jot 
\halign{
$\displaystyle\hfil#$&$\displaystyle{}#\hfil$&$\displaystyle#\hfil$\cr
\|  f_\ell^{({\rm II};r,s)} \|_{1- d_{r-1} - 2\delta_r} &\le&
\sum_{j=0}^{\lfloor s/r \rfloor}\frac{1}{j!}  \|  \Lie{\chi_1^{(r)}}^j  f_{\ell+j}^{({\rm I}; r,s-jr)}  \|_{1- d_{r-1} - 2\delta_r}
\cr
&\le &\sum_{j=0}^{\lfloor s/r \rfloor} M^{(3r-1)j}\, \Vscr\left(\Gscr_1^{(r)}\right)^j \left(\nu_{r,r}^{({\rm I})} \right)^j
 \Ebarra \frac{M^{3(s-jr)- \zeta_{\ell+j}}}{2^{l+j}} \Vscr\left(\Hscr_{\ell+j}^{({\rm I}; r, s-jr)}\right)\nu_{r,s-jr}^{({\rm I})} \, ;
\cr
}}
$$
then, observing that $1/2^j <1$, $M^{j+\zeta_{l+j}} \ge M^{\zeta_l}$ and using the fourth definition of the list of indices in~\eqref{frm:H_l^(r,s)-lemmone} and that for $\nu_{r,r}^{(\rmII)}$ in~\eqref{frm:seqnu}, the inequality~\eqref{frm:f_l^(I,II,III;r,s)-lemmone} follows, in the considered cases.

\noindent
For $2r\le s <3r$,
$$
\vcenter{\openup1\jot 
\halign{
$\displaystyle\hfil#$&$\displaystyle{}#\hfil$&$\displaystyle#\hfil$\cr
\|  f_0^{({\rm II};r,s)} \|_{1- d_{r-1} - 2\delta_r} &\le&
\|   f_0^{({\rm I}; r,s)}  \|_{1- d_{r-1} - 2\delta_r} + \|  \Lie{\chi_1^{(r)}} f_1^{({\rm I}; r,s-r)}  \|_{1- d_{r-1} - 2\delta_r} 
\cr
& \le & \sum_{j=0}^1 M^{(3r-1)j} \,\Vscr\left(\Gscr_1^{(r)}\right)^j (\nu_{r,r}^{({\rm I})})^j \bar E\, \frac{M^{3(s-jr)-\zeta_j}}{2^j}\, \Vscr\left(\Hscr_j^{({\rm I}; r, s-jr)}\right) \nu_{r,s-jr}^{({\rm I})}\,;
\cr
}}
$$
again, using $\nu_{r,s}^{(\rm II)} \ge \nu_{r,s}^{(\rm I)}+\nu_{r,r}^{(\rm I)}\, \nu_{r,s-r}^{(\rm I)}$ and the fourth definition of the lists of indices in~\eqref{frm:H_l^(r,s)-lemmone}, the inequality~\eqref{frm:f_l^(I,II,III;r,s)-lemmone} is fully justified for what concerns the upper index {\rm II}.

\noindent
The upper bound for $\chi_2^{(r)}$ follows as for the other generating functions, using these last formul{\ae}.

\noindent
We have now to discuss the different cases for $f_\ell^{(r,s)}$.

\noindent
For $\ell \ge 3 $ and $s \ge 0$, 
$$
\vcenter{\openup1\jot 
\halign{
$\displaystyle\hfil#$&$\displaystyle{}#\hfil$&$\displaystyle#\hfil$\cr
\|  f_\ell^{(r,s)} \|_{1- d_{r-1} - 3\delta_r} &\le&
\sum_{j=0}^{\lfloor s/r \rfloor} \frac{1}{j!}  \| \Lie{\chi_2^{(r)}}^j  f_\ell^{({\rm II}; r,s-jr)}  \|_{1- d_{r-1} - 3\delta_r}
\cr
&\le &\sum_{j=0}^{\lfloor s/r \rfloor} M^{3rj} \, \Vscr\left(\Gscr_2^{(r)}\right)^j \left(\nu_{r,r}^{({\rm II})} \right)^j
 \Ebarra \, \frac{M^{3(s-jr)- \zeta_\ell}}{2^\ell}\Vscr\left(\Hscr_l^{({\rm II}; r, s-jr)}\right)\nu_{r,s-jr}^{({\rm II})}\,,
\cr
}}
$$ then the inequality~\eqref{frm:f_l^(r,s)-lemmone} follows in view
of the same arguments we previously used.  For $\ell=0,1$ and $s>r$,
the formula is the same of the previous case, with the only difference
that the sum is up to $\lfloor s/r \rfloor -1$ and this is considered
in the definition of $\Hscr_\ell^{(r,kr+m)}$, which appears in
formula~\eqref{frm:H_l^(r,s)-lemmone}.  We can estimate $f_2^{(r,s)}$
for $s \ge r$ with $\sum_{j=0}^{\lfloor s/r \rfloor -1} \frac
{1}{j!}\Lie{\chi_2^{(r)}}^j f_2^{(\rmII;r,s-jr)}$ and complete the
proof of inequality~\eqref{frm:f_l^(r,s)-lemmone} in a way that is
similar to what has been done before.

\noindent
Finally, we have to show the inequalities in~\eqref{frm:stima-variazione-frequenze}. They can be proved recalling the change of frequencies described in~\eqref{frm:chgfreq.toro} and~\eqref{frm:chgfreq.trasv} and using the inequalities already proved for  $f_2^{({\rmII}; r,r)}$ (observing that $ \| Z^{(r)}\| \le \| f_2^{({\rmII}; r,r)}\|$) jointly with Cauchy estimates.
Using the inequalities $\delta_r<1$, $\frac{1}{r^\tau} \le 1$ and $\frac{4 \pi^4 \bar E}{\sigma\rho} \le \gamma M$ (see definitions~\eqref{frm:dr-and-deltar} and~\eqref{frm:def-M}), the following estimates hold true
\begin{align*}
\frac{1}{\sigma}|\omega_j^{(r)} - \omega_j^{(r-1)}| 
\le& \frac{\epsilon^r}{\sigma} \Big|\frac{\partial \langle f_2^{({\rm II}; r,r)}\rangle_q}{\partial p_j}\Big|
 \le \frac{\epsilon^r}{\sigma \delta_r \rho} \| f_2^{({\rm II}; r,r)}\| _{1-d_{r-1}-2\delta_r}
\\
\le & \frac{\epsilon^r}{\sigma \delta_r \rho} \frac{\bar E}{4}M^{3r-1}\Vscr(\Hscr_2^{(\rmII;r,r)}) \nu_{r,r}^{(\rmII)}
\le   \frac{\epsilon^r}{\sigma \rho} \bar E \delta_r \left(\frac{4 \pi^4 r^4}{r^\tau}\right) r^\tau M^{3r-1}\Vscr(\Hscr_2^{(\rmII;r,r)}) \nu_{r,r}
\\
\le &  \frac{\epsilon^r}{\sigma \rho} \bar E \delta_r (4 \pi^4) M^{3r-1}\Vscr(\Gscr_2^{(r)}) \nu_{r,r} \le \gamma \epsilon^r M^{3r}\Vscr(\Gscr_2^{(r)}) \nu_{r,r} \,,
\end{align*}
where the last definitions in both~\eqref{frm:seqnu} and~\eqref{frm:H_l^(r,s)-lemmone} have been taken into account.

\noindent
By proceeding in a similar way and using that $\frac{4 \pi^4 \bar E}{R^2} \le \gamma M$, we can prove the other inequality making part of formula~\eqref{frm:stima-variazione-frequenze}:
\begin{align*}
|\Omega_j^{(r)} - \Omega_j^{(r-1)}| 
\le& \epsilon^r \Big|\frac{\partial ^2 Z^{(r)}}{\partial z_j \partial \imunit z_j}\Big|
 \le \frac{2\epsilon^r}{\delta_r^2 R^2} \| f_2^{(\rmII; r,r)}\|_{1-d_{r-1}-2\delta_r}
\le  \frac{2\epsilon^r}{\delta_r^2 R^2} \frac{\bar E}{4}M^{3r-1}\Vscr(\Hscr_2^{(\rmII;r,r)}) \nu_{r,r}^{(\rmII)}
\\ 
\le &  \frac{\epsilon^r}{R^2} \bar E \left(\frac{4 \pi^4 r^4}{r^\tau}\right) r^\tau M^{3r-1}\Vscr(\Hscr_2^{(\rmII;r,r)}) \nu_{r,r}
\\
\le &  \frac{\epsilon^r}{R^2} \bar E (4 \pi^4) M^{3r-1}\Vscr(\Gscr_2^{(r)}) \nu_{r,r} \le \gamma \epsilon^r M^{3r}\Vscr(\Gscr_2^{(r)}) \nu_{r,r} \,.
\end{align*}
This concludes the proof.

\section{Proof of lemma~\ref{lem:small-div-acc}}
\label{app:accumulodiv}
The proof proceeds by induction on $r$.
For $r=0$, $\Hscr_\ell^{(0,s)} = \emptyset$, thus~\eqref{frm:Nk_Ham} are obviously satisfied.

\noindent
Let us suppose that~\eqref{frm:Nk_Ham} are true for $r-1$ and prove they are still valid for $r$. 

\noindent
We start showing the inequalities in~\eqref{frm:Nk_gen} for $j=0$.  By
the definition of $\Gscr_0^{(r)}$,
$$
\Nscr_k\big(\Gscr_0^{(r)}\big) =  \Nscr_k\big(\Hscr_0^{(r-1,r)}\big) + \Nscr_k(\{ r \})\,.
$$ Therefore, we have to distinguish between different cases:\\ (i) if
$k< \lfloor \log_2 r \rfloor$, the second addend is null and we can
use the inductive hypotesis to conclude that $
\Nscr_k\big(\Hscr_0^{(r-1,r)}\big) \le 3 \left\lfloor \frac{r}{2^k}
\right\rfloor$;\\ (ii) if $k =\lfloor\log_2 r \rfloor$, only the
second term gives a contribution and, then,
$\Nscr_k\big(\Gscr_0^{(r)}\big) =1 $;\\ (iii) the case with $k >
\lfloor \log_2 r \rfloor$ is trivial (thereafter it will be
omitted).\\

In order to prove~\eqref{frm:Nk_Ham}, it is useful to show that the same inequalities are satisfied by the intermediate list of indices defined in~\eqref{frm:H_l^(r,s)-lemmone}, i.e.,
\begin{equation}
\Nscr_k\left(\Hscr_\ell^{({\rm I};r,s)}\right) , \,
\Nscr_k\left(\Hscr_\ell^{({\rm II};r,s)}\right) \le \left\{
\vcenter{\openup1\jot \halign{
    $\displaystyle#\hfil$&\quad$\displaystyle{}#\hfil$&\quad$\displaystyle#\hfil$\cr
    3 \left\lfloor \frac {s}{2^k} \right\rfloor & {\rm for} \ k <
    \lfloor \log_2{r} \rfloor \>, \cr 3 \left\lfloor \frac {s}{2^k}
    \right\rfloor -3 +\ell & {\rm for} \ \ell \le 3, \, k = \lfloor
    \log_2{r} \rfloor \>, \cr 3 \left\lfloor \frac {s}{2^k}
    \right\rfloor & {\rm for} \ \ell > 3, \, k = \lfloor \log_2{r}
    \rfloor \>, \cr 0 & {\rm for} \ k > \lfloor \log_2{r} \rfloor \>.
    \cr }} \right .
\label{frm:Nk_Ham-intermedie}
\end{equation}
\noindent
Therefore, we start to prove the estimate for $\Hscr_\ell^{({\rm I}; \, r,s)}$, showing that every argument of the operator $\Mscr \Ascr \Xscr$ cannot exceed such an upper bound.
For $\ell= 0$, $s \ge 2r$ and $l \ge 1$, $s \ge 0$, by the definition of $\Hscr_\ell^{(\rmI; r,s)}$ in~\eqref{frm:H_l^(r,s)-lemmone}, we need an estimate for
$$
\Nscr_k\left(\Hscr_\ell^{({\rm I}; \, r,s)}\right) \le j \Nscr_k\left(\Gscr_0^{(r)}\right) + \Nscr_k\left(\Hscr_{\ell+2j}^{(r-1,s-jr)}\right)\,.
$$ Therefore,\\ (i) if $k< \lfloor \log_2 r \rfloor$, $\forall \ 0\le
j \le \lfloor \frac{s}{r} \rfloor$, the estimate
in~\eqref{frm:Nk_Ham-intermedie} is satisfied since $j
\Nscr_k\left(\Gscr_0^{(r)}\right) +
\Nscr_k\left(\Hscr_{\ell+2j}^{(r-1,s-jr)}\right) \le 3j\left\lfloor
\frac{r}{2^k} \right\rfloor + 3\left\lfloor \frac{s-jr}{2^k}
\right\rfloor - \zeta_{\ell+2j} \le 3\left\lfloor \frac{s}{2^k}
\right\rfloor - \zeta_{\ell+2j} \le 3\left\lfloor \frac{s}{2^k}
\right\rfloor \, ; $\\ (ii) if $r=2^k $, the second terms is null;
then, the estimate is proved, provided that $j \le 3\left\lfloor
\frac{s}{r} \right\rfloor - \zeta_\ell $, i.e., $2\left\lfloor
\frac{s}{r} \right\rfloor \ge \zeta_\ell$, that holds true for
$\ell\ge 1, \ s \ge r$ or for $\ell =0, \ s\ge 2r$;\\ (iii) if $k =
\lfloor \log_2 r \rfloor$ and $2^k< r$, we need an estimate similar to
that in the case~(i) above; $\forall \ 0 \le j \le \lfloor
s/r\rfloor$, we have $ \Nscr_k\left(\Gscr_0^{(r)}\right) +
\Nscr_k\left(\Hscr_{\ell+2j}^{(r-1,s-jr)}\right) \le j+ 3\lfloor
\frac{s-jr}{2^k}\rfloor - \zeta_{l+2j}= j\lfloor \frac{r}{2^k}\rfloor+
3\lfloor \frac{s-jr}{2^k}\rfloor - \zeta_{l+2j}\le 3 \lfloor
\frac{s}{2^k}\rfloor -2j - \zeta_{l+2j} \le 3 \lfloor
\frac{s}{2^k}\rfloor - \zeta_l.$\\

\noindent
For $r <s<2r$, $\Hscr_0^{({\rm I}; r,s)} = \Hscr_0^{(r-1,s)}$. Since $\lfloor \frac{s}{r} \rfloor = 1$, using the inductive hypothesis the inequality is satisfied and the estimate for $\Nscr_k\left(\Hscr_\ell^{({\rm I}; \, r,s)}\right)$ is proved.

Let us prove now the inequality for $\Nscr_k\big(\Gscr_1^{(r)}\big)$. By definition, 
$$
\Nscr_k\big(\Gscr_1^{(r)}\big) =  \Nscr_k\big(\Hscr_1^{({\rm I}; r,r)}\big) + \Nscr_k(\{ r \})\,;
$$ therefore, using the inequalities proved for
$\Nscr_k\big(\Hscr_1^{({\rm I}; r,r)}\big)$, the estimate for
$\Nscr_k\big(\Gscr_1^{(r)}\big) $ is valid. In fact,\\ (i) if $k<
\lfloor \log_2 r \rfloor$, the estimate is an immediate consequence
of~\eqref{frm:Nk_Ham-intermedie};\\ (ii) if $k =\lfloor\log_2 r
\rfloor $, \ $\Nscr_k\big(\Gscr_1^{(r)}\big) = 3\lfloor \frac{r}{2^k}
\rfloor -2 +1 = 2$\,.\\

\noindent
We procede similarly for the other cases.

\noindent
For $\ell=0$, $s\ge 3r$ or $\ell=1$, $s \ge 2r$ or $l\ge 2$, $s \ge r$ and $\forall \ 0 \le j \le \lfloor \frac{s}{r}  \rfloor$, $\Hscr_\ell^{(\rmII; r,s)}$ satisfies the following equality:
$$
\Nscr_k\left(\Hscr_\ell^{({\rm II}; r, s)}\right) = j  \Nscr_k\left(\Gscr_1^{(r)}\right) +  \Nscr_k\left(\Hscr_{\ell+j}^{({\rm I}; r, s-jr)}\right) \,.
$$ As a consequence,\\ (i) if $k< \left\lfloor \log_2 r \right\rfloor
$, $ \Nscr_k\big(\Hscr_\ell^{({\rm II}; r, s)}\big) \le 3j
\left\lfloor \frac {r}{2^k} \right\rfloor + 3 \left\lfloor \frac
             {s-jr}{2^k}\right\rfloor \le 3 \left\lfloor \frac
             {s}{2^k}\right\rfloor $;\\ (ii) if $k= \left\lfloor
             \log_2 r \right\rfloor $, $\Nscr_k\big(\Hscr_\ell^{({\rm
                 II}; r, s)}\big) \le 2j + 3 \left\lfloor \frac
                 {s-jr}{2^k}\right\rfloor - \zeta_{\ell+j} = 3
                 \left\lfloor \frac {s}{2^k}\right\rfloor -j -
                 \zeta_{\ell+j} \le 3 \left\lfloor \frac
                      {s}{2^k}\right\rfloor - \zeta_\ell \, .  $\\ For
                      $\ell=0$, $ r< s <3r$ or $\ell=1$, $r<s < 2r$,
                      the only difference is due to the fact that $j
                      \le \lfloor \frac{s}{r} \rfloor -1$, but the
                      proof is the same as in the previous case.

\noindent
The list of indices of the generating function $\chi_2^{(r)}$
satisfies the following equality:
$$
\Nscr_k\left(\Gscr_2^{(r)}\right) =  \Nscr_k\left(\Hscr_2^{({\rm II}; r,r)}\right) + \Nscr_k(\{ r \}) \,.
$$
Then,\\
(i) if $k< \lfloor \log_2 r \rfloor$, the estimate is an immediate consequence of~\eqref{frm:Nk_Ham-intermedie};\\
(ii) if $k =\lfloor\log_2 r \rfloor $, $\Nscr_k\big(\Gscr_2^{(r)}\big)  = 3\left\lfloor \frac{r}{2^k} \right\rfloor -1 +1 = 3\,$ .\\

\noindent
We come to prove the formula in~\eqref{frm:Nk_Ham}.

\noindent
For $\ell \ge 3$, $s \ge 0$ and $j \le \lfloor \frac{s}{r}  \rfloor$ , $\Hscr_\ell^{(r, s)}$ satisfies the following equality:
$$ \Nscr_k\left(\Hscr_\ell^{(r, s)}\right) = j
\Nscr_k\left(\Gscr_2^{(r)}\right) + \Nscr_k\left(\Hscr_\ell^{({\rm
    II}; r, s-jr)}\right) \, .
$$
Then,\\
(i) if $k< \left\lfloor \log_2 r \right\rfloor $ , $\Nscr_k\big(\Hscr_\ell^{(r, s)}\big) \le 3j \left\lfloor \frac {r}{2^k} \right\rfloor + 3 \left\lfloor \frac {s-jr}{2^k}\right\rfloor \le 3 \left\lfloor \frac {s}{2^k}\right\rfloor
$;\\
(ii) if $k= \left\lfloor \log_2 r \right\rfloor $, $
\Nscr_k\big(\Hscr_\ell^{(r, s)}\big)  \le 3j + 3 \left\lfloor \frac {s-jr}{r}\right\rfloor = 3 \left\lfloor \frac {s}{r}\right\rfloor
$. \\

\noindent
For $\ell=0,1,2$, $s\ge1$, $0<m<r$ and $j \le s-1$ , 
$$
\Nscr_k\left(\Hscr_\ell^{(r, sr+m)}\right) = j  \Nscr_k\left(\Gscr_2^{(r)}\right) +  \Nscr_k\left(\Hscr_\ell^{({\rm II}; r, (s-j)r+m)}\right).
$$ The two cases are\\ (i) if $k< \left\lfloor \log_2 r \right\rfloor
$, $ \Nscr_k\big(\Hscr_\ell^{(r, sr+m)}\big) \le 3j \left\lfloor \frac
     {r}{2^k} \right\rfloor + 3 \left\lfloor \frac
     {(s-j)r+m}{2^k}\right\rfloor \le 3 \left\lfloor \frac
     {sr+m}{2^k}\right\rfloor $;\\ (ii) if $k= \left\lfloor \log_2 r
     \right\rfloor $, $ \Nscr_k\big(\Hscr_\ell^{(r, sr+m)}\big) \le 3j
     + 3 \left\lfloor \frac {(s-j)r+m}{r}\right\rfloor - \zeta_\ell =
     3 \left\lfloor \frac {sr+m}{r}\right\rfloor - \zeta_\ell \,.  $\\

\noindent 
For $\ell=0,1$, $s\ge2$ and $j \le s-2$, $\Hscr_\ell^{(r, sr)}$
satisfies the following equality:
$$ \Nscr_k\left(\Hscr_\ell^{(r, sr)}\right) = j
\Nscr_k\left(\Gscr_2^{(r)}\right) + \Nscr_k\left(\Hscr_\ell^{({\rm
    II}; r, (s-j)r)}\right),
$$ and, repeating the previous computations, the inequality
follows. Moreover, the estimate for $\Hscr_2^{(r, sr)}$ follows from
analogous computations.

\end{document}